\documentclass[a4paper,12pt]{article} 
\usepackage{amssymb}
\usepackage{amscd}
\usepackage{amsmath}
\usepackage{graphicx}
\usepackage{latexsym} 
\usepackage{enumerate}

\usepackage{theorem}
\theoremstyle{plain}
\theorembodyfont{\itshape}
\newtheorem{theorem}{Theorem}[section]

\newtheorem{lemma}[theorem]{Lemma}

\theorembodyfont{\rmfamily}

\newtheorem{definition}[theorem]{Definition}
\newtheorem{remark}[theorem]{Remark}

\newtheorem{example}[theorem]{Example}

\def\Aff{\mathrm{Aff}}
\def\Aut{\mathrm{Aut}}

\def\det{\mathrm{det}}
\def\End{\mathrm{End}}

\def\id{\mathrm{id}}

\def\PVI{\mathrm{P}_{\mathrm{VI}}}
\def\Per{\mathrm{Per}} 
\def\PS{\mathrm{PS}}

\def\Res{\mathrm{Res}}
\def\RH{\mathrm{RH}} 
\def\rh{\mathrm{rh}}

\def\Tr{\mathrm{Tr}}

\def\E{\mathcal{E}}

\def\K{\mathcal{K}}

\def\M{\mathcal{M}}
\def\O{\mathcal{O}}
\def\R{\mathcal{R}}
\def\Y{\mathcal{Y}}
\def\Sol{\mathcal{S}}

\def\C{\mathbb{C}}

\def\N{\mathbb{N}}
\def\P{\mathbb{P}}

\def\bR{\mathbb{R}}
\def\Z{\mathbb{Z}}
\def\Wall{\mathbf{Wall}} 
\def\a{\alpha}
\def\b{\beta}
\def\ga{\gamma}
\def\k{\kappa}
\def\l{\lambda}
\def\si{\sigma}
\def\th{\theta}
\def\vD{\varDelta}
\def\vG{\varGamma}

\def\Th{\Theta}
\def\ve{\varepsilon}
\def\carl{\circlearrowleft}
\def\car{\curvearrowright}

\def\ra{\rightarrow}

\def\ol{\overline}

\def\ci{\circ}
\def\dfrac#1#2{{\displaystyle\frac{#1}{#2}}}

\def\la{\langle}
\def\ra{\rangle}

\def\wt{\widetilde} 
\def\-{\phantom{-}}
\title{\bf An Ergodic Study of Painlev\'e VI} 
\author{Katsunori Iwasaki and Takato Uehara \\ \\
Graduate School of Mathematics, Kyushu University \\
6-10-1 Hakozaki, Higashi-ku, Fukuoka 812-8581 
Japan\thanks{E-mail addresses: {\tt iwasaki@math.kyushu-u.ac.jp} \ 
and \ {\tt ma205003@math.kyushu-u.ac.jp}}} 
\date{Dedicated to Professor Masuo Hukuhara on his $100$th birthday} 
\begin{document}
\maketitle
\begin{abstract} 
An ergodic study of Painlev\'e VI is developed. 
The chaotic nature of its Poincar\'e return map is 
established for almost all loops. 
The exponential growth of the numbers of periodic 
solutions is also shown. 
Principal ingredients of the arguments are a moduli-theoretical 
formulation of Painlev\'e VI, a Riemann-Hilbert correspondence, 
the dynamical system of a birational map on a cubic surface, and 
the Lefschetz fixed point formula. 
\end{abstract} 
\section{Introduction}  \label{sec:intro} 
Painlev\'e equations have been investigated actively in 
recent years. 
However most researches have been done from the 
viewpoint of integrable systems and little attention has 
been paid to the ergodic and chaotic aspects of their dynamics. 
In this paper we develop an ergodic study of the sixth 
Painlev\'e equation $\PVI(\k)$ and explore the chaotic behavior 
of its global dynamics, namely, that of its Poincar\'e return 
map. 
The aim of this paper is to show that the Poincar\'e 
return map is chaotic along almost all loops in the space 
of independent variable 
\[
Z = \P^1-\{0,1,\infty\}. 
\]
The exponential growth of the number of periodic 
solutions along those loops is also established. 
\par
The sixth Painlev\'e equation $\PVI(\k)$ is a Hamiltonian 
system of differential equations 
\begin{equation} \label{eqn:PVI}
\dfrac{d q}{d z} = \dfrac{\partial H(\k)}{\partial p}, 
\qquad 
\dfrac{d p}{d z} = -\dfrac{\partial H(\k)}{\partial q}, 
\end{equation}
with an independent variable $z \in Z$ and unknown functions 
$q = q(z)$ and $p = p(z)$, depending on complex parameters 
$\k = (\k_0,\k_1,\k_2,\k_3,\k_4)$ in a $4$-dimensional affine 
space 
\begin{equation} \label{eqn:K}
\K := \{\, \k = (\k_0,\k_1,\k_2,\k_3,\k_4) \in \C^5 \,:\, 
2 \k_0 + \k_1 + \k_2 + \k_3 +\k_4 = 1 \,\},  
\end{equation}
where the Hamiltonian $H(\k) = H(q,p,z;\k)$ is given by 
\[
\begin{array}{rcl}
z(z-1) H(\k) &=&  
(q_0q_1q_z) p^2 - \{\k_1q_1q_z 
+ (\k_2-1)q_0q_1 + \k_3q_0q_z \} p 
+ \k_0(\k_0+\k_4) q_z, 
\end{array}
\]
with $q_{\nu} = q - \nu$ for $\nu \in \{0, 1, z \}$. 
More intrinsically, $\PVI(\k)$ can be formulated as a 
holomorphic uniform foliation on a fibration of certain 
smooth quasi-projective rational surfaces 
\begin{equation} \label{eqn:phase}
\pi_{\k} : M(\k) \to Z, 
\end{equation} 
which is transversal to each fiber of the fibration. 
Equation (\ref{eqn:PVI}) is just a coordinate expression of 
the foliation in terms of a natural coordinate system on an 
affine open subset of the phase space $M(\k)$. 
See \cite{AL,IIS1,IIS2,Okamoto,STT,Sakai} for various 
construction of the space $M(\k)$. 
Especially the papers \cite{IIS1,IIS2} give a comprehensive 
description of it as a moduli space of stable parabolic 
connections. 
The fiber $M_z(\k)$ over $z \in Z$ is called the space of 
initial conditions at time $z$. 
\par 
Since the foliation is uniform (Painlev\'e property), 
each loop $\ga \in \pi_1(Z,z)$ admits global horizontal lifts 
along the foliation and induces an automorphism 
$\ga_* : M_z(\k) \to M_z(\k)$, called the holonomy or the 
Poincar\'e return map along the loop $\ga$. 
Then the global structure of the foliation is described by 
the holonomy representation 
\begin{equation} \label{eqn:PSz} 
\PS_z(\k) \,\, : \,\, \pi_1(Z,z) \to \Aut\, M_z(\k), 
\qquad \ga \mapsto \ga_*. 
\end{equation}
which is referred to as the Poincar\'e section of the 
Painlev\'e dynamical system $\PVI(\k)$. 
Here and hereafter a {\sl loop} means the homotopy class 
of a loop without further comment. 
\par 
In this paper we are interested in the dynamics of the 
Poincar\'e return map $\ga_* : M_z(\k) \carl$ for each 
individual loop $\ga \in \pi_1(Z,z)$. 
One of our main results will state that $\ga_*$ always exhibits 
a chaotic behavior as long as $\ga$ is a non-elementary loop 
(see Theorem \ref{thm:chaos}),  
where the adjective ``chaotic" and the words ``non-elementary 
loop" are used in the following senses. 
\begin{definition} \label{def:chaos} 
The dynamical system of a holomorphic map $f : S \to S$ on a 
complex surface $S$ (in our case, $S = M_z(\k)$ and 
$f = \ga_*$) is said to be {\sl chaotic} if there exists an 
$f$-invariant Borel probability measure $\mu$ on $S$ such 
that the following conditions are satisfied: 
\begin{enumerate}
\item[(C1)] $f$ has a positive entropy $h_{\mu}(f) > 0$ with 
respect to the measure $\mu$.  
\item[(C2)] $f$ is mixing with respect to the measure $\mu$, 
that is, $\mu(f^{-n}(A) \cap B) \to \mu(A) \mu(B)$ as $n \to 
\infty$ for any Borel subsets $A$ and $B$ of $S$. 
In particular, $f$ is ergodic with respect to $\mu$. 
\item[(C3)] $\mu$ is a hyperbolic measure of saddle type, 
that is, the two Lyapunov exponents $L_{\pm}(f)$ of $f$ 
with respect to the ergodic measure $\mu$ satisfy 
$L_-(f) < 0 < L_+(f)$. 
Moreover, $\mu$ has product structure with respect to local 
stable and unstable manifolds. 
\item[(C4)] hyperbolic periodic points of $f$ are dense in the 
support of $\mu$. 
\end{enumerate}
\end{definition} 
For the basic terminology used here, we refer to the standard 
textbooks \cite{KH,Walters} on dynamical systems and ergodic 
theory. 
While there are many possible definitions of a chaotic dynamical 
system \cite{Devaney}, the definition adopted here is a typical 
one possessing the three ingredients usually required for 
a ``chaos": (i) unpredictability, that is, the sensitive 
dependence on initial values represented by conditions 
(C1) and (C3); 
(ii) indecomposability, that is, ergodicity or a related 
property as in (C2); (iii) an element of regularity, that is, 
the existence of periodic points which are dense in a dynamically 
interesting subset as in (C4). 
We also remark that conditions (C2) and (C3) imply that the 
dynamical system $f$ with invariant measure $\mu$ is Bernoulli, 
namely, it is measurably conjugate to a Bernoulli shift \cite{OW}. 
\par
\begin{figure}[t] 
\begin{center}
\unitlength 0.1in
\begin{picture}(29.80,29.80)(3.60,-32.10)
%
\special{pn 13}%
\special{ar 1850 1720 1490 1490  0.0000000 6.2831853}%
%
\special{pn 20}%
\special{ar 1844 1012 324 324  0.0000000 6.2831853}%
%
\special{pn 20}%
\special{sh 0.600}%
\special{ar 1850 1726 47 47  0.0000000 6.2831853}%
%
\special{pn 20}%
\special{pa 1844 1336}%
\special{pa 1844 1684}%
\special{fp}%
%
\special{pn 20}%
\special{ar 1292 2197 324 324  0.0000000 6.2831853}%
%
\special{pn 20}%
\special{pa 1539 1987}%
\special{pa 1804 1762}%
\special{fp}%
%
\special{pn 20}%
\special{ar 2414 2192 324 324  0.0000000 6.2831853}%
%
\special{pn 20}%
\special{pa 2165 1985}%
\special{pa 1898 1762}%
\special{fp}%
%
\special{pn 20}%
\special{sh 0.600}%
\special{ar 1838 1006 47 47  0.0000000 6.2831853}%
%
\special{pn 20}%
\special{sh 0.600}%
\special{ar 1292 2200 47 47  0.0000000 6.2831853}%
%
\special{pn 20}%
\special{sh 0.600}%
\special{ar 2408 2188 47 47  0.0000000 6.2831853}%
\put(18.0000,-12.1000){\makebox(0,0)[lb]{$0$}}%
\put(12.5000,-24.1000){\makebox(0,0)[lb]{$1$}}%
\put(23.4000,-23.8000){\makebox(0,0)[lb]{$\infty$}}%
\put(17.9000,-19.2400){\makebox(0,0)[lb]{$z$}}%
\put(17.6400,-29.4800){\makebox(0,0)[lb]{$Z$}}%
%
\special{pn 13}%
\special{ar 1842 1010 381 381  4.3704392 5.1338082}%
%
\special{pn 13}%
\special{pa 1728 644}%
\special{pa 1686 668}%
\special{fp}%
\special{sh 1}%
\special{pa 1686 668}%
\special{pa 1754 652}%
\special{pa 1732 642}%
\special{pa 1734 618}%
\special{pa 1686 668}%
\special{fp}%
%
\special{pn 13}%
\special{ar 1286 2207 381 381  2.0165887 2.7791331}%
%
\special{pn 13}%
\special{pa 1108 2546}%
\special{pa 1154 2559}%
\special{fp}%
\special{sh 1}%
\special{pa 1154 2559}%
\special{pa 1095 2522}%
\special{pa 1103 2544}%
\special{pa 1084 2560}%
\special{pa 1154 2559}%
\special{fp}%
%
\special{pn 13}%
\special{ar 2423 2187 381 381  0.2772048 1.0389390}%
%
\special{pn 13}%
\special{pa 2788 2305}%
\special{pa 2792 2258}%
\special{fp}%
\special{sh 1}%
\special{pa 2792 2258}%
\special{pa 2766 2323}%
\special{pa 2787 2311}%
\special{pa 2806 2326}%
\special{pa 2792 2258}%
\special{fp}%
\put(11.7000,-26.9000){\makebox(0,0)[lb]{$\ga_2$}}%
\put(23.0400,-26.9000){\makebox(0,0)[lb]{$\ga_3$}}%
\put(22.3200,-10.5800){\makebox(0,0)[lb]{$\ga_1$}}%
\end{picture}%
\end{center}
\caption{Three basic loops in $Z = \P^1-\{0,1,\infty\}$} 
\label{fig:elementary} 
\end{figure}
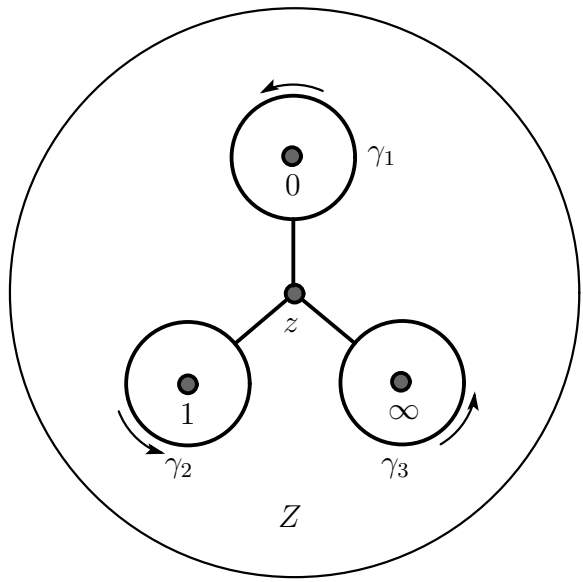

 
Next we explain what we mean by the words ``non-elementary loop". 
Treating the three fixed singular points $0$, $1$, $\infty$ 
of $\PVI(\k)$ symmetrically, we put 
\[
z_1 = 0, \qquad z_2 = 1, \qquad z_3 = \infty. 
\]
For each $i \in \{1,2,3\}$ with $\{i,j,k\} = \{1,2,3\}$, 
let $\ga_i \in \pi_1(Z,z)$ be a loop surrounding the points 
$z_i$ once anti-clockwise, leaving the 
remaining points $z_j$ and $z_k$ outside, as in 
Figure \ref{fig:elementary}. 
Then the fundamental group $\pi_1(Z,z)$ is generated 
by $\ga_1$, $\ga_2$, $\ga_3$, having a defining relation 
\begin{equation} \label{eqn:relation}
\ga_1 \ga_2 \ga_3 = 1. 
\end{equation} 
\begin{definition} \label{def:elementary} 
A loop $\ga \in \pi_1(Z,z)$ is said to be {\sl elementary} 
if $\ga$ is conjugate to the loop $\ga_i^m$ for some 
$i \in \{1,2,3\}$ and $m \in \Z$. 
Otherwise, $\ga$ is said to be {\sl non-elementary}. 
\end{definition} 
\par
The second issue to be discussed in this paper is the number 
of periodic solutions to $\PVI(\k)$. 
Given a loop $\ga \in \pi_1(Z,z)$ and a positive integer 
$N \in \N$, we are interested in the number of periodic 
solutions to $\PVI(\k)$ of period $N$ along the loop $\ga$. 
To be more precise, we wish to count the number of all initial 
conditions at time $z$ that come back to the original positions 
after the $N$-th iterate of the Poincar\'e return map along 
$\ga$, namely, the cardinality of the set 
\begin{equation} \label{eqn:per}
\Per_N(\ga;\k) = \{\, Q \in M_z(\k)\,:\, \ga_*^N Q = Q \,\}. 
\end{equation} 
It will be shown that for any non-elementary loop $\ga$, the 
cardinality is finite for every period $N \in \N$ and grows 
exponentially as the period $N$ tends to infinity (see 
Theorem \ref{thm:per}). 
We shall also give an algorithm to count the number exactly 
as well as to determine its exponential growth rate 
explicitly (see Theorem \ref{thm:algorithm}). 
The logarithm of this rate will give the entropy of the 
Poincar\'e return map $\ga_*$. 
Recently several authors 
\cite{Boalch1,Boalch2,DM,Hitchin1,Hitchin2,Kitaev,Mazzocco} 
have been interested in finding algebraic solutions, which 
must have only finitely many branches under the 
analytic continuations along {\sl all} loops in $Z$. 
On the other hand, in this article we will be concerned with 
those solutions which are finitely many-valued along a fixed 
{\sl single} loop. 
\par
Painlev\'e equations and dynamical systems on complex manifolds 
are two subjects of mathematics which have attracted much attention 
in recent years. 
In this paper we shall demonstrate a substantial interplay between 
them by presenting a fruitful application to the former subject 
of the latter. 
On the former side, algebraic geometry of Painlev\'e equations, 
especially a moduli-theoretical formulation of 
Painlev\'e dynamical systems \cite{IIS1,IIS2} is an 
essential ingredient of our discussion. 
On the latter side, recent advances in complex surface dynamics, 
especially some deep ergodic studies of birational maps of complex 
surfaces \cite{BD,DF,DS,Dujardin} are another basis of our analysis. 
These two stuffs are combined fruitfully via a Riemann-Hilbert 
correspondence to reveal the chaotic nature of the sixth 
Painlev\'e dynamics. 
\section{Main Results} \label{sec:result}
\begin{figure}[t]
\begin{center}
\unitlength 0.1in
\begin{picture}(14.83,13.55)(8.95,-14.45)
%
\special{pn 20}%
\special{ar 1795 797 52 52  6.1412883 6.2831853}%
\special{ar 1795 797 52 52  0.0000000 6.1180366}%
%
\special{pn 20}%
\special{ar 1255 257 52 52  6.1412883 6.2831853}%
\special{ar 1255 257 52 52  0.0000000 6.1180366}%
%
\special{pn 20}%
\special{ar 2317 284 52 52  6.1412883 6.2831853}%
\special{ar 2317 284 52 52  0.0000000 6.1180366}%
%
\special{pn 20}%
\special{ar 1246 1337 52 52  6.1412883 6.2831853}%
\special{ar 1246 1337 52 52  0.0000000 6.1180366}%
%
\special{pn 20}%
\special{ar 2326 1346 52 52  6.1412883 6.2831853}%
\special{ar 2326 1346 52 52  0.0000000 6.1180366}%
%
\special{pn 20}%
\special{pa 1300 311}%
\special{pa 1750 761}%
\special{fp}%
%
\special{pn 20}%
\special{pa 1840 842}%
\special{pa 2299 1301}%
\special{fp}%
%
\special{pn 20}%
\special{pa 2281 302}%
\special{pa 1831 752}%
\special{fp}%
%
\special{pn 20}%
\special{pa 1759 833}%
\special{pa 1291 1301}%
\special{fp}%
\put(18.0400,-10.2600){\makebox(0,0){$w_0$}}%
\put(26.7700,-1.1700){\makebox(0,0)[rt]{$w_2$}}%
\put(24.9700,-15.3000){\makebox(0,0){$w_4$}}%
\put(11.2000,-15.0300){\makebox(0,0){$w_3$}}%
\put(9.7600,-0.9000){\makebox(0,0)[lt]{$w_1$}}%
\end{picture}%
\end{center}
\caption{Dynkin diagram of type $D_4^{(1)}$}
\label{fig:dynkin} 
\end{figure}
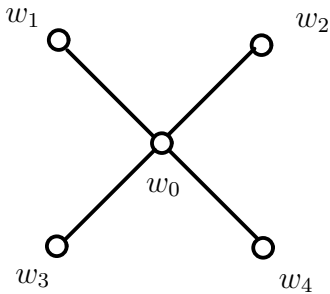

Let us describe our main results in more detail. 
In this paper we make a certain generic assumption on the 
parameters $\k \in \K$ to avoid a technical 
difficulty (see Remark \ref{rem:generic}). 
To this end we recall an affine Weyl group structure 
of the parameter space $\K$ \cite{IIS1,Iwasaki3}. 
In view of formula (\ref{eqn:K}), the affine space $\K$ can 
be identified with the linear space $\C^4$ by the isomorphism 
\[
\K \to \C^4, \quad \k = (\k_0,\k_1,\k_2,\k_3,\k_4) 
\mapsto (\k_1,\k_2,\k_3,\k_4), 
\]
where the latter space $\C^4$ is equipped with the 
standard (complex) Euclidean inner product. 
For each $i \in \{0,1,2,3,4\}$, let $w_i : \K \to \K$ be 
the orthogonal reflection having $\{\, \k \in \K\,:\, \k_i =0\}$ 
as its reflecting hyperplane with respect to the inner 
product mentioned above. 
Then the group generated by $w_0$, $w_1$, $w_2$, $w_3$, 
$w_4$ is an affine Weyl group of type $D_4^{(1)}$, 
\[
W(D_4^{(1)}) = \la w_0, w_1, w_2, w_3, w_4 \ra 
\car \K.   
\]
corresponding to the Dynkin diagram in Figure 
\ref{fig:dynkin}. 
The reflecting hyperplanes of all reflections in the 
group $W(D_4^{(1)})$ are given by affine linear relations 
\[
\k_i = m, \qquad \k_1 \pm \k_2 \pm \k_3 \pm \k_4 = 2m+1 
\qquad (i \in \{1,2,3,4\}, \, m \in \Z), 
\]
where the signs $\pm$ may be chosen arbitrarily. 
Let $\Wall$ be the union of all these hyperplanes. 
Then the generic condition to be imposed on parameters is 
that $\k$ should lie outside $\Wall$; this is a necessary 
and sufficient condition for $\PVI(\k)$ to admit no 
Riccati solutions \cite{IIS1}. 
\par 
The first main theorem of this paper is concerned with 
the chaotic behavior of $\PVI(\k)$. 
\begin{theorem} \label{thm:chaos} 
Assume that $\k \in \K - \Wall$. 
For any non-elementary loop $\ga \in \pi_1(Z,z)$, the 
Poincar\'{e} return map $\ga_* : M_z(\k) \carl$ along 
the loop $\ga$ is chaotic, that is, there exists a 
$\ga_*$-invariant Borel probability measure $\mu_{\ga}$ 
such that the conditions of Definition $\ref{def:chaos}$ 
are satisfied. 
Moreover there exists an algorithm to calculate the entropy 
$h(\ga) := h_{\mu_{\ga}}(\ga_*)$ of the map $\ga_*$ with 
respect to the measure $\mu_{\ga}$ in terms of a reduced 
word for the loop $\ga$ $($see Theorem $\ref{thm:algorithm})$.  
\end{theorem}
\par
\begin{figure}[t]
\begin{center}
\unitlength 0.1in
\begin{picture}(56.52,11.04)(1.50,-14.45)
%
\special{pn 20}%
\special{sh 0.600}%
\special{ar 3822 893 40 40  0.0000000 6.2831853}%
%
\special{pn 20}%
\special{sh 0.600}%
\special{ar 5255 893 41 41  0.0000000 6.2831853}%
\put(37.7400,-10.7300){\makebox(0,0)[lb]{$z_i$}}%
\put(52.2000,-10.9000){\makebox(0,0)[lb]{$z_j$}}%
%
\special{pn 20}%
\special{pa 4175 953}%
\special{pa 4734 953}%
\special{fp}%
%
\special{pn 20}%
\special{pa 3839 533}%
\special{pa 3791 533}%
\special{fp}%
\special{sh 1}%
\special{pa 3791 533}%
\special{pa 3858 553}%
\special{pa 3844 533}%
\special{pa 3858 513}%
\special{pa 3791 533}%
\special{fp}%
%
\special{pn 20}%
\special{pa 3462 869}%
\special{pa 3462 922}%
\special{fp}%
\special{sh 1}%
\special{pa 3462 922}%
\special{pa 3482 855}%
\special{pa 3462 869}%
\special{pa 3442 855}%
\special{pa 3462 922}%
\special{fp}%
%
\special{pn 20}%
\special{pa 3798 1258}%
\special{pa 3851 1258}%
\special{fp}%
\special{sh 1}%
\special{pa 3851 1258}%
\special{pa 3784 1238}%
\special{pa 3798 1258}%
\special{pa 3784 1278}%
\special{pa 3851 1258}%
\special{fp}%
%
\special{pn 20}%
\special{pa 4472 953}%
\special{pa 4508 953}%
\special{fp}%
\special{sh 1}%
\special{pa 4508 953}%
\special{pa 4441 933}%
\special{pa 4455 953}%
\special{pa 4441 973}%
\special{pa 4508 953}%
\special{fp}%
%
\special{pn 20}%
\special{pa 5231 1265}%
\special{pa 5298 1258}%
\special{fp}%
\special{sh 1}%
\special{pa 5298 1258}%
\special{pa 5230 1245}%
\special{pa 5245 1264}%
\special{pa 5234 1285}%
\special{pa 5298 1258}%
\special{fp}%
%
\special{pn 20}%
\special{pa 5634 917}%
\special{pa 5634 869}%
\special{fp}%
\special{sh 1}%
\special{pa 5634 869}%
\special{pa 5614 936}%
\special{pa 5634 922}%
\special{pa 5654 936}%
\special{pa 5634 869}%
\special{fp}%
%
\special{pn 20}%
\special{pa 5286 526}%
\special{pa 5250 526}%
\special{fp}%
\special{sh 1}%
\special{pa 5250 526}%
\special{pa 5317 546}%
\special{pa 5303 526}%
\special{pa 5317 506}%
\special{pa 5250 526}%
\special{fp}%
%
\special{pn 20}%
\special{pa 3846 1445}%
\special{pa 3803 1438}%
\special{fp}%
\special{sh 1}%
\special{pa 3803 1438}%
\special{pa 3866 1468}%
\special{pa 3856 1447}%
\special{pa 3872 1429}%
\special{pa 3803 1438}%
\special{fp}%
%
\special{pn 20}%
\special{pa 3275 929}%
\special{pa 3275 869}%
\special{fp}%
\special{sh 1}%
\special{pa 3275 869}%
\special{pa 3255 936}%
\special{pa 3275 922}%
\special{pa 3295 936}%
\special{pa 3275 869}%
\special{fp}%
%
\special{pn 20}%
\special{pa 3803 346}%
\special{pa 3851 346}%
\special{fp}%
\special{sh 1}%
\special{pa 3851 346}%
\special{pa 3784 326}%
\special{pa 3798 346}%
\special{pa 3784 366}%
\special{pa 3851 346}%
\special{fp}%
%
\special{pn 20}%
\special{pa 5238 365}%
\special{pa 5303 358}%
\special{fp}%
\special{sh 1}%
\special{pa 5303 358}%
\special{pa 5235 345}%
\special{pa 5250 364}%
\special{pa 5239 385}%
\special{pa 5303 358}%
\special{fp}%
%
\special{pn 20}%
\special{pa 5795 862}%
\special{pa 5802 929}%
\special{fp}%
\special{sh 1}%
\special{pa 5802 929}%
\special{pa 5815 861}%
\special{pa 5796 876}%
\special{pa 5775 865}%
\special{pa 5802 929}%
\special{fp}%
%
\special{pn 20}%
\special{pa 5286 1433}%
\special{pa 5231 1433}%
\special{fp}%
\special{sh 1}%
\special{pa 5231 1433}%
\special{pa 5298 1453}%
\special{pa 5284 1433}%
\special{pa 5298 1413}%
\special{pa 5231 1433}%
\special{fp}%
\put(44.8200,-12.4600){\makebox(0,0)[lb]{$\wp$}}%
%
\special{pn 20}%
\special{ar 5262 898 537 537  3.0378003 6.2831853}%
\special{ar 5262 898 537 537  0.0000000 2.9271389}%
%
\special{pn 20}%
\special{pa 4355 1006}%
\special{pa 4734 1006}%
\special{fp}%
%
\special{pn 20}%
\special{pa 4175 838}%
\special{pa 4907 838}%
\special{fp}%
%
\special{pn 20}%
\special{ar 3822 898 358 358  0.1545646 6.1148227}%
%
\special{pn 20}%
\special{ar 3822 898 547 547  0.1999199 6.2831853}%
%
\special{pn 20}%
\special{pa 4374 898}%
\special{pa 4902 898}%
\special{fp}%
%
\special{pn 20}%
\special{ar 5267 898 372 372  3.3067413 6.2831853}%
\special{ar 5267 898 372 372  0.0000000 3.1415927}%
%
\special{pn 20}%
\special{pa 4614 838}%
\special{pa 4578 838}%
\special{fp}%
\special{sh 1}%
\special{pa 4578 838}%
\special{pa 4645 858}%
\special{pa 4631 838}%
\special{pa 4645 818}%
\special{pa 4578 838}%
\special{fp}%
%
\special{pn 20}%
\special{pa 4619 1006}%
\special{pa 4578 1006}%
\special{fp}%
\special{sh 1}%
\special{pa 4578 1006}%
\special{pa 4645 1026}%
\special{pa 4631 1006}%
\special{pa 4645 986}%
\special{pa 4578 1006}%
\special{fp}%
%
\special{pn 20}%
\special{pa 4470 898}%
\special{pa 4511 893}%
\special{fp}%
\special{sh 1}%
\special{pa 4511 893}%
\special{pa 4442 881}%
\special{pa 4458 899}%
\special{pa 4447 921}%
\special{pa 4511 893}%
\special{fp}%
%
\special{pn 20}%
\special{sh 0.600}%
\special{ar 675 878 38 38  0.0000000 6.2831853}%
%
\special{pn 20}%
\special{sh 0.600}%
\special{ar 2036 878 38 38  0.0000000 6.2831853}%
\put(6.3000,-10.6000){\makebox(0,0)[lb]{$z_i$}}%
\put(20.0000,-10.8000){\makebox(0,0)[lb]{$z_j$}}%
%
\special{pn 20}%
\special{pa 2020 376}%
\special{pa 2082 370}%
\special{fp}%
\special{sh 1}%
\special{pa 2082 370}%
\special{pa 2014 357}%
\special{pa 2029 375}%
\special{pa 2018 396}%
\special{pa 2082 370}%
\special{fp}%
%
\special{pn 20}%
\special{pa 2549 848}%
\special{pa 2556 912}%
\special{fp}%
\special{sh 1}%
\special{pa 2556 912}%
\special{pa 2569 844}%
\special{pa 2550 859}%
\special{pa 2529 848}%
\special{pa 2556 912}%
\special{fp}%
%
\special{pn 20}%
\special{pa 2066 1391}%
\special{pa 2013 1391}%
\special{fp}%
\special{sh 1}%
\special{pa 2013 1391}%
\special{pa 2080 1411}%
\special{pa 2066 1391}%
\special{pa 2080 1371}%
\special{pa 2013 1391}%
\special{fp}%
%
\special{pn 20}%
\special{pa 690 370}%
\special{pa 651 370}%
\special{fp}%
\special{sh 1}%
\special{pa 651 370}%
\special{pa 718 390}%
\special{pa 704 370}%
\special{pa 718 350}%
\special{pa 651 370}%
\special{fp}%
%
\special{pn 20}%
\special{pa 690 1400}%
\special{pa 751 1393}%
\special{fp}%
\special{sh 1}%
\special{pa 751 1393}%
\special{pa 682 1381}%
\special{pa 698 1399}%
\special{pa 687 1420}%
\special{pa 751 1393}%
\special{fp}%
%
\special{pn 20}%
\special{pa 150 890}%
\special{pa 157 957}%
\special{fp}%
\special{sh 1}%
\special{pa 157 957}%
\special{pa 170 889}%
\special{pa 151 904}%
\special{pa 130 893}%
\special{pa 157 957}%
\special{fp}%
%
\special{pn 20}%
\special{pa 1050 523}%
\special{pa 1643 1196}%
\special{fp}%
%
\special{pn 20}%
\special{pa 1643 569}%
\special{pa 1050 1230}%
\special{fp}%
%
\special{pn 20}%
\special{ar 674 877 516 516  0.7237480 5.6219733}%
%
\special{pn 20}%
\special{ar 2042 877 502 502  3.7720777 6.2831853}%
\special{ar 2042 877 502 502  0.0000000 2.4534819}%
\end{picture}%
\end{center}
\caption{An eight-loop (left) and a Pochhammer loop (right)} 
\label{fig:loop} 
\end{figure}
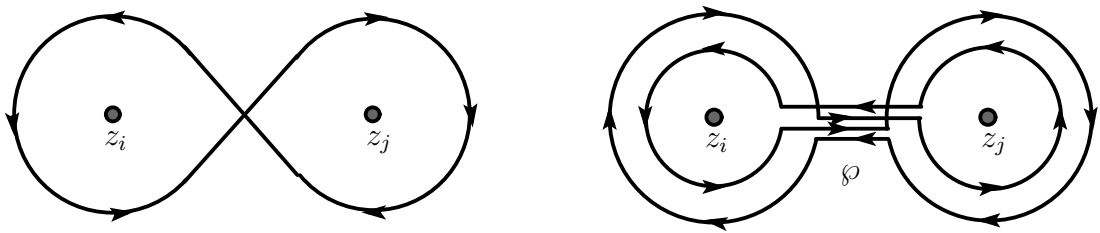

 
The second main theorem is about the periodic 
solutions to $\PVI(\k)$ along a given loop. 
\begin{theorem} \label{thm:per} 
Assume that $\k \in \K - \Wall$. 
For any non-elementary loop $\ga \in \pi_1(Z,z)$, 
the cardinality of the set $\Per_N(\ga ; \k)$ is finite for 
every period $N \in \N$ and grows exponentially 
as $N$ tends to infinity. 
There is an algorithm to count the cardinality exactly as 
well as to determine its exponential growth rate in terms of 
a reduced word for the loop $\ga$ 
$($see Theorem $\ref{thm:algorithm})$. 
\end{theorem} 
\begin{example} \label{ex:per} 
We illustrate Theorems \ref{thm:chaos} and \ref{thm:per} by 
presenting two examples. 
\begin{enumerate} 
\item An {\sl eight-loop} is a loop conjugate to 
$\ga_i \ga_j^{-1}$ for some indices $\{i,j,k\} = \{1,2,3\}$ 
as in Figure \ref{fig:loop} (left). 
For an eight-loop $\ga$ we have 
\[
h(\ga) = \log(3+2\sqrt{2}), \qquad \mathrm{\#} \, \Per_N(\ga ; \k) = 
(3+2\sqrt{2})^N +(3+2\sqrt{2})^{-N} + 4. 
\]
\item A {\sl Pochhammer loop} is a loop conjugate to the 
commutator 
$[\ga_i, \ga_j^{-1}] = \ga_i \ga_j^{-1} \ga_i^{-1} \ga_j$ 
for some indices $\{i,j,k\} = \{1,2,3\}$ as in 
Figure \ref{fig:loop} (right). 
For a Pochhammer loop $\wp$ we have 
\[
h(\wp) = \log(9+4\sqrt{5}), \qquad 
\mathrm{\#} \, \Per_N(\wp ; \k) = 
(9+4\sqrt{5})^N +(9+4\sqrt{5})^{-N} + 4. 
\]
\end{enumerate} 
\end{example}
\par 
As is mentioned in Theorems \ref{thm:chaos} and \ref{thm:per}, 
there are algorithms to calculate the entropy and to count the 
number of periodic solutions exactly. 
In order to describe them we need some preparations concerning 
reduced words for representing loops in terms of the standard 
generators $\ga_1$, $\ga_2$, $\ga_3$. 
\begin{definition} \label{def:reduced} 
For any nontrivial loop $\ga \in \pi_1(Z,z)$, there exists 
an expression 
\begin{equation} \label{eqn:reduced} 
\ga = \ga_{i_1}^{\ve_{i_1}} \ga_{i_2}^{\ve_{i_2}} \cdots 
\ga_{i_m}^{\ve_{i_m}}, 
\end{equation} 
with some positive number $m \in \N$, some indices 
$(i_1, \dots, i_m) \in \{1,2,3\}^m$ and some signs 
$(\ve_{i_1}, \dots, \ve_{i_m}) \in \{\pm1 \}^m$. 
Such an expression is not unique and its length $m$ 
may be reduced by using the relation (\ref{eqn:relation}). 
The expression (\ref{eqn:reduced}) is said to be {\sl reduced} 
if its length $m$ is minimal among all feasible expressions. 
The {\sl length} $\ell_{\pi_1}(\ga)$ of the loop $\ga$ is 
defined to be the length $m$ of a reduced 
expression (\ref{eqn:reduced}) for $\ga$. 
By convention the length of the trivial loop is zero. 
\end{definition} 
\begin{remark} \label{rem:conjugacy}
At this stage we should notice that relevant to our 
discussion is not a loop itself but the conjugacy class 
of a loop. 
Indeed, if two loops $\ga$ and $\ga'$ are conjugate to 
each other, say, $\ga' = \delta \ga \delta^{-1}$ 
for some loop $\delta$, then the corresponding 
Poincar\'{e} return maps are also conjugate to each other as 
$\ga'_* = \delta_* \ga_* \delta_*^{-1}$, and 
hence have the same dynamical properties. 
If $\mu_{\ga}$ is a $\ga_*$-invariant measure asserted in 
Theorem \ref{thm:chaos}, then the push-forward 
$\mu_{\ga'} = (\delta_*)_* \mu_{\ga}$ of the measure 
$\mu_{\ga}$ by the map $\delta_*$ is a desired 
invariant measure for $\ga'_*$. 
As for the sets of periodic points, the loop $\delta$ 
induces a bijection $\delta_* : \Per_N(\ga; \k) \to 
\Per_N(\ga'; \k)$ and hence an equality 
$\mathrm{\#} \Per_N(\ga; \k) = 
\mathrm{\#} \Per_N(\ga'; \k)$. 
So what is relevant is only the conjugacy 
class of a loop. 
\end{remark} \par
This remark leads us to the following definition. 
\begin{definition} \label{def:minimal} 
A loop $\ga \in \pi_1(Z,z)$ is said to be {\sl minimal} if it 
has the minimal length among all loops conjugate to $\ga$, 
namely, if 
$\ell_{\pi_1}(\ga) = \min \{\, \ell_{\pi_1}(\ga')\,:\, \ga' 
\in \pi_1(Z,z) \,\, \mbox{is conjugate to} \,\, \ga \,\}$. 
\end{definition} 
In what follows we may and shall consider minimal loops only 
by replacing a given loop with a minimal 
representative for the conjugacy class of the loop, if it is 
not a minimal loop. 
\par 
In order to give the algorithm, we shall identify $\pi_1(Z,z)$ 
with an index-two subgroup of the universal Coxeter group of 
rank three, that is, the free product of three copies of 
$\Z/2\Z$, 
\[
G = \la\, \si_1, \,\si_2,\, \si_3 \,| \, 
\si_1^2 = \si_2^2 = \si_3^2 = 1 \, \ra 
\simeq (\Z/2\Z) * (\Z/2\Z) * (\Z/2\Z). 
\] 
Any element $\si \in G$ other than the unit element is uniquely 
represented in the form 
\begin{equation} \label{eqn:reduced2}
\si = \si_{i_1} \si_{i_2} \cdots \si_{i_n}, 
\end{equation} 
for some $n \in \N$ and some $n$-tuple of indices 
$(i_1,\dots, i_n) \in \{1,2,3\}^n$ such that every 
neighboring indices $i_{\nu}$ and $i_{\nu+1}$ are distinct. 
The expression (\ref{eqn:reduced2}) is called the {\sl reduced} 
expression of $\si$ and the number $\ell_G(\si) = n$ is 
called the {\sl length} of $\si$, where the unit element is of 
length zero by convention. 
An element of even length is called an {\sl even} element. 
Let $G(2)$ be the subgroup of all even elements in $G$. 
Then there exists an isomorphism of groups 
\begin{equation} \label{eqn:isom} 
\pi_1(Z,z) \to G(2) 
\end{equation} 
sending the basic loops and their inverses as 
\begin{equation} \label{eqn:transl}
\begin{array}{rrr}
\ga_1^{\phantom{-1}} \mapsto \si_1\si_2, \qquad & 
\ga_2^{\phantom{-1}} \mapsto \si_2\si_3, \qquad &
\ga_3^{\phantom{-1}} \mapsto \si_3\si_1, \\[2mm]
\ga_1^{-1} \mapsto \si_2 \si_1, \qquad &
\ga_2^{-1} \mapsto \si_3 \si_2, \qquad &
\ga_3^{-1} \mapsto \si_1 \si_3. 
\end{array}
\end{equation}
Given an expression $(\ref{eqn:reduced})$ of a loop 
$\ga \in \pi_1(Z,z)$, make the replacement of alphabets 
\[
\{\ga_1^{\pm1}, \ga_2^{\pm1}, \ga_3^{\pm1}\} \to 
\{\si_1, \si_2, \si_3\}
\]
according to the rule $(\ref{eqn:transl})$. 
If the expression (\ref{eqn:reduced}) is reduced in 
$\pi_1(Z,z)$, then the resulting word is also reduced 
in $G$. 
In particular the reduced expression (\ref{eqn:reduced}) 
is unique for a given loop $\ga$ and one has 
$\ell_G(\si) = 2 \ell_{\pi_1}(\ga)$, where $\si \in G(2)$ 
is the element corresponding to the loop $\ga$. 
\par 
Recall that any Coxeter group admits its geometric 
representation \cite{Humphreys}. 
We apply this construction to our particular group $G$. 
Let $V$ be the $3$-dimensional vector space spanned by 
basis vectors $e_1$, $e_2$, $e_3$, endowed with a 
nondegenerate symmetric bilinear form 
\begin{equation} \label{eqn:B}
B(e_i, e_j) = \left\{ 
\begin{array}{ll}
\phantom{-}1 \qquad & (i = j), \\[1mm]
-1 \qquad & (i \neq j). 
\end{array}
\right. 
\end{equation} 
For each $i \in \{1,2,3\}$ we can define an orthogonal 
reflection $r_i : V \to V$ by the rule
\begin{equation} \label{eqn:ri}
r_i(v) := v - 2 B(e_i,v)\, e_i \qquad (v \in V). 
\end{equation}
Note that $r_i$ sends $e_i$ to its negative $-e_i$ while 
fixing all the vectors orthogonal to $e_i$ relative to the 
bilinear form $B$. 
It is known that there is a unique injective homomorphism 
$\mathrm{GR} : G \to \mathrm{O}_B(V)$ such that 
$\mathrm{GR}(\si_i) = r_i$ for $i \in \{1,2,3\}$, 
where $\mathrm{O}_B(V)$ is the group of orthogonal 
transformations on $(V,B)$. 
Identified with its image $\mathrm{GR}(G)$, the group 
$G$ can be thought of as a reflection group acting on $(V, B)$. 
The faithful representation $\mathrm{GR} : G \to 
\mathrm{O}_B(V)$ is called the geometric representation of $G$. 
For each $i \in \{1,2,3\}$ we define an endomorphism 
$s_i : V \to V$ by 
\begin{equation} \label{eqn:esi}
s_i(v) := \dfrac{v+r_i(v)}{2} = v- B(e_i, v) \, e_i 
\qquad (v \in V), 
\end{equation} 
and make the following definition. 
\begin{definition} \label{def:alpha}
Given a loop $\ga \in \pi_1(Z,z)$, choose a minimal 
representative for the conjugacy class of $\ga$ and call it 
$\ga$ again. 
Take the reduced expression of $\ga$ as in 
(\ref{eqn:reduced}). 
Make the change of alphabets $\{\ga_1^{\pm1}, \ga_2^{\pm1}, 
\ga_3^{\pm1}\} \to \{\si_1, \si_2, \si_3\}$ according to 
the rule (\ref{eqn:transl}) to obtain the corresponding 
element $\si \in G(2)$, together with its reduced 
expression as in (\ref{eqn:reduced2}). 
To the indices $(i_1,\dots,i_n)$ in (\ref{eqn:reduced2}), 
associate an endomorphism 
$s_{\ga} := s_{i_n} \cdots s_{i_2} s_{i_1} \in \End \, V$. 
Finally, take its trace 
\begin{equation} \label{eqn:alpha} 
\a(\ga) = \Tr[\,s_{\ga} : V \to V\,]. 
\end{equation}
\end{definition} 
\par 
We are now in a position to give the algorithm to calculate 
the entropy and to count the number of periodic points, 
which complete the statements of Theorems \ref{thm:chaos} and 
\ref{thm:per}. 
\begin{theorem} \label{thm:algorithm} 
Assume that $\k \in \K-\Wall$ and let $\ga \in \pi_1(Z,z)$ be 
any non-elementary loop.  
Then the number $\a(\ga)$ defined in $(\ref{eqn:alpha})$ is an 
even integer not smaller than $6$, with the equality 
$\a(\ga) = 6$ if and only if $\ga$ is an eight-loop as in 
Example $\ref{ex:per}$. 
Put 
\begin{equation} \label{eqn:lambda}
\l(\ga) := \dfrac{1}{2}
\left\{\a(\ga) + \sqrt{\a(\ga)^2-4}\right\}. 
\end{equation}
\begin{enumerate}
\item The measure-theoretic entropy 
$h(\ga) := h_{\mu_{\ga}}(\ga_*)$ of the Poincar\'e return map 
$\ga_* : M_z(\k) \carl$ with respect to the invariant 
measure $\mu_{\ga}$ mentioned in Theorem $\ref{thm:chaos}$ 
is given by 
\[
h(\ga) = \log \lambda(\ga). 
\]
\item The cardinality of the set $\Per_N(\ga;\k)$ is given by 
\[
\mathrm{\#} \Per_N(\ga;\k) = \l(\ga)^N + 
\l(\ga)^{-N} + 4 \qquad (N \in \N).  
\]
In particular its exponential growth rate is given by 
$\l(\ga)$. 
\end{enumerate} 
\end{theorem}
\begin{remark} \label{rem:mostelem} 
Theorem \ref{thm:algorithm} implies that 
for any non-elementary loop $\ga \in \pi_1(Z,z)$, we have 
\[
\l(\ga) \ge 3 + 2 \sqrt{2}, \qquad 
h(\ga) \ge \log(3+2\sqrt{2}), 
\] 
with the equalities if and only if $\ga$ is an eight-loop. 
In this sense the eight-loops are the most ``elementary" loops 
among all non-elementary loops in $Z$. 
On the other hand, one may ask what happens with the Poincar\'e 
return map $\ga_* : M_z(\k) \carl$ when the loop $\ga$ is 
elementary. 
In this case it turns out that $\ga_*$ preserves a certain 
analytic fibration $M_z(\k) \to \C$ and exhibits an essentially 
$1$-dimensional dynamical behavior. 
Hence $\ga_*$ is not so interesting or too elementary 
from the standpoint of chaotic dynamical systems. 
See Remark \ref{rem:proof} for more information. 
\end{remark} 
\begin{remark} \label{rem:lyapunov} 
There exists a standard complex area form $\omega_z(\k)$ on 
$M_z(\k)$ such that the Poincar\'e return map $\ga_*$ is 
area-preserving for every loop $\ga \in \pi_1(Z,z)$, where 
we refer to Remark \ref{rem:area} for the description of 
$\omega_z(\k)$. 
Hence the Lyapunov exponents $L_{\pm}(\ga)$ of $\ga_*$ 
satisfy the relation $L_-(\ga) = - L_+(\ga)$. 
Moreover the positive Lyapunov exponent admits an 
estimate $L_+(\ga) \ge \frac{1}{8}\log \l(\ga)$. 
We refer to Remark \ref{rem:proof} for the derivation of 
this estimate. 
\end{remark} 
\begin{remark} \label{rem:generic} 
In this article we restrict our attention to the 
generic case $\k \in \K-\Wall$ only, leaving 
the nongeneric case $\k \in \Wall$ untouched. 
The difference between the generic case and the nongeneric 
case lies in the fact that the Riemann-Hilbert correspondence 
to be used in the proof becomes a biholomorphism in the former 
case, while it gives an analytic minimal resolution of Klein 
singularities in the latter case (see Remark \ref{rem:RH}). 
The presence of singularities would make the treatment of 
the nongeneric case more complicated. 
However it is expected that the basic strategy developed in this 
article will be effective also in the nongeneric case. 
The relevant discussion will be made elsewhere. 
\end{remark} 
\par 
The plan of this article is as follows: 
$\PVI(\k)$ is formulated as a flow, {\sl Painlev\'e flow}, on 
a moduli space of stable parabolic connections in \S\ref{sec:mspc}. 
It is conjugated to an isomonodromic flow on a moduli space of 
monodromy representations via a Riemann-Hilbert correspondence 
in \S\ref{sec:RHC}. 
The moduli space of monodromy representations is identified 
with an affine cubic surface and each Poincar\'e return map 
for $\PVI(\k)$ is conjugated to a biregular automorphism of 
the affine cubic in \S\ref{sec:cubic}. 
This map is extended to a birational map on the compactified 
projective cubic surface and some basic properties of it 
are studied in \S\ref{sec:dynamics}. 
The induced cohomological action of the birational map 
is investigated in \S\ref{sec:cohomology}. 
After these preliminaries, the ergodic properties of our 
dynamical system are established by applying some recent deep 
results from birational surface dynamics in \S\ref{sec:ergodic}. 
Moreover the number of periodic points of the birational map 
is counted by using the Lefschetz fixed point formula in 
\S\ref{sec:periodic}. 
Then, back to the original phase space of $\PVI(\k)$ in 
\S\ref{sec:back}, we arrive at our final goals, that is, at 
the ergodic properties of the Poincar\'e return map and the 
exact number of periodic solutions to $\PVI(\k)$ of any 
period along a given loop. 
\section{Moduli Space of Stable Parabolic Connections} 
\label{sec:mspc}
In order to describe the fibration (\ref{eqn:phase}), we first 
construct an auxiliary fibration 
$\pi_{\k} : \M(\k) \to T$ over the configuration 
space of mutually distinct, ordered, three points in $\C$, 
\[
T = \{\, t=(t_1,t_2,t_3) \in \C^3\,:\, t_i \neq t_j \,\, 
\text{for} \,\, i \neq j \, \},  
\]
and then reduce it to the original fibration (\ref{eqn:phase}). 
We put the fourth point $t_4$ at infinity. 
Given any $(t,\k) \in T \times \K$, a $(t,\k)$-parabolic 
connection is a quadruple $Q = (E,\nabla,\psi,l)$ such that
\begin{enumerate} 
\item $E$ is a rank $2$ algebraic vector bundle of degree 
$-1$ over $\P^1$, 
\item $\nabla : E \to E \otimes \Omega^1_{\P^1}(D_t)$ is a Fuchsian 
connection with pole divisor $D_t = t_1 + t_2 + t_3 + t_4$ and 
Riemann scheme as in Table \ref{tab:riemann}, where 
$t_4 = \infty$ as mentioned above, 
\item $\psi : \det \, E \to \O_{\P^1}(-t_4)$ is a horizontal 
isomorphism called a determinantal structure, where 
$\O_{\P^1}(-t_4)$ is equipped with the connection induced from 
$d : \O_{\P^1} \to \Omega^1_{\P^1}$, 
\item $l = (l_1,l_2,l_3,l_4)$ is a parabolic structure, namely, 
$l_i$ is an eigenline of $\Res_{t_i}(\nabla) \in \End(E_{t_i})$ 
corresponding to eigenvalue $\l_i$ (whose minus is the first 
exponent $-\l_i$ in Table \ref{tab:riemann}). 
\end{enumerate} 
\begin{table}[t]
\begin{center} 
\begin{tabular}{|c||c|c|c|c|}
\hline
\vspace{-4mm} & & & & \\
singularities & $t_1$ & $t_2$ & $t_3$ & $t_4$ \\[1mm]
\hline
\vspace{-4mm} & & & & \\
first exponent & $-\l_1$ & $-\l_2$ & $-\l_3$ & $-\l_4$ \\[1mm]
\hline
\vspace{-4mm} & & & & \\
second exponent & $\l_1$ & $\l_2$  & $\l_3$ & $\l_4-1$ \\[1mm]
\hline
\vspace{-4mm} & & & & \\
difference  & $\k_1$ & $\k_2$ & $\k_3$ & $\k_4$ \\[1mm]  
\hline
\end{tabular}
\end{center}
\caption{Riemann scheme: $\k_i$ is the difference of the second 
exponent from the first.} 
\label{tab:riemann}
\end{table}
There exists a concept of stability for parabolic connections, 
with which the geometric invariant theory \cite{Mumford} can 
be worked out to establish the following theorem \cite{IIS1,IIS2}.  
\begin{theorem} \label{thm:moduli} 
For any $(t,\k) \in T \times \K$ there exists a fine moduli scheme 
$\M_t(\k)$ of stable $(t,\k)$-parabolic connections.  
The moduli space $\M_t(\k)$ is a smooth, irreducible, 
quasi-projective surface. 
As a relative setting over $T$, for any $\k \in \K$, there exists 
a family of moduli spaces 
\begin{equation} \label{eqn:family} 
\pi_{\k} : \M(\k) \rightarrow T 
\end{equation}
such that the projection $\pi_{\k}$ is a smooth morphism 
having fiber $\M_t(\k)$ over $t \in T$. 
\end{theorem} 
\par 
In \cite{IIS1,IIS2} the moduli space $\M_t(\k)$ is compactified 
into a moduli space of stable parabolic phi-connections. 
Given any $(t,\k) \in T \times \K$, a {\sl parabolic phi-connection} 
is roughly speaking a sextuple of data 
$Q = (E_1,E_2,\phi,\nabla,\psi,l)$ consisting of 
\begin{enumerate}
\item a variant of connection $\nabla : E_1 \to E_2 \otimes 
\Omega^1_{\P^1}(D_t)$ over rank $2$, degree $-1$ bundles on $\P^1$, 
\item an $\O_{\P^1}$-homomorphism $\phi : E_1 \to E_2$ 
(called a {\sl phi-operator}), which may be {\sl degenerate} 
or non-isomorphic, satisfying a generalized Leibniz rule 
\[
\nabla(fs) = \phi(s) \otimes df + f \nabla(s), 
\qquad (s \in E_1, \,  f \in \O_{\P^1}), 
\]
\item extra data of a determinantal structure $\psi$ and a 
parabolic structure $l$. 
\end{enumerate} 
We refer to \cite{IIS1,IIS2} for the complete definition. 
Very roughly the idea of compactification is as follows: 
If a parabolic connection is regarded as a ``matrix-valued 
Schr\"odinger operator", then a parabolic phi-connection may 
be thought of as a matrix-valued Schr\"odinger operator 
with a ``matrix-valued Planck constant" $\phi$ which may 
be degenerate, namely, may be semi-classical. 
Then the moduli space $\M_t(\k)$ can be compactified by 
adding some semi-classical objects, that is, some parabolic 
phi-connections with degenerate phi-operator $\phi$. 
\par 
There exists a concept of {\sl stability} for parabolic 
phi-connections, with which geometric invariant theory 
can be worked out to establish the following theorem 
\cite{IIS1,IIS2}. 
\begin{theorem} \label{thm:sppc} 
For any $(t,\k) \in T \times \K$ there exists a coarse moduli 
scheme $\ol{\M}_t(\k)$ of stable 
parabolic phi-connections.  
The moduli space $\ol{\M}_t(\k)$ is a smooth, irreducible, 
projective surface, having a unique effective anti-canonical 
divisor $\Y_t(\k)$. 
Under the natural embedding 
\[
\M_t(\k) \hookrightarrow \ol{\M}_t(\k), \qquad 
(E,\nabla,\psi,l) \mapsto (E,E,\id,\nabla,\psi,l), 
\]
the space $\M_t(\k)$ is exactly the locus of $\ol{\M}_t(\k)$ 
where the phi-operator $\phi$ is isomorphic, and so 
\[
\M_t(\k) = \ol{\M}_t(\k) - \Y_t(\k). 
\]
\end{theorem} 
\par
\begin{figure}[t]
\begin{center}
\unitlength 0.1in
\begin{picture}(40.82,36.49)(1.30,-38.09)
%
\special{pn 20}%
\special{pa 3114 698}%
\special{pa 3114 2642}%
\special{fp}%
%
\special{pn 20}%
\special{pa 1365 689}%
\special{pa 1365 2633}%
\special{fp}%
%
\special{pn 20}%
\special{pa 801 689}%
\special{pa 3707 689}%
\special{pa 3707 2642}%
\special{pa 801 2642}%
\special{pa 801 689}%
\special{fp}%
%
\special{pn 20}%
\special{pa 801 1097}%
\special{pa 3717 1097}%
\special{fp}%
%
\special{pn 20}%
\special{pa 1763 1282}%
\special{pa 2142 1680}%
\special{dt 0.054}%
\special{pa 2142 1680}%
\special{pa 2142 1680}%
\special{dt 0.054}%
%
\special{pn 20}%
\special{pa 2337 1272}%
\special{pa 2716 1670}%
\special{dt 0.054}%
\special{pa 2716 1670}%
\special{pa 2716 1670}%
\special{dt 0.054}%
%
\special{pn 20}%
\special{pa 2930 1291}%
\special{pa 3309 1690}%
\special{dt 0.054}%
\special{pa 3309 1690}%
\special{pa 3309 1690}%
\special{dt 0.054}%
%
\special{pn 20}%
\special{pa 1170 1282}%
\special{pa 1550 1680}%
\special{dt 0.054}%
\special{pa 1550 1680}%
\special{pa 1550 1680}%
\special{dt 0.054}%
%
\special{pn 20}%
\special{pa 1559 2059}%
\special{pa 1170 2438}%
\special{dt 0.054}%
\special{pa 1170 2438}%
\special{pa 1170 2438}%
\special{dt 0.054}%
%
\special{pn 20}%
\special{pa 2133 2059}%
\special{pa 1744 2438}%
\special{dt 0.054}%
\special{pa 1744 2438}%
\special{pa 1744 2438}%
\special{dt 0.054}%
%
\special{pn 20}%
\special{pa 2726 2050}%
\special{pa 2337 2429}%
\special{dt 0.054}%
\special{pa 2337 2429}%
\special{pa 2337 2429}%
\special{dt 0.054}%
%
\special{pn 20}%
\special{pa 3299 2040}%
\special{pa 2910 2419}%
\special{dt 0.054}%
\special{pa 2910 2419}%
\special{pa 2910 2419}%
\special{dt 0.054}%
\put(37.9500,-18.0600){\makebox(0,0)[lb]{$\M_t(\k)$}}%
%
\special{pn 20}%
\special{sh 0.600}%
\special{ar 3114 1486 34 34  0.0000000 6.2831853}%
%
\special{pn 20}%
\special{sh 0.600}%
\special{ar 3124 2225 33 33  0.0000000 6.2831853}%
%
\special{pn 20}%
\special{sh 0.600}%
\special{ar 2531 2244 35 35  0.0000000 6.2831853}%
%
\special{pn 20}%
\special{sh 0.600}%
\special{ar 2541 1495 35 35  0.0000000 6.2831853}%
%
\special{pn 20}%
\special{sh 0.600}%
\special{ar 1958 1505 35 35  0.0000000 6.2831853}%
%
\special{pn 20}%
\special{sh 0.600}%
\special{ar 1958 2244 35 35  0.0000000 6.2831853}%
%
\special{pn 20}%
\special{sh 0.600}%
\special{ar 1365 1495 35 35  0.0000000 6.2831853}%
%
\special{pn 20}%
\special{sh 0.600}%
\special{ar 1365 2244 35 35  0.0000000 6.2831853}%
%
\special{pn 8}%
\special{pa 149 3799}%
\special{pa 3444 3799}%
\special{fp}%
%
\special{pn 8}%
\special{pa 898 2953}%
\special{pa 4212 2953}%
\special{fp}%
\put(10.5400,-3.3000){\makebox(0,0)[lb]{$\Y_t(\k)$ : vertical leaves}}%
%
\special{pn 8}%
\special{pa 1480 390}%
\special{pa 1916 846}%
\special{fp}%
\special{sh 1}%
\special{pa 1916 846}%
\special{pa 1884 784}%
\special{pa 1879 807}%
\special{pa 1855 812}%
\special{pa 1916 846}%
\special{fp}%
\put(36.2900,-31.7600){\makebox(0,0)[lb]{$T$}}%
%
\special{pn 20}%
\special{sh 0.600}%
\special{ar 2084 3138 52 52  0.0000000 6.2831853}%
\put(20.9000,-33.2000){\makebox(0,0)[lb]{$t$}}%
%
\special{pn 8}%
\special{pa 898 2953}%
\special{pa 130 3799}%
\special{fp}%
%
\special{pn 8}%
\special{pa 4203 2953}%
\special{pa 3435 3809}%
\special{fp}%
%
\special{pn 20}%
\special{ar 2375 3313 457 253  4.3193680 6.2831853}%
\special{ar 2375 3313 457 253  0.0000000 3.9349651}%
%
\special{pn 20}%
\special{pa 2220 3089}%
\special{pa 2162 3100}%
\special{fp}%
\special{sh 1}%
\special{pa 2162 3100}%
\special{pa 2231 3107}%
\special{pa 2214 3090}%
\special{pa 2224 3068}%
\special{pa 2162 3100}%
\special{fp}%
%
\special{pn 20}%
\special{pa 2337 3556}%
\special{pa 2424 3556}%
\special{fp}%
\special{sh 1}%
\special{pa 2424 3556}%
\special{pa 2357 3536}%
\special{pa 2371 3556}%
\special{pa 2357 3576}%
\special{pa 2424 3556}%
\special{fp}%
\put(27.4500,-36.1400){\makebox(0,0)[lb]{$\beta$}}%
\put(38.1400,-11.5600){\makebox(0,0)[lb]{$E_0$}}%
\put(12.3800,-6.4000){\makebox(0,0)[lb]{$E_1$}}%
\put(18.6000,-6.4900){\makebox(0,0)[lb]{$E_2$}}%
\put(24.4400,-6.4900){\makebox(0,0)[lb]{$E_3$}}%
\put(30.1700,-6.4900){\makebox(0,0)[lb]{$E_4$}}%
%
\special{pn 20}%
\special{ar 2240 1787 827 204  5.6711699 6.2831853}%
\special{ar 2240 1787 827 204  0.0000000 3.9825888}%
%
\special{pn 20}%
\special{sh 0.600}%
\special{ar 1695 1642 52 52  0.0000000 6.2831853}%
%
\special{pn 20}%
\special{sh 0.600}%
\special{ar 2793 1631 53 53  0.0000000 6.2831853}%
%
\special{pn 20}%
\special{pa 2158 1991}%
\special{pa 2314 1991}%
\special{fp}%
\special{sh 1}%
\special{pa 2314 1991}%
\special{pa 2247 1971}%
\special{pa 2261 1991}%
\special{pa 2247 2011}%
\special{pa 2314 1991}%
\special{fp}%
%
\special{pn 20}%
\special{pa 2968 1690}%
\special{pa 2880 1651}%
\special{fp}%
\special{sh 1}%
\special{pa 2880 1651}%
\special{pa 2933 1696}%
\special{pa 2929 1673}%
\special{pa 2949 1660}%
\special{pa 2880 1651}%
\special{fp}%
\special{pa 2880 1651}%
\special{pa 2880 1651}%
\special{fp}%
\put(21.0300,-19.2200){\makebox(0,0)[lb]{$\beta_*$}}%
%
\special{pn 20}%
\special{pa 1954 707}%
\special{pa 1954 1895}%
\special{fp}%
%
\special{pn 20}%
\special{pa 1954 2027}%
\special{pa 1954 2639}%
\special{fp}%
%
\special{pn 20}%
\special{pa 2530 695}%
\special{pa 2530 1907}%
\special{fp}%
%
\special{pn 20}%
\special{pa 2530 2027}%
\special{pa 2530 2639}%
\special{fp}%
%
\special{pn 8}%
\special{pa 1480 400}%
\special{pa 1660 1070}%
\special{fp}%
\special{sh 1}%
\special{pa 1660 1070}%
\special{pa 1662 1000}%
\special{pa 1646 1018}%
\special{pa 1623 1011}%
\special{pa 1660 1070}%
\special{fp}%
\end{picture}%
\end{center}
\caption{Poincare section on the space of initial conditions}
\label{fig:ps} 
\end{figure}
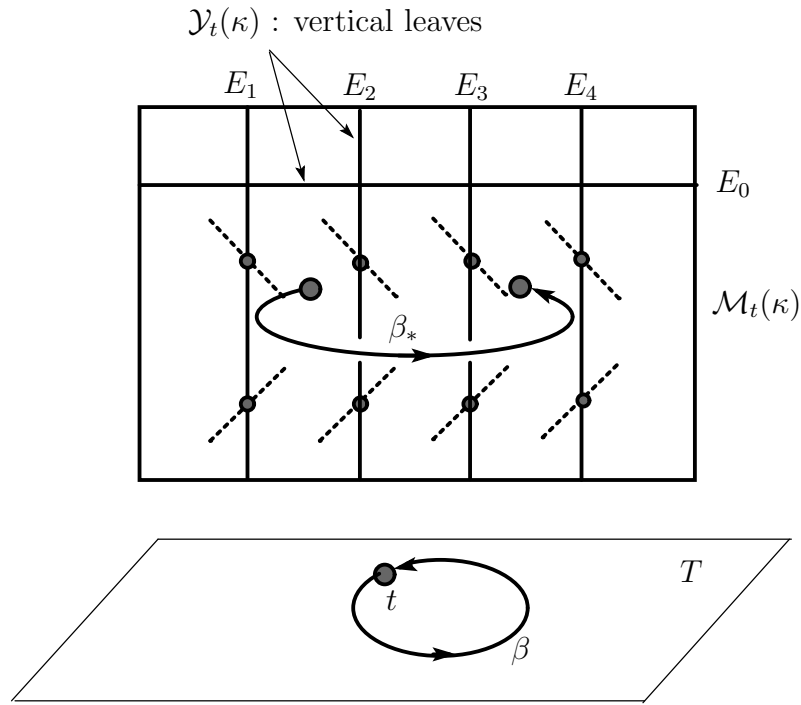

 
The divisor $\Y_t(\k)$ on $\ol{\M}_t(\k)$ is called the 
{\sl vertical leaves} at time $t$. 
There is the following realization of our moduli spaces 
\cite{IIS1,IIS2} (see Figure \ref{fig:ps}). 
\begin{theorem} \label{thm:hirzebruch} 
The compactified moduli space $\ol{\M}_t(\k)$ is isomorphic to an 
$8$-point blow-up of the Hirzebruch surface 
${\mit\Sigma}_2 \to \P^1$ of degree $2$, blown up at certain two 
points on each fiber over the points 
$t_1, \, t_2, \, t_3, \, t_4 \in \P^1$. 
The unique effective anti-canonical divisor on 
$\ol{\M}_t(\k)$ is given by 
\[
\Y_t(\k) = 2 E_0 + E_1 + E_2 + E_3 + E_4, 
\]
where $E_0$ is the strict transform of the section at infinity 
of the fibration $\varSigma_2 \to \P^1$, while $E_1$, $E_2$, 
$E_3$, $E_4$ are the strict transforms of the fibers over 
$t_1$, $t_2$, $t_3$, $t_4$, respectively. 
\end{theorem}
\begin{remark} \label{rem:area} 
There is a meromorphic $2$-form $\omega_t(\k)$ on $\ol{\M}_t(\k)$, 
holomorphic and nondegenerate on $\M_t(\k)$, whose pole divisor 
is given by the vertical leaves $\Y_t(\k)$ \cite{IIS1,IIS2,STT,Sakai}. 
It is unique up to constant multiples. 
This complex area form is just what we have mentioned in 
Remark \ref{rem:lyapunov}. 
A further description of the area form $\omega_t(\k)$ will be 
given in Remark \ref{rem:area2}. 
\end{remark} 
\par
Now the fibration (\ref{eqn:phase}) is defined to be the pull-back 
of the fibration (\ref{eqn:family}) by an injection 
\[
\iota : Z \hookrightarrow T, \quad z \mapsto (0,z,1),  
\]
The group $\Aff(\C)$ of affine linear transformations on 
$\C$ acts diagonally on the configuration space $T$ 
and the quotient space $T/\Aff(\C)$ is isomorphic 
to $Z$, with the quotient map given by 
\begin{equation} \label{eqn:reduction} 
r : T \to Z, \quad t = (t_1,t_2,t_3) \mapsto 
z = \dfrac{t_2-t_1}{t_3-t_1}. 
\end{equation}
The map $r$ yields a trivial $\Aff(\C)$-bundle 
structure of $T$ over $Z$ and the fibration (\ref{eqn:family}) 
is in turn the pull-back of the fibration (\ref{eqn:phase}) 
by the map $r$. 
Hence we have a commutative diagram 
\begin{equation} \label{cd:reduction} 
\begin{CD}
\M(\k) @>>> M(\k) \\
 @V \pi_{\k} VV  @VV \pi_{\k} V \\
T @>> r > Z. 
\end{CD}
\end{equation}
In \cite{IIS1,IIS2} the Painlev\'e dynamical system $\PVI(\k)$ is 
formulated as a holomorphic uniform foliations on the fibration 
(\ref{eqn:family}) which is compatible with the diagram 
(\ref{cd:reduction}). 
Thus the Poincar\'e section (\ref{eqn:PSz}) is reformulated as a 
group homomorphism
\begin{equation} \label{eqn:PSt} 
\PS_t(\k) \,\,:\,\, \pi_1(T,t) \to \Aut\,\M_t(\k), 
\end{equation}
a visual image of which is given in Figure \ref{fig:ps}. 
\par
\begin{figure}[t]
\begin{center}
\unitlength 0.1in
\begin{picture}(63.50,12.85)(1.20,-14.60)
%
\special{pn 20}%
\special{pa 391 403}%
\special{pa 2990 403}%
\special{fp}%
%
\special{pn 20}%
\special{pa 383 1180}%
\special{pa 2981 1180}%
\special{fp}%
%
\special{pn 20}%
\special{pa 2457 1170}%
\special{pa 2457 423}%
\special{fp}%
\special{sh 1}%
\special{pa 2457 423}%
\special{pa 2437 490}%
\special{pa 2457 476}%
\special{pa 2477 490}%
\special{pa 2457 423}%
\special{fp}%
%
\special{pn 20}%
\special{pa 1880 529}%
\special{pa 1880 423}%
\special{fp}%
\special{sh 1}%
\special{pa 1880 423}%
\special{pa 1860 490}%
\special{pa 1880 476}%
\special{pa 1900 490}%
\special{pa 1880 423}%
\special{fp}%
%
\special{pn 20}%
\special{pa 933 1180}%
\special{pa 955 1159}%
\special{pa 978 1139}%
\special{pa 1000 1118}%
\special{pa 1023 1097}%
\special{pa 1046 1077}%
\special{pa 1069 1057}%
\special{pa 1093 1037}%
\special{pa 1117 1017}%
\special{pa 1142 998}%
\special{pa 1167 978}%
\special{pa 1193 959}%
\special{pa 1219 941}%
\special{pa 1246 923}%
\special{pa 1274 905}%
\special{pa 1303 888}%
\special{pa 1333 871}%
\special{pa 1363 855}%
\special{pa 1395 839}%
\special{pa 1428 824}%
\special{pa 1462 810}%
\special{pa 1496 796}%
\special{pa 1531 782}%
\special{pa 1566 768}%
\special{pa 1600 754}%
\special{pa 1634 739}%
\special{pa 1667 725}%
\special{pa 1699 709}%
\special{pa 1729 693}%
\special{pa 1757 675}%
\special{pa 1783 657}%
\special{pa 1806 637}%
\special{pa 1826 615}%
\special{pa 1843 592}%
\special{pa 1856 567}%
\special{pa 1866 540}%
\special{pa 1871 511}%
\special{pa 1872 479}%
\special{pa 1871 452}%
\special{sp}%
%
\special{pn 20}%
\special{pa 1871 1190}%
\special{pa 1850 1165}%
\special{pa 1829 1140}%
\special{pa 1807 1115}%
\special{pa 1785 1091}%
\special{pa 1763 1068}%
\special{pa 1740 1046}%
\special{pa 1716 1025}%
\special{pa 1691 1005}%
\special{pa 1666 987}%
\special{pa 1639 971}%
\special{pa 1611 956}%
\special{pa 1582 943}%
\special{pa 1552 932}%
\special{pa 1522 921}%
\special{pa 1491 911}%
\special{pa 1460 901}%
\special{pa 1450 898}%
\special{sp}%
%
\special{pn 20}%
\special{pa 959 500}%
\special{pa 933 432}%
\special{fp}%
\special{sh 1}%
\special{pa 933 432}%
\special{pa 938 501}%
\special{pa 952 482}%
\special{pa 975 487}%
\special{pa 933 432}%
\special{fp}%
%
\special{pn 20}%
\special{pa 959 491}%
\special{pa 971 522}%
\special{pa 982 552}%
\special{pa 995 582}%
\special{pa 1009 610}%
\special{pa 1025 638}%
\special{pa 1042 663}%
\special{pa 1062 688}%
\special{pa 1084 710}%
\special{pa 1108 731}%
\special{pa 1133 751}%
\special{pa 1159 770}%
\special{pa 1186 789}%
\special{pa 1214 806}%
\special{pa 1242 824}%
\special{pa 1252 830}%
\special{sp}%
\put(17.5100,-3.5500){\makebox(0,0)[lb]{$t_i'$}}%
\put(17.8000,-14.7000){\makebox(0,0)[lb]{$t_j$}}%
\put(23.9600,-14.7100){\makebox(0,0)[lb]{$t_k$}}%
\put(4.3000,-8.8000){\makebox(0,0)[lb]{$\b_i$}}%
\put(8.5600,-3.4500){\makebox(0,0)[lb]{$t_j'$}}%
\put(23.8800,-3.4500){\makebox(0,0)[lb]{$t_k'$}}%
\put(8.4700,-14.8100){\makebox(0,0)[lb]{$t_i$}}%
%
\special{pn 13}%
\special{pa 3760 805}%
\special{pa 6068 795}%
\special{fp}%
%
\special{pn 20}%
\special{ar 4550 805 430 430  1.5946014 6.2831853}%
\special{ar 4550 805 430 430  0.0000000 1.5707963}%
%
\special{pn 20}%
\special{sh 0.600}%
\special{ar 4980 810 36 36  0.0000000 6.2831853}%
%
\special{pn 20}%
\special{sh 0.600}%
\special{ar 4120 800 35 35  0.0000000 6.2831853}%
%
\special{pn 20}%
\special{sh 0.600}%
\special{ar 5710 790 36 36  0.0000000 6.2831853}%
%
\special{pn 20}%
\special{pa 4960 935}%
\special{pa 4970 885}%
\special{fp}%
\special{sh 1}%
\special{pa 4970 885}%
\special{pa 4937 946}%
\special{pa 4960 937}%
\special{pa 4977 954}%
\special{pa 4970 885}%
\special{fp}%
%
\special{pn 20}%
\special{pa 4590 380}%
\special{pa 4530 380}%
\special{fp}%
\special{sh 1}%
\special{pa 4530 380}%
\special{pa 4597 400}%
\special{pa 4583 380}%
\special{pa 4597 360}%
\special{pa 4530 380}%
\special{fp}%
%
\special{pn 20}%
\special{pa 4530 1235}%
\special{pa 4600 1235}%
\special{fp}%
\special{sh 1}%
\special{pa 4600 1235}%
\special{pa 4533 1215}%
\special{pa 4547 1235}%
\special{pa 4533 1255}%
\special{pa 4600 1235}%
\special{fp}%
\put(56.6000,-10.7000){\makebox(0,0)[lb]{$t_k$}}%
\put(50.2000,-10.8000){\makebox(0,0)[lb]{$t_j$}}%
\put(39.5000,-10.4500){\makebox(0,0)[lb]{$t_i$}}%
\put(60.7000,-4.3000){\makebox(0,0)[lb]{$\hat{\C}$}}%
%
\special{pn 13}%
\special{pa 3550 180}%
\special{pa 6470 180}%
\special{pa 6470 1460}%
\special{pa 3550 1460}%
\special{pa 3550 180}%
\special{fp}%
\put(61.5000,-8.7000){\makebox(0,0)[lb]{$\hat{\bR}$}}%
\put(1.2000,-4.6000){\makebox(0,0)[lb]{$T$}}%
\put(1.3000,-12.4000){\makebox(0,0)[lb]{$T$}}%
%
\special{pn 20}%
\special{pa 4130 680}%
\special{pa 4120 730}%
\special{fp}%
\special{sh 1}%
\special{pa 4120 730}%
\special{pa 4153 669}%
\special{pa 4130 678}%
\special{pa 4113 661}%
\special{pa 4120 730}%
\special{fp}%
\end{picture}%
\end{center}
\caption{Basic braid $\b_i$ in $T$ and the corresponding 
movement of $t$ in $\hat{\C}$}
\label{fig:braid}
\end{figure}
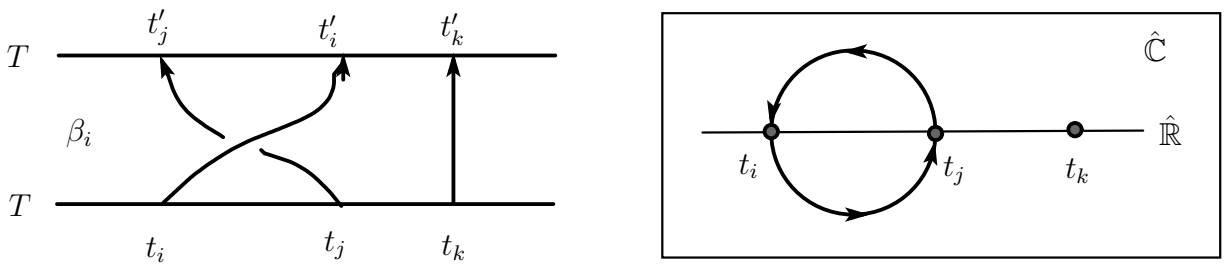

 
Let us describe the fundamental group $\pi_1(T,t)$ in terms of 
a braid group \cite{Birman}. 
We take a base point $t = (t_1,t_2,t_3) \in T$ in such a manner 
that the three points lie on the real line in an increasing order 
$t_1 < t_2 < t_3$.  
To treat them symmetrically, we denote them by $t_i$, $t_j$, $t_k$, 
where $(i,j,k)$ is a cyclic permutation of $(1,2,3)$, and think of 
them as cyclically ordered three points on the equator 
$\hat{\bR} = \bR \cup \{\infty\}$ of the Riemann sphere 
$\hat{\C} = \C \cup \{\infty\}$. 
Let $\b_i$ be a braid on three strings as in 
Figure \ref{fig:braid} (left) along which $t_i$ and $t_j$ make a 
half-turn, with $t_i$ moving in the southern hemisphere and 
$t_j$ in the northern hemisphere, while $t_k$ is kept fixed as in 
Figure \ref{fig:braid} (right). 
Then the braid group on three strings is the group generated by 
$\b_i$, $\b_j$, $\b_k$, and the pure braid group $P_3$ is the 
normal subgroup of $B_3$ generated by their squares 
$\b_i^2$, $\b_j^2$, $\b_k^2$, 
\[
P_3 = \la \b_i^2, \b_j^2, \b_k^2 \ra \triangleleft 
B_3 = \la \b_i, \b_j, \b_k \ra. 
\]
The generators of $B_3$ satisfy relations 
$\b_i\b_j\b_i = \b_j\b_i\b_j$ and $\b_k = \b_i \b_j \b_i^{-1}$, 
so that $B_3$ is generated by $\b_i$ and $\b_j$ only. 
The fundamental group $\pi_1(T,t)$ can be identified with 
the pure braid group $P_3$. 
The reduction map (\ref{eqn:reduction}) induces a group 
homomorphism $r_* : P_3 = \pi_1(T,t) \to \pi_1(Z,z)$. 
It is easy to see that this homomorphism sends 
the three basic pure braids in $P_3$ to the three basic loops in 
$\pi_1(Z,z)$ (see Figure \ref{fig:elementary}) 
in such a manner that 
\begin{equation} \label{eqn:r*} 
r_* : \b_i^2 \mapsto \ga_i \qquad (i = 1,2,3).  
\end{equation} \par 
It is sometimes convenient to lift the Poincar\'e section 
(\ref{eqn:PSt}), which makes sense for pure braids, to the 
``half-Poincar\'e section" for ordinary braids. 
Now let us construct this lift. 
The symmetric group $S_3$ acts on $T$ by permuting the entries 
of $t = (t_1,t_2,t_3)$ and the quotient space 
$T/S_3$ is the configuration space of mutually distinct, 
unordered, three points in $\C$. 
The fundamental group $\pi(T/S_3, s)$ with base point 
$s = \{t_1,t_2,t_3\}$ is identified with the ordinary braid group 
$B_3$ and there exists a short exact sequence of groups 
\[
\begin{CD}
1 @>>> \pi_1(T,t) @>>> \pi_1(T/S_3,s) @>>> S_3 @>>> 1 \\
@.       @|                    @|           @|    @.  \\
1 @>>> P_3        @>>> B_3              @>>> S_3 @>>> 1.
\end{CD}
\]
Then the Poincar\'e section (\ref{eqn:PSt}) naturally lifts to 
a collection of isomorphisms 
\[
\b_* : \M_t(\k) \to \M_{\tau(t)}(\tau(\k)), \qquad (\b \in B_3) 
\]
which should be called the {\sl half-Poincar\'e section} of 
$\PVI(\k)$, where $\tau \in S_3$ denotes the permutation 
corresponding to $\b \in B_3$. 
Note that $\tau \in S_3$ acts on $\k \in \K$ by permuting 
the entries of $(\k_1,\k_2,\k_3)$ in the same manner as it 
does on $t = (t_1,t_2,t_3) \in T$, since 
$\k_i$ is loaded on $t_i$. 
Now the permutation corresponding to the basic braid $\b_i$ is 
the substitution $\tau_i = (i,j)$ that exchanges $t_i$ and 
$t_j$ while keeping $t_k$ fixed. 
Thus there are three basic half-Poincar\'e maps: 
\begin{equation} \label{eqn:PSi}
\b_{i*} : \M_t(\k) \to \M_{\tau_i(t)}(\tau_i(\k)), 
\qquad (i = 1,2,3).  
\end{equation}
\section{Riemann-Hilbert Correspondence} \label{sec:RHC} 
It is very difficult or rather hopeless to deal with the 
Painlev\'e flow directly, since it is a highly transcendental 
dynamical system on the moduli space of stable parabolic 
connections. 
A good idea is to recast it to a more tractable dynamical 
system, called an isomonodromic flow, on a moduli space of 
monodromy representations via a Riemann-Hilbert correspondence. 
We review the construction of such a Riemann-Hilbert 
correspondence in the sequel. 
\par
\begin{figure}[t]
\begin{center}
\unitlength 0.1in
\begin{picture}(27.93,16.78)(2.10,-18.78)
%
\special{pn 20}%
\special{ar 483 480 273 273  3.1165979 4.6873942}%
%
\special{pn 20}%
\special{pa 210 480}%
\special{pa 210 1600}%
\special{fp}%
%
\special{pn 20}%
\special{ar 483 1600 273 273  1.5458015 3.1415927}%
%
\special{pn 20}%
\special{ar 1603 1040 287 287  0.0000000 6.2831853}%
%
\special{pn 20}%
\special{ar 2450 1047 287 287  0.0000000 6.2831853}%
%
\special{pn 20}%
\special{pa 490 200}%
\special{pa 2730 200}%
\special{fp}%
%
\special{pn 20}%
\special{pa 483 1873}%
\special{pa 2730 1873}%
\special{fp}%
%
\special{pn 20}%
\special{ar 2730 1600 273 273  6.2831853 6.2831853}%
\special{ar 2730 1600 273 273  0.0000000 1.5964317}%
%
\special{pn 20}%
\special{pa 763 200}%
\special{pa 763 760}%
\special{fp}%
%
\special{pn 20}%
\special{pa 1610 200}%
\special{pa 1610 753}%
\special{fp}%
%
\special{pn 20}%
\special{pa 2443 200}%
\special{pa 2443 753}%
\special{fp}%
%
\special{pn 20}%
\special{pa 3003 487}%
\special{pa 3003 1607}%
\special{fp}%
%
\special{pn 20}%
\special{ar 2723 480 274 274  4.7123890 6.2831853}%
\special{ar 2723 480 274 274  0.0000000 0.0256354}%
%
\special{pn 20}%
\special{ar 763 1054 287 287  0.0000000 6.2831853}%
%
\special{pn 13}%
\special{sh 0.600}%
\special{ar 2450 1047 41 41  0.0000000 6.2831853}%
%
\special{pn 20}%
\special{sh 0.600}%
\special{ar 1610 1040 41 41  0.0000000 6.2831853}%
%
\special{pn 20}%
\special{sh 0.600}%
\special{ar 763 1047 41 41  0.0000000 6.2831853}%
%
\special{pn 13}%
\special{pa 707 340}%
\special{pa 707 620}%
\special{fp}%
\special{sh 1}%
\special{pa 707 620}%
\special{pa 727 553}%
\special{pa 707 567}%
\special{pa 687 553}%
\special{pa 707 620}%
\special{fp}%
%
\special{pn 13}%
\special{pa 1554 333}%
\special{pa 1554 613}%
\special{fp}%
\special{sh 1}%
\special{pa 1554 613}%
\special{pa 1574 546}%
\special{pa 1554 560}%
\special{pa 1534 546}%
\special{pa 1554 613}%
\special{fp}%
%
\special{pn 13}%
\special{pa 2380 333}%
\special{pa 2380 613}%
\special{fp}%
\special{sh 1}%
\special{pa 2380 613}%
\special{pa 2400 546}%
\special{pa 2380 560}%
\special{pa 2360 546}%
\special{pa 2380 613}%
\special{fp}%
%
\special{pn 13}%
\special{pa 826 613}%
\special{pa 826 340}%
\special{fp}%
\special{sh 1}%
\special{pa 826 340}%
\special{pa 806 407}%
\special{pa 826 393}%
\special{pa 846 407}%
\special{pa 826 340}%
\special{fp}%
%
\special{pn 13}%
\special{pa 1673 606}%
\special{pa 1673 333}%
\special{fp}%
\special{sh 1}%
\special{pa 1673 333}%
\special{pa 1653 400}%
\special{pa 1673 386}%
\special{pa 1693 400}%
\special{pa 1673 333}%
\special{fp}%
%
\special{pn 13}%
\special{pa 2506 606}%
\special{pa 2506 333}%
\special{fp}%
\special{sh 1}%
\special{pa 2506 333}%
\special{pa 2486 400}%
\special{pa 2506 386}%
\special{pa 2526 400}%
\special{pa 2506 333}%
\special{fp}%
%
\special{pn 13}%
\special{pa 728 1341}%
\special{pa 791 1341}%
\special{fp}%
\special{sh 1}%
\special{pa 791 1341}%
\special{pa 724 1321}%
\special{pa 738 1341}%
\special{pa 724 1361}%
\special{pa 791 1341}%
\special{fp}%
%
\special{pn 13}%
\special{pa 1568 1327}%
\special{pa 1638 1320}%
\special{fp}%
\special{sh 1}%
\special{pa 1638 1320}%
\special{pa 1570 1307}%
\special{pa 1585 1325}%
\special{pa 1574 1347}%
\special{pa 1638 1320}%
\special{fp}%
%
\special{pn 13}%
\special{pa 2422 1334}%
\special{pa 2478 1341}%
\special{fp}%
\special{sh 1}%
\special{pa 2478 1341}%
\special{pa 2414 1313}%
\special{pa 2425 1334}%
\special{pa 2409 1353}%
\special{pa 2478 1341}%
\special{fp}%
%
\special{pn 13}%
\special{pa 1673 1873}%
\special{pa 1610 1873}%
\special{fp}%
\special{sh 1}%
\special{pa 1610 1873}%
\special{pa 1677 1893}%
\special{pa 1663 1873}%
\special{pa 1677 1853}%
\special{pa 1610 1873}%
\special{fp}%
\put(7.6000,-14.8000){\makebox(0,0){$C_1$}}%
\put(16.3800,-14.8100){\makebox(0,0){$C_2$}}%
\put(24.7100,-14.8800){\makebox(0,0){$C_3$}}%
\put(16.4500,-17.4000){\makebox(0,0){$C_4$}}%
\put(7.6300,-11.8000){\makebox(0,0){$t_1$}}%
\put(16.1700,-11.8700){\makebox(0,0){$t_2$}}%
\put(24.7100,-11.8000){\makebox(0,0){$t_3$}}%
\end{picture}%
\end{center}
\caption{Four loops in $\P^1-D_t$; 
the fourth point $t_4$ is outside $C_4$, invisible.} 
\label{fig:loop2}
\end{figure}
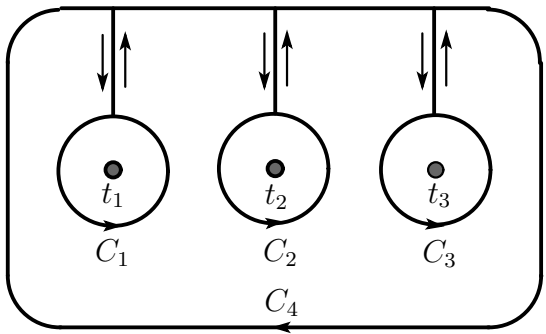

 
Let $A := \C^4$ be the complex $4$-space with 
coordinates $a = (a_1,a_2,a_3,a_4)$, called the space of 
local monodromy data. 
Given $(t,a) \in T \times A$, let $\R_t(a)$ be the moduli 
space of Jordan equivalence classes of representations 
$\rho : \pi_1(\P^1-D_t,*) \to SL_2(\C)$
such that $\Tr\, \rho(C_i) = a_i$ for $i \in \{1,2,3,4\}$, 
where the divisor $D_t = t_1+t_2+t_3+t_4$ is identified with 
the point set $\{t_1,t_2,t_3,t_4\}$ and $C_i$ is a loop 
surrounding $t_i$ as in Figure \ref{fig:loop2}. 
Any stable parabolic connection $Q = (E,\nabla,\psi,l) 
\in \M_t(\k)$, when restricted to $\P^1-D_t$, induces a flat 
connection 
\[
\nabla|_{\P^1-D_t} : E|_{\P^1-D_t} \to (E|_{\P^1-D_t}) 
\otimes \Omega_{\P^1-D_t}^1,
\]
and one can speak of the Jordan equivalence class $\rho$ of 
its monodromy representations. 
Then the Riemann-Hilbert correspondence at $t \in T$ is 
defined by  
\begin{equation} \label{eqn:RHtk}
\RH_{t,\k} : \M_t(\k) \to \R_t(a), \quad Q \mapsto \rho, 
\end{equation} 
where in view of the Riemann scheme in Table \ref{tab:riemann}, 
the local monodromy data $a \in A$ is given by 
\begin{equation} \label{eqn:a}
a_i = \left\{
\begin{array}{ll}
{\-}2 \cos \pi \k_i \qquad & (i = 1,2,3), \\[2mm] 
-2 \cos \pi \k_4 \qquad & (i = 4). 
\end{array}
\right.
\end{equation}
\par 
As a relative setting over $T$, let $\pi_a : \R(a) \to T$ be 
the family of moduli spaces of monodromy representations with 
fiber $\R_t(a)$ over $t \in T$. 
Then the relative version of Riemann-Hilbert correspondence 
is formulated as the commutative diagram 
\begin{equation} \label{eqn:RHk}
\begin{CD} 
\M(\k) @> \RH_{\k} >> \R(a)    \\
@V \pi_{\k} VV     @VV \pi_a V \\
T @=  T, 
\end{CD}
\end{equation}
whose fiber over $t \in T$ is given by (\ref{eqn:RHtk}). 
Then we have the following theorem \cite{IIS1,IIS2}. 
\begin{theorem} \label{thm:RH} 
If $\k \in \K-\Wall$, then $\R(a)$ as well as each fiber 
$\R_t(a)$ is smooth and the Riemann-Hilbert correspondence 
$\RH_{\k}$ in $(\ref{eqn:RHk})$ is a biholomorphism. 
\end{theorem}
\begin{remark} \label{rem:RH} 
If $\k \in \Wall$, then $\R_t(a)$ is not a smooth surface 
but a surface with Klein singularities and (\ref{eqn:RHtk}) 
yields an analytic minimal resolution of singularities, 
so that (\ref{eqn:RHk}) gives a family of resolutions of 
singularities. 
We refer to \cite{IIS1} for a detailed description of these 
singularity structures. 
As is mentioned in Remark \ref{rem:generic}, this fact makes 
the treatment of the nongeneric case more involved and 
we leave this case in another occasion. 
\end{remark}
\section{Cubic Surface and the 27 Lines} \label{sec:cubic}
In this section, following the construction in \cite{IIS1}, 
we shall realize the moduli space $\R_t(a)$ of monodromy 
representations as an affine cubic surface $\Sol(\th)$ and 
describe the braid group action on $\R_t(a)$ explicitly in 
terms of $\Sol(\th)$. 
Moreover we discuss some materials from the geometry of a 
cubic surface, including the 27 lines on it, as a preliminary 
to the later sections. 
\par 
Given $\th = (\th_1,\th_2,\th_3,\th_4) \in \Th := \C^4_{\th}$, 
we consider an affine cubic surface 
\[
\Sol(\th) = \{\, x = (x_1,x_2,x_3) \in \C^3_x \,:\, 
f(x,\th) = 0 \, \}, 
\]
where the cubic polynomial $f(x,\th)$ of $x$ with parameter 
$\th$ is given by 
\[
f(x,\th) = x_1x_2x_3 + x_1^2 + x_2^2 + x_3^2 
- \th_1 x_1 - \th_2 x_2 - \th_3 x_3 + \th_4. 
\]
Then there exists an isomorphism of affine algebraic 
surfaces, $\R_t(a) \to \Sol(\th)$, $\rho \mapsto x$, 
where 
\[
x_i = \Tr\,\rho(C_jC_k), \qquad 
\text{for} \quad \{i,j,k\} = \{1,2,3\}, 
\]
together with a correspondence of parameters, 
$A \to \Th$, $a \mapsto \th$, given by 
\begin{equation} \label{eqn:th}
\th_i = \left\{
\begin{array}{ll}
a_i a_4 + a_j a_k \qquad & (\{i,j,k\}=\{1,2,3\}), \\[2mm]
a_1 a_2 a_3 a_4 + a_1^2 + a_2^2 + a_3^2 + a_4^2 - 4 
\qquad & (i = 4). 
\end{array}
\right.
\end{equation}
\par 
With this identification, the Riemann-Hilbert correspondence 
(\ref{eqn:RHtk}) is reformulated as a map
\begin{equation} \label{eqn:RHtk2}
\RH_t(\k) : \M_t(\k) \to \Sol(\th), 
\qquad \text{with} \quad \th = \rh(\k), 
\end{equation}
where $\rh : \K \to \Th$ is the composition of the maps 
$\K \to A$ and $A \to \Th$ defined by (\ref{eqn:a}) and 
(\ref{eqn:th}), and is referred to as the Riemann-Hilbert 
correspondence in the parameter level. 
Through the reformulated Riemann-Hilbert correspondence 
(\ref{eqn:RHtk2}), the $i$-th basic half-Poincar\'e map 
$\b_{i*}$ in (\ref{eqn:PSi}) is conjugated to a map 
$g_i : \Sol(\th) \to \Sol(\th')$, 
$(x,\th) \mapsto (x',\th')$, defined by 
\begin{equation} \label{eqn:gi}
g_i \quad : \quad 
(x_i',x_j',x_k',\th_i',\th_j',\th_k',\th_4') = 
(\th_j-x_j-x_kx_i, x_i, x_k, \th_j, \th_i, \th_k, \th_4), 
\end{equation}
where $(i,j,k)$ is a cyclic permutation of $(1,2,3)$. 
A derivation of this formula can be found in \cite{Iwasaki3} 
(see also \cite{Boalch2,DM,Goldman,Iwasaki2,Jimbo}). 
The map (\ref{eqn:gi}) is strictly conjugate to the map 
(\ref{eqn:PSi}), since (\ref{eqn:RHtk2}) is biholomorphic 
by Theorem \ref{thm:RH}. 
We can easily check the relations $g_ig_jg_i = g_jg_ig_j$ and 
$g_k = g_i g_j g_i^{-1}$, which are just parallel to those 
for the braids $\b_i$, $\b_j$, $\b_k$. 
\begin{remark} \label{rem:area2} 
The affine cubic surface $\Sol(\th)$ admits a natural complex 
area form 
\begin{equation} \label{eqn:poincare} 
\omega(\th) = \dfrac{dx_1 \wedge dx_2 \wedge dx_3}{d_xf(x,\th)}, 
\end{equation} 
the Poincar\'e residue for the surface $\Sol(\th)$. 
The transformations $g_i$ are area-preserving with 
respect to $\omega(\th)$. 
It is known \cite{IIS1,Iwasaki1,Iwasaki3} that the standard 
area form $\omega_t(\k)$ on the moduli space $\M_t(\k)$ in 
Remark \ref{rem:area} is the pull-back of $\omega(\th)$ 
by the Riemann-Hilbert correspondence (\ref{eqn:RHtk2}). 
\end{remark} 
\par
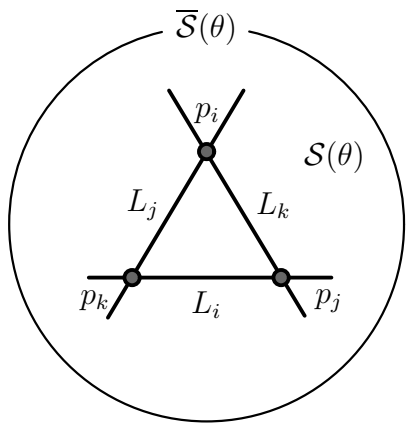
\begin{figure}[t] 
\begin{center}
\unitlength 0.1in
\begin{picture}(20.68,21.44)(14.16,-27.54)
%
\special{pn 13}%
\special{ar 2450 1720 1034 1034  4.9299431 6.2831853}%
\special{ar 2450 1720 1034 1034  0.0000000 4.4841152}%
%
\special{pn 20}%
\special{pa 1834 2000}%
\special{pa 3104 2000}%
\special{fp}%
%
\special{pn 20}%
\special{pa 2644 1020}%
\special{pa 1934 2210}%
\special{fp}%
%
\special{pn 20}%
\special{pa 2254 1020}%
\special{pa 2964 2200}%
\special{fp}%
\put(22.9000,-7.8000){\makebox(0,0)[lb]{$\ol{\Sol}(\th)$}}%
\put(29.6000,-14.5000){\makebox(0,0)[lb]{$\Sol(\th)$}}%
\put(23.7000,-22.1000){\makebox(0,0)[lb]{$L_i$}}%
\put(20.3000,-17.0000){\makebox(0,0)[lb]{$L_j$}}%
\put(27.1000,-16.9000){\makebox(0,0)[lb]{$L_k$}}%
\put(30.2000,-21.8000){\makebox(0,0)[lb]{$p_j$}}%
\put(23.9000,-11.9000){\makebox(0,0)[lb]{$p_i$}}%
\put(17.9000,-21.7000){\makebox(0,0)[lb]{$p_k$}}%
%
\special{pn 20}%
\special{sh 0.600}%
\special{ar 2060 2000 40 40  0.0000000 6.2831853}%
%
\special{pn 20}%
\special{sh 0.600}%
\special{ar 2450 1340 40 40  0.0000000 6.2831853}%
%
\special{pn 20}%
\special{sh 0.600}%
\special{ar 2840 2000 40 40  0.0000000 6.2831853}%
\end{picture}%
\end{center}
\caption{Tritangent lines at infinity on $\ol{\Sol}(\th)$} 
\label{fig:cubic1} 
\end{figure}

 
In order to utilize standard techniques from algebraic geometry 
and complex geometry, we need to compactify the affine cubic 
surface $\Sol(\th)$ by a standard embedding 
\[
\Sol(\th) \hookrightarrow \ol{\Sol}(\th) \subset \P^3, \qquad 
x = (x_1,x_2,x_3) \mapsto [1:x_1:x_2:x_3], 
\]
where the compactified surface $\ol{\Sol}(\th)$ is defined 
by $\ol{\Sol}(\th) = \{\,X \in \P^3\,:\, F(X,\th)=0 \,\}$ with 
\[
F(X,\th) = X_1X_2X_3+X_0(X_1^2+X_2^2+X_3^2)-
X_0^2(\th_1X_1+\th_2X_2+\th_3X_3)+\th_4X_0^3. 
\]
It is obtained from the affine surface $\Sol(\th)$ by adding 
three lines at infinity, 
\begin{equation} \label{eqn:lines}
L_i = \{\, X \in \P^3 \,:\, X_0 = X_i = 0 \,\} \qquad 
(i = 1,2,3). 
\end{equation}
The union $L = L_1 \cup L_2 \cup L_3$ is called the 
tritangent lines at infinity and the intersection point of 
$L_j$ and $L_k$ is denoted by $p_i$ as in Figure \ref{fig:cubic1}. 
Note that 
\[
p_1 = [0:1:0:0], \qquad p_2 = [0:0:1:0], \qquad 
p_3 = [0:0:0:1]. 
\]
For $i \in \{1,2,3\}$ we put 
$U_i = \{\, X \in \P^3 \,:\, X_i \neq 0 \,\}$ and take 
inhomogeneous coordinates of $\P^3$ as 
\begin{equation} \label{eqn:coord}
\begin{array}{rclcll}
u &=& (u_0,u_j,u_k) &=& (X_0/X_i, X_j/X_i, X_k/X_i) \qquad & 
\text{on} \quad U_i, \\[1mm]
v &=& (v_0,v_i,v_k) &=& (X_0/X_j, X_i/X_j, X_k/X_j) \qquad & 
\text{on} \quad U_j, \\[1mm]
w &=& (w_0,w_i,w_j) &=& (X_0/X_k, X_i/X_k, X_j/X_k) \qquad & 
\text{on} \quad U_k, 
\end{array}
\end{equation}
where $\{i,j,k\} = \{1,2,3\}$. 
In terms of these coordinates we shall find local coordinates 
and local equations of $\ol{\Sol}(\th)$ around $L$. 
Since $L \subset U_1 \cup U_2 \cup U_3$, we can divide $L$ 
into three components $L \cap U_i$, $i = 1,2,3$, and make a 
further decomposition 
\[
L \cap U_i = \{p_i\} \cup (L_j-\{p_i,p_k\}) \cup (L_k-\{p_i,p_j\}) 
\qquad (\{i,j,k\} = \{1,2,3\}) 
\]
into a total of nine pieces. 
Then a careful inspection of equation $F(X,\th) = 0$ implies 
that around those pieces we can take local coordinates and 
local equations as in Table \ref{tab:local}, where 
$O_m(u_j, u_k) = O((|u_j|+|u_k|)^m)$ denotes a small term of 
order $m$ as $(u_j, u_k) \to (0,0)$. 
\begin{table}[t] 
\begin{center}
\begin{tabular}{|c|c|l|}
\hline
\vspace{-3mm} & & \\
coordinates & valid around & local equation \\[2mm]
\hline
\hline
\vspace{-3mm} & & \\
$(u_j,u_k)$ & $p_i$ & 
$u_0 = -(u_j u_k) \{1-(u_j^2+\th_i u_j u_k+u_k^2)+O_3(u_j,u_k)\}$ \\[2mm]
\hline
\vspace{-3mm} & & \\
$(u_0,u_k)$ & $L_j-\{ p_i,p_k \}$ & 
$u_j = -(u_k+1 / u_k)u_0+(\th_k+\th_i / u_k) u_0^2 + O(u_0^3)$ \\[2mm]
\hline
\vspace{-3mm} & & \\
$(u_0,u_j)$ & $L_k-\{ p_i,p_j \} $ & 
$u_k = -(u_j+1 / u_j)u_0+(\th_j+\th_i / u_j) u_0^2 + O(u_0^3)$ \\[2mm]
\hline
\hline
\vspace{-3mm} & & \\
$(v_i,v_k)$ & $p_j$ & 
$v_0 = -(v_i v_k) \{1-(v_i^2+\th_j v_i v_k+v_k^2)+O_3(v_i,v_k)\}$ \\[2mm]
\hline
\vspace{-3mm} & & \\
$(v_0,v_i)$ & $L_k - \{ p_i,p_j \}$ & 
$v_k = -(v_i+1 / v_i)v_0+(\th_i+\th_j / v_i) v_0^2 + O(v_0^3)$ \\[2mm]
\hline
\vspace{-3mm} & & \\
$(v_0,v_k)$ & $L_i - \{ p_j,p_k \}$ & 
$v_i = -(v_k+1 / v_k)v_0+(\th_k+\th_j / v_k) v_0^2 + O(v_0^3)$ \\[2mm]
\hline
\hline
\vspace{-3mm} & & \\
$(w_i,w_j)$ & $p_k$ & 
$w_0 = -(w_i w_j) \{1-(w_i^2+\th_k w_i w_j+w_j^2)+O_3(w_i,w_j)\}$ \\[2mm]
\hline
\vspace{-3mm} & & \\
$(w_0,w_j)$ & $L_i -\{ p_j,p_k\}$ & 
$w_i = -(w_j+1 / w_j)w_0+(\th_j+\th_k / w_j) w_0^2 + O(w_0^3)$ \\[2mm]
\hline
\vspace{-3mm} & & \\
$(w_0,w_i)$ & $L_j - \{ p_i,p_k \}$ & 
$w_j = -(w_i+1 / w_i)w_0+(\th_i+\th_k / w_i) w_0^2 + O(w_0^3)$ \\[2mm]
\hline
\end{tabular}
\end{center}
\caption{Local coordinates and local equations of $\ol{\Sol}(\th)$}
\label{tab:local}
\end{table}
\begin{lemma} \label{lem:smooth} 
As to the smoothness of the surface $\ol{\Sol}(\th)$, 
the following hold. 
\begin{enumerate} 
\item For any $\th \in \Th$, the surface $\ol{\Sol}(\th)$ is 
smooth in a neighborhood of $L$.  
\item If $\th = \rh(\k)$ with $\k \in \K$, the surface 
$\ol{\Sol}(\th)$ is smooth everywhere if and only if 
$\k \in \K-\Wall$. 
\end{enumerate} 
\end{lemma}
{\it Proof}. 
In terms of the inhomogeneous coordinates $u$ in (\ref{eqn:coord}), 
we have 
\[
\ol{\Sol}(\th) \cap U_i \cong 
\{\, u = (u_0,u_j,u_k) \in \C^3 \,:\, f_i(u,\th) = 0 \,\},
\]
where the defining equation $f_i(u,\th)$ is given by 
\[
f_i(u,\th) = u_j u_k +u_0 (1 + u_j^2 + u_k^2) - u_0^2 
(\th_i + \th_j u_j + \th_k u_k) + \th_4 u_0^3. 
\]
The partial derivatives of $f_i = f_i(u,\th)$ with respect 
to $u = (u_0, u_j, u_k)$ are calculated as 
\begin{eqnarray*}
\frac{\partial f_i}{\partial u_0} 
&=& (1 + u_j^2 + u_k^2) 
- 2 u_0 (\th_i + \th_j u_j + \th_k u_k) + 3 \th_4 u_0^2 \\[2mm]
\frac{\partial f_i}{\partial u_j} 
&=& u_k + 2 u_0 u_j - \th_j u_0^2 \\[2mm]
\frac{\partial f_i}{\partial u_k} 
&=& u_j + 2 u_0 u_k - \th_k u_0^2 .
\end{eqnarray*}
Restricted to the set 
$L \cap U_i = (L_j \cap U_i) \cup (L_k \cap U_i)$, 
these derivatives become 
\[
\begin{array}{rclrclrcll}
\dfrac{\partial f_i}{\partial u_0} &=& 1 + u_k^2, \qquad & 
\dfrac{\partial f_i}{\partial u_j} &=& u_k,       \qquad & 
\dfrac{\partial f_i}{\partial u_k} &=& 0,         \qquad & 
\text{on} \quad L_j \cap U_i, \\[4mm]
\dfrac{\partial f_i}{\partial u_0} &=& 1 + u_j^2, \qquad & 
\dfrac{\partial f_i}{\partial u_j} &=& 0,         \qquad & 
\dfrac{\partial f_i}{\partial u_k} &=& u_j,       \qquad & 
\text{on} \quad L_k \cap U_i. 
\end{array}
\]
Hence the exterior derivative $d_u f_i$ does not vanish on 
$L \cap U_i$, and the implicit function theorem implies that 
$\ol{\Sol}(\th)$ is smooth in a neighborhood of $L$. 
This proves assertion (1). 
In order to show assertion (2) we recall that the affine 
surface $\Sol(\th)$ is smooth if and only if $\th = \rh(\k)$ 
with $\k \in \K-\Wall$ (see \cite{IIS1}). 
Then assertion (2) readily follows from assertion (1). 
\hfill $\Box$ 
\par\medskip 
Now let us review some basic facts about smooth cubic surfaces 
in $\P^3$ (see e.g. \cite{GH}).
It is well known that every smooth cubic surface $S$ 
in $\P^3$ can be obtained by blowing up $\P^2$ at six points 
$P_1, \dots, P_6$, no three colinear and not all six on a conic, 
and embedding the blow-up surface into $\P^3$ by the proper 
transform of the linear system of cubics passing through the 
six points $P_1, \dots, P_6$. 
It is also well known that there are exactly $27$ 
lines on the smooth cubic surface $S$, each of which has 
self-intersection number $-1$. 
Explicitly, they are given by
\[
E_a \quad (a = 1,\dots,6); \qquad F_{ab} \quad 
(1 \le a < b \le 6); \qquad G_a \quad (a = 1,\dots,6), 
\]
\begin{enumerate}
\item $E_a$ is the exceptional curve over the point $P_a$, 
\item $F_{ab}$ is the strict transform of the line in 
$\P^2$ through the two points $P_a$ and $P_b$, 
\item $G_a$ is the strict transform of the conic in $\P^2$ 
through the five points $P_1,\dots,\hat{P}_a,\dots,P_6$. 
\end{enumerate}
Here the index $a$ should not be confused with the 
local monodromy data $a \in A$.  
All the intersection relations among the $27$ lines 
with {\sl nonzero} intersection numbers are listed as 
\[
\begin{array}{rl}
(E_a, E_a) = (G_a, G_a) = (F_{ab}, F_{ab}) = -1 \qquad & 
(\forall \, a,b), \\[2mm]
(E_a, F_{bc}) = (G_a, F_{bc}) = {\-} 1  \qquad & 
(a \in \{b,c\}), \\[2mm] 
(E_a, G_b) = {\-} 1 \qquad & 
(a \neq b),  \\[2mm]
(F_{ab}, F_{cd}) = {\-} 1 \qquad & 
(\{a,b\} \cap \{c,d\} = \emptyset).
\end{array}
\]
\par 
\begin{figure}[t] 
\begin{center}
\unitlength 0.1in
\begin{picture}(48.10,35.30)(11.10,-35.70)
\put(59.2000,-31.7000){\makebox(0,0)[lb]{$L_1=F_{14}$}}%
\put(11.8000,-37.4000){\makebox(0,0)[lb]{$L_2=F_{25}$}}%
\put(29.8000,-2.1000){\makebox(0,0)[lb]{$L_3=F_{36}$}}%
%
\special{pn 20}%
\special{pa 1170 3080}%
\special{pa 5790 3080}%
\special{fp}%
%
\special{pn 20}%
\special{pa 3790 270}%
\special{pa 1380 3510}%
\special{fp}%
\special{pa 3180 270}%
\special{pa 5570 3470}%
\special{fp}%
%
\special{pn 13}%
\special{pa 2090 2920}%
\special{pa 2470 3520}%
\special{fp}%
\special{pa 2480 2920}%
\special{pa 2090 3510}%
\special{fp}%
\special{pa 2890 2920}%
\special{pa 3280 3510}%
\special{fp}%
\special{pa 3280 2920}%
\special{pa 2890 3510}%
\special{fp}%
\special{pa 3680 2910}%
\special{pa 4080 3520}%
\special{fp}%
\special{pa 4090 2920}%
\special{pa 3680 3520}%
\special{fp}%
\special{pa 4480 2930}%
\special{pa 4870 3520}%
\special{fp}%
\special{pa 4880 2930}%
\special{pa 4490 3520}%
\special{fp}%
%
\special{pn 13}%
\special{pa 3500 840}%
\special{pa 4187 768}%
\special{fp}%
\special{pa 3726 1150}%
\special{pa 3959 473}%
\special{fp}%
\special{pa 3963 1475}%
\special{pa 4648 1418}%
\special{fp}%
\special{pa 4189 1784}%
\special{pa 4422 1108}%
\special{fp}%
\special{pa 4412 2108}%
\special{pa 5118 2046}%
\special{fp}%
\special{pa 4657 2426}%
\special{pa 4887 1728}%
\special{fp}%
\special{pa 4891 2731}%
\special{pa 5575 2672}%
\special{fp}%
\special{pa 5122 3047}%
\special{pa 5355 2371}%
\special{fp}%
%
\special{pn 13}%
\special{pa 2970 520}%
\special{pa 3242 1155}%
\special{fp}%
\special{pa 2741 827}%
\special{pa 3457 849}%
\special{fp}%
\special{pa 2501 1150}%
\special{pa 2758 1787}%
\special{fp}%
\special{pa 2272 1457}%
\special{pa 2987 1479}%
\special{fp}%
\special{pa 2029 1766}%
\special{pa 2297 2422}%
\special{fp}%
\special{pa 1798 2094}%
\special{pa 2533 2108}%
\special{fp}%
\special{pa 1576 2407}%
\special{pa 1835 3044}%
\special{fp}%
\special{pa 1342 2722}%
\special{pa 2057 2745}%
\special{fp}%
\put(56.2000,-27.5000){\makebox(0,0)[lb]{$E_3$}}%
\put(54.0000,-23.9000){\makebox(0,0)[lb]{$G_6$}}%
\put(51.6000,-21.1000){\makebox(0,0)[lb]{$E_6$}}%
\put(49.2000,-17.2000){\makebox(0,0)[lb]{$G_3$}}%
\put(47.0000,-14.8000){\makebox(0,0)[lb]{$F_{12}$}}%
\put(44.4000,-11.4000){\makebox(0,0)[lb]{$F_{45}$}}%
\put(42.2000,-8.2000){\makebox(0,0)[lb]{$F_{15}$}}%
\put(39.9000,-5.2000){\makebox(0,0)[lb]{$F_{24}$}}%
\put(28.2000,-5.0000){\makebox(0,0)[lb]{$E_2$}}%
\put(25.4000,-9.1000){\makebox(0,0)[lb]{$G_5$}}%
\put(23.8000,-11.3000){\makebox(0,0)[lb]{$E_5$}}%
\put(20.6000,-15.3000){\makebox(0,0)[lb]{$G_2$}}%
\put(18.5000,-17.5000){\makebox(0,0)[lb]{$F_{13}$}}%
\put(15.4000,-21.5000){\makebox(0,0)[lb]{$F_{46}$}}%
\put(14.6000,-23.8000){\makebox(0,0)[lb]{$F_{34}$}}%
\put(11.1000,-28.0000){\makebox(0,0)[lb]{$F_{16}$}}%
\put(19.9000,-37.1000){\makebox(0,0)[lb]{$E_1$}}%
\put(23.9000,-37.3000){\makebox(0,0)[lb]{$G_4$}}%
\put(28.2000,-37.3000){\makebox(0,0)[lb]{$E_4$}}%
\put(32.1000,-37.2000){\makebox(0,0)[lb]{$G_1$}}%
\put(36.2000,-37.1000){\makebox(0,0)[lb]{$F_{23}$}}%
\put(40.3000,-37.1000){\makebox(0,0)[lb]{$F_{56}$}}%
\put(44.4000,-37.1000){\makebox(0,0)[lb]{$F_{26}$}}%
\put(48.4000,-37.2000){\makebox(0,0)[lb]{$F_{35}$}}%
\end{picture}%
\end{center}
\caption{The $27$ lines on $\ol{\Sol}(\th)$ 
viewed from the tritangent plane at infinity} 
\label{fig:cubic3} 
\end{figure}
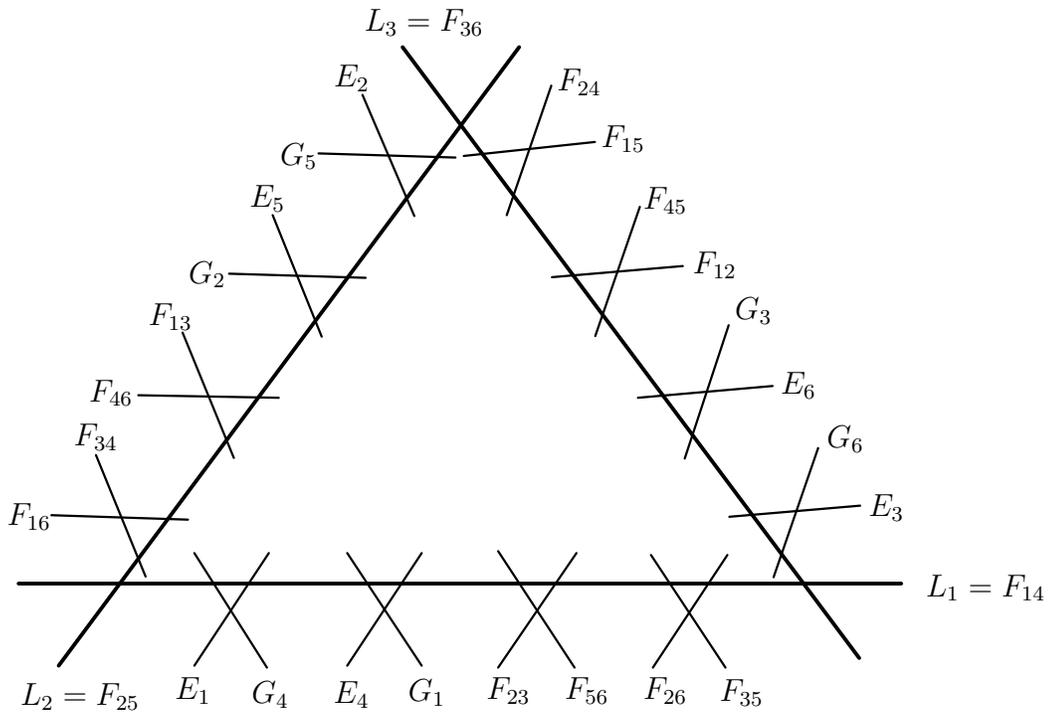
Moreover there are exactly $45$ tritangent planes that 
cut out a triplet of lines on $S$. 
In our case $S = \ol{\Sol}(\th)$, the plane at infinity 
$\{\, X \in \P^3 \,:\, X_0 = 0 \,\}$ is an instance of 
tritangent plane, which cuts out the lines in 
(\ref{eqn:lines}). 
The arrangement of the $27$ lines viewed from the tritangent 
plane at infinity is shown in Figure \ref{fig:cubic3} and 
the lines at infinity are given by 
\begin{equation} \label{eqn:lineinf}
L_1 = F_{14},  \qquad L_2 = F_{25}, \qquad L_3 = F_{36}. 
\end{equation}
Each line at infinity is intersected by exactly eight lines 
and this fact enables us to divide the 27 lines into three 
groups of nine lines labeled by lines at infinity. 
{\sl Caution:} only the intersection relations among the 
lines of the same group are indicated and no other 
intersection relations are depicted in Figure \ref{fig:cubic3}. 
\par 
If $E_0$ is the strict transform of a line in $\P^2$ not 
passing through $P_1, \dots, P_6$ 
relative to the $6$-point blow-up $S \to \P^2$, then the 
second cohomology group of $S = \ol{\Sol}(\th)$ is 
expressed as 
\begin{equation} \label{eqn:basis}
H^2(\ol{\Sol}(\th),\Z) = 
\Z E_0 \oplus \Z E_1 \oplus \Z E_2 \oplus \Z E_3 \oplus 
\Z E_4 \oplus \Z E_5 \oplus \Z E_6, 
\end{equation}
where a divisor is identified with the cohomology class 
it represents. 
It is a Lorentzian lattice of rank $7$ with intersection 
numbers 
\begin{equation} \label{eqn:in1}
(E_a, E_b) = \left\{ 
\begin{array}{ll}
{\-}1 \quad & (a = b = 0), \\[1mm] 
 -1 \quad & (a = b \neq 0), \\[1mm]
{\-}0 \quad & (\text{otherwise}). 
\end{array}\right. 
\end{equation}
In terms of the basis in (\ref{eqn:basis}) 
the lines $F_{ab}$ and $G_a$ are represented as 
\begin{equation} \label{eqn:FG}
F_{ab} = E_0 - E_a - E_b, \qquad 
G_a = 2 E_0 - (E_1 + \cdots + \widehat{E}_a + \cdots + E_6). 
\end{equation}
\par
\begin{table}[t] 
\begin{center} 
\begin{tabular}{|c|c|c|}
\hline
\vspace{-4mm} & & \\
$1$ & $L_i(b_i,b_4;b_j,b_k)$ & $L_i(1/b_i,1/b_4;b_j,b_k)$ \\[1mm]
\hline
\hline
\vspace{-4mm} & & \\
$2$ & $L_i(b_j,b_k;b_i,b_4)$ & $L_i(1/b_j,1/b_k;b_i,b_4)$ \\[1mm]
\hline
\hline
\vspace{-4mm} & & \\
$3$ & $L_i(1/b_i,b_4;b_j,b_k)$ & $L_i(b_i,1/b_4;b_j,b_k)$ \\[1mm]
\hline
\hline
\vspace{-4mm} & & \\
$4$ & $L_i(1/b_j,b_k;b_i,b_4)$ & $L_i(b_j,1/b_k;b_i,b_4)$ \\[1mm]
\hline
\end{tabular}
\end{center} 
\caption{Eight lines intersecting the line $L_i$ at infinity, 
divided into four pairs} 
\label{tab:line}
\end{table}
We shall describe the $27$ lines on our cubic surface 
$\ol{\Sol}(\th)$ under the condition that $\ol{\Sol}(\th)$ 
is smooth, namely, $\th = \rh(\k)$ with $\k \in \K-\Wall$. 
To this end we introduce new parameters 
$b = (b_1,b_2,b_3,b_4) \in B := (\C_b^{\times})^4$ in such a 
manner that $b$ is expressed as 
\[ 
b_i = \left\{
\begin{array}{rl}
\exp(\sqrt{-1}\pi\k_i) \qquad  & (i = 1,2,3), \\[2mm]
-\exp(\sqrt{-1}\pi\k_4) \qquad & (i = 4), 
\end{array}
\right.
\]
as a function of $\k \in \K$. 
Then the Riemann scheme in Table \ref{tab:riemann} implies 
that $b_i$ is an eigenvalue of the monodromy matrix $\rho(C_i)$ 
around the point $t_i$ and formula (\ref{eqn:a}) implies that 
$a_i = b_i + b_i^{-1}$. 
Here parameters $b \in B$ should not be confused with the 
index $b$ above. 
In terms of the parameters $b \in B$, the discriminant $\vD(\th)$ 
of the cubic surfaces $\Sol(\th)$ factors as 
\begin{equation} \label{eqn:vD}
\vD(\th) = \prod_{l=1}^4(b_l-b_l^{-1})^2 
\prod_{\ve \in \{\pm1\}^4}(b^{\ve}-1), 
\end{equation}
where 
$b^{\ve} = b_1^{\ve_1}b_2^{\ve_2}b_3^{\ve_3}b_4^{\ve_4}$ 
for each quadruple sign 
$\ve = (\ve_1,\ve_2,\ve_3,\ve_4) \in \{\pm1\}^4$. 
Formula (\ref{eqn:vD}) clearly shows for which parameters 
$b \in B$ the cubic surface $\ol{\Sol}(\th)$ is smooth 
or singular. 
\par
Assume that $\ol{\Sol}(\th)$ is smooth, namely, $\vD(\th) \neq 0$. 
Then, as is mentioned earlier, for each index $i \in \{1,2,3\}$ 
with $\{i,j,k\} = \{1,2,3\}$, there are exactly eight lines on 
$\ol{\Sol}(\th)$ intersecting $L_i$, but not intersecting the 
remaining two lines at infinity, $L_j$ and $L_k$. 
They are just $\{E_i, G_{i+3}\}$, $\{E_{i+3}, G_i\}$, 
$\{F_{jk}, F_{j+3,k+3}\}$, 
$\{F_{j,k+3}, F_{j+3,k}\}$ as in Figure \ref{fig:cubic3}, where 
two lines from the same pair intersect, but ones from 
different pairs are disjoint. 
In terms of the parameters $b \in B$ introduced above, those 
eight lines are given as in Table \ref{tab:line}, where  
$L_i(b_i,b_4;b_j,b_k)$ stands for the line in $\P^3$ defined by 
the system of linear equations 
\begin{equation} \label{eqn:line2} 
X_i = (b_i b_4 + b_i^{-1}b_4^{-1}) X_0, \qquad 
X_j + (b_ib_4) X_k =  
\{ b_i (b_k+b_k^{-1}) + b_4 (b_j+b_j^{-1}) \} X_0. 
\end{equation}
\section{Dynamical System on Cubic Surface} \label{sec:dynamics}
The affine cubic surface $\Sol(\th)$ is a $(2,2,2)$-surface, 
that is, its defining equation $f(x,\th) = 0$ is a quadratic 
equation in each variable $x_i$, $i = 1,2,3$. 
Therefore the line through a point $x \in \Sol(\th)$ parallel 
to the $x_i$-axis passes through a unique second point 
$x' \in \Sol(\th)$ (see Figure \ref{fig:involution}). 
This defines an involution $\si_i : \Sol(\th) \to \Sol(\th)$, 
$x \mapsto x'$, which is explicitly given by 
\begin{equation} \label{eqn:si} 
\si_i : \qquad 
(x_i',x_j', x_k') = (\th_i-x_i-x_jx_k, x_j, x_k), 
\qquad (i = 1,2,3).  
\end{equation}
It is easy to see that the involution $\si_i$ preserves the 
Poincar\'e residue $\omega(\th)$ in (\ref{eqn:poincare}). 
\par 
The automorphism $\si_i$ of the affine surface 
$\Sol(\th)$ extends to a birational map of the projective surface 
$\ol{\Sol}(\th)$, which will also be denoted by $\si_i$. 
In terms of the homogeneous coordinates $X$ of $\P^3$, the 
birational map $\si_i : X \mapsto X'$ is expressed as 
\[
[X_0':X_i':X_j':X_k'] = 
[X_0^2 : \th_i X_0^2 - X_0 X_i - X_j X_k : X_0 X_j : X_0 X_k]. 
\]
Let $G = \la \si_1, \si_2, \si_3 \ra$ be the group of birational 
transformations on $\ol{\Sol}(\th)$ generated by the involutions 
$\si_1$, $\si_2$, $\si_3$. 
It will turn out that $G$ is a universal Coxeter group of rank 
three with generators $\si_1$, $\si_2$, $\si_3$ 
(see Theorem \ref{thm:coxeter}). 
We are interested in the dynamics of the $G$-action on 
$\ol{\Sol}(\th)$. 
Usually the dynamics of a group action is more involved than 
that of a single transformation; more techniques and tools have 
been developed for the latter rather than for the former. 
So in this article we pick up each individual transformation 
from the group $G$ and study its dynamics, leaving the 
interaction of plural transformations in another occasion. 
\par
In order to study the dynamics of any element $\si \in G$, 
we begin with investigating the basic elements $\si_i$, 
$i = 1,2,3$, especially their behaviors in a neighborhood 
of the tritangent lines $L$ at infinity. 
To this end let us introduce the following three points 
\[
q_1 = [0:0:1:1], \qquad 
q_2 = [0:1:0:1], \qquad 
q_3 = [0:1:1:0],  
\]
where $q_i$ may be thought of as the ``mid-point" of $p_j$ and $p_k$ 
on the line $L_i$. 
\begin{figure}[t]
\begin{center}
\unitlength 0.1in
\begin{picture}(33.24,17.19)(4.90,-19.59)
%
\special{pn 13}%
\special{pa 490 1014}%
\special{pa 1111 1014}%
\special{fp}%
%
\special{pn 13}%
\special{pa 2641 1023}%
\special{pa 3397 1023}%
\special{fp}%
%
\special{pn 20}%
\special{ar 1858 1095 1053 855  3.3124258 6.1241532}%
%
\special{pn 20}%
\special{ar 1849 1104 1062 855  6.2694533 6.2831853}%
\special{ar 1849 1104 1062 855  0.0000000 3.1735968}%
\put(17.1400,-16.6200){\makebox(0,0)[lb]{$\Sol(\th)$}}%
%
\special{pn 13}%
\special{pa 1201 1014}%
\special{pa 2524 1023}%
\special{dt 0.045}%
\special{pa 2524 1023}%
\special{pa 2523 1023}%
\special{dt 0.045}%
\put(10.6000,-9.0000){\makebox(0,0)[lb]{$x$}}%
\put(25.8000,-9.1000){\makebox(0,0)[lb]{$x'$}}%
%
\special{pn 20}%
\special{sh 0.600}%
\special{ar 1102 1005 44 44  0.0000000 6.2831853}%
%
\special{pn 20}%
\special{sh 0.600}%
\special{ar 2632 1014 44 44  0.0000000 6.2831853}%
\put(18.1000,-6.1000){\makebox(0,0)[lb]{$\si_i$}}%
%
\special{pn 13}%
\special{pa 3040 1295}%
\special{pa 3814 1295}%
\special{fp}%
\special{sh 1}%
\special{pa 3814 1295}%
\special{pa 3747 1275}%
\special{pa 3761 1295}%
\special{pa 3747 1315}%
\special{pa 3814 1295}%
\special{fp}%
\put(31.8400,-15.2900){\makebox(0,0)[lb]{$x_i$-axis}}%
%
\special{pn 13}%
\special{ar 1880 1850 1161 1161  4.1675251 5.2764727}%
%
\special{pn 13}%
\special{pa 2380 800}%
\special{pa 2510 870}%
\special{fp}%
\special{sh 1}%
\special{pa 2510 870}%
\special{pa 2461 821}%
\special{pa 2463 845}%
\special{pa 2442 856}%
\special{pa 2510 870}%
\special{fp}%
\end{picture}%
\end{center}
\caption{Involutions of a $(2,2,2)$-surface} 
\label{fig:involution}
\end{figure}
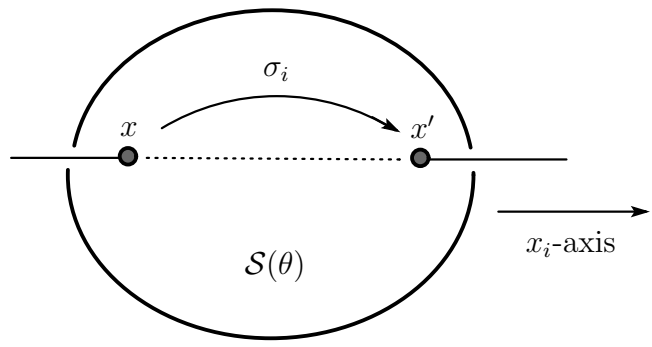
\begin{lemma} \label{lem:L} 
The birational map $\si_i$ has the following properties 
$($see Figure $\ref{fig:cubic2}$$)$. 
\begin{enumerate}
\item $\si_i$ blows down the line $L_i$ to the point $p_i$,  
\item $\si_i$ restricts to the automorphism of $L_j$ that fixes 
$q_j$ and exchanges $p_i$ and $p_k$, 
\item $\si_i$ restricts to the automorphism 
of $L_k$ that fixes $q_k$ and exchanges $p_i$ and $p_j$, 
\item $p_i$ is the unique indeterminacy point of $\si_i$, 
\end{enumerate}
\end{lemma} 
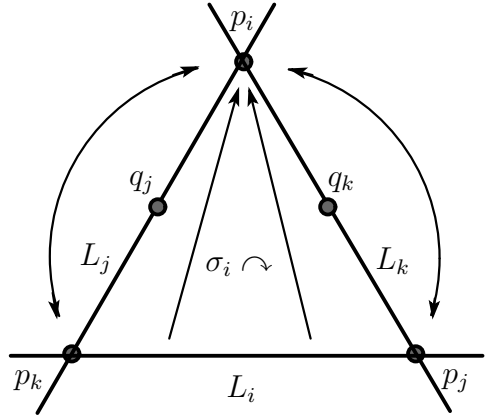
\begin{figure}[t] 
\begin{center}
\unitlength 0.1in
\begin{picture}(24.70,21.70)(6.00,-24.70)
\put(17.3000,-24.1000){\makebox(0,0)[lb]{$L_i$}}%
\put(9.6000,-17.4000){\makebox(0,0)[lb]{$L_j$}}%
\put(25.1000,-17.4000){\makebox(0,0)[lb]{$L_k$}}%
\put(28.6000,-23.6000){\makebox(0,0)[lb]{$p_j$}}%
\put(17.5000,-4.7000){\makebox(0,0)[lb]{$p_i$}}%
\put(6.2000,-23.5000){\makebox(0,0)[lb]{$p_k$}}%
\put(16.2000,-17.5000){\makebox(0,0)[lb]{$\si_i \car$}}%
\put(22.6000,-13.1000){\makebox(0,0)[lb]{$q_k$}}%
\put(12.2000,-13.2000){\makebox(0,0)[lb]{$q_j$}}%
%
\special{pn 20}%
\special{sh 0.600}%
\special{ar 1820 630 40 40  0.0000000 6.2831853}%
%
\special{pn 20}%
\special{sh 0.600}%
\special{ar 2720 2160 40 40  0.0000000 6.2831853}%
%
\special{pn 20}%
\special{sh 0.600}%
\special{ar 920 2160 40 40  0.0000000 6.2831853}%
%
\special{pn 20}%
\special{pa 1640 330}%
\special{pa 2900 2460}%
\special{fp}%
%
\special{pn 20}%
\special{pa 1980 330}%
\special{pa 740 2470}%
\special{fp}%
%
\special{pn 20}%
\special{pa 600 2170}%
\special{pa 3070 2170}%
\special{fp}%
%
\special{pn 20}%
\special{sh 0.600}%
\special{ar 1370 1390 40 40  0.0000000 6.2831853}%
%
\special{pn 20}%
\special{sh 0.600}%
\special{ar 2260 1390 40 40  0.0000000 6.2831853}%
%
\special{pn 13}%
\special{ar 1830 1660 1015 1015  4.9362935 6.2831853}%
\special{ar 1830 1660 1015 1015  0.0000000 0.2970642}%
%
\special{pn 13}%
\special{ar 1820 1650 1016 1016  2.8289604 4.4722180}%
%
\special{pn 13}%
\special{pa 2820 1870}%
\special{pa 2810 1960}%
\special{fp}%
\special{sh 1}%
\special{pa 2810 1960}%
\special{pa 2837 1896}%
\special{pa 2816 1907}%
\special{pa 2797 1892}%
\special{pa 2810 1960}%
\special{fp}%
%
\special{pn 13}%
\special{pa 820 1840}%
\special{pa 840 1960}%
\special{fp}%
\special{sh 1}%
\special{pa 840 1960}%
\special{pa 849 1891}%
\special{pa 831 1907}%
\special{pa 809 1898}%
\special{pa 840 1960}%
\special{fp}%
%
\special{pn 13}%
\special{pa 2180 700}%
\special{pa 2070 670}%
\special{fp}%
\special{sh 1}%
\special{pa 2070 670}%
\special{pa 2129 707}%
\special{pa 2121 684}%
\special{pa 2140 668}%
\special{pa 2070 670}%
\special{fp}%
%
\special{pn 13}%
\special{pa 1450 700}%
\special{pa 1570 660}%
\special{fp}%
\special{sh 1}%
\special{pa 1570 660}%
\special{pa 1500 662}%
\special{pa 1519 677}%
\special{pa 1513 700}%
\special{pa 1570 660}%
\special{fp}%
%
\special{pn 13}%
\special{pa 1430 2080}%
\special{pa 1790 780}%
\special{fp}%
\special{sh 1}%
\special{pa 1790 780}%
\special{pa 1753 839}%
\special{pa 1776 831}%
\special{pa 1791 850}%
\special{pa 1790 780}%
\special{fp}%
%
\special{pn 13}%
\special{pa 2170 2080}%
\special{pa 1850 780}%
\special{fp}%
\special{sh 1}%
\special{pa 1850 780}%
\special{pa 1847 850}%
\special{pa 1863 832}%
\special{pa 1885 840}%
\special{pa 1850 780}%
\special{fp}%
\end{picture}%
\end{center}
\caption{The birational map $\si_i$ restricted to $L$}
\label{fig:cubic2} 
\end{figure}
{\it Proof}. 
In order to investigate $\si_i$, we use the inhomogeneous 
coordinates of $\P^3$ in (\ref{eqn:coord}) and local 
coordinates and local equations of $\ol{\Sol}(\th)$ in 
Table \ref{tab:local}, with target coordinates 
being dashed. 
\par 
In terms of the inhomogeneous coordinates $v$ and $u'$ of $\P^3$, 
the map $\si_i : v \mapsto u'$ is given by 
\begin{equation} \label{eqn:vu'}
u_0' = \frac{v_0^2}{\th_i v_0^2 - v_0 v_i - v_k}, \qquad  
u_j' = \frac{v_0}{\th_i v_0^2 - v_0 v_i - v_k},   \qquad 
u_k' = \frac{v_0 v_k}{\th_i v_0^2 - v_0 v_i - v_k}. 
\end{equation}
In a neighborhood of $L_i-\{p_j,p_k\}$ in $\ol{\Sol}(\th)$, 
using $v_i = O(v_0)$, we observe that 
\[
\th_i v_0^2-v_0 v_i-v_k = -v_k \{1+O(v_0^2)\}, 
\]
which is substituted into (\ref{eqn:vu'}) to yield 
\[
u_j' = -\frac{v_0}{v_k \{1+O(v_0^2)\}} = 
-\frac{v_0}{v_k}\{1 + O(v_0^2)\}, 
\qquad 
u_k' = -\frac{v_0 v_k}{v_k \{1+O(v_0^2)\}} = -v_0\{1 + O(v_0^2)\}.
\]
In particular putting $v_0 = 0$ leads to $u_j' = u_k' = 0$. 
This means that $\si_i$ maps a neighborhood of $L_i - \{p_j,p_k \}$ 
to a neighborhood of $p_i$, collapsing $L_i - \{p_j,p_k \}$ 
to the single point $p_i$. 
\par 
In a similar manner, in a neighborhood of $p_j$ in $\ol{\Sol}(\th)$ 
we observe that 
\[
v_0 = -(v_i v_k) \{1+O_2(v_i,v_k)\}, \qquad 
\th_i v_0^2-v_0 v_i-v_k = -v_k \{1+O_2(v_i, v_k)\}, 
\] 
which are substituted into (\ref{eqn:vu'}) to yield 
\[
u_j' = v_i \{1+O_2(v_i,v_k)\}, 
\qquad 
u_k' = (v_i v_k) \{1+O_2(v_i,v_k)\}. 
\]
In particular putting $v_i = 0$ leads to $u_j' = u_k' = 0$. 
This means that $\si_i$ maps a neighborhood of $p_j$ to 
a neighborhood of $p_i$, collapsing a neighborhood in $L_i$ 
of $p_j$ to the single point $p_i$. 
Using $(w, u')$ in place of $(v, u')$, we can argue similarly 
in a neighborhood of $p_k$.  
Therefore $\si_i$ blows down $L_i$ to the point $p_i$, which 
proves assertion (1). 
Moreover it is clear from the argument that there is no 
indeterminacy point on the line $L_i$. 
\par 
In terms of the inhomogeneous coordinates $u$ and $u'$ of 
$\P^3$ the map $\si_i : u \mapsto u'$ is given by 
\begin{equation} \label{eqn:siu}
u_0' = \frac{u_0^2}{\th_i u_0^2 - u_0 -u_j u_k},   \quad  
u_j' = \frac{u_0 u_j}{\th_i u_0^2 - u_0 -u_j u_k}, \quad 
u_k' = \frac{u_0 u_k}{\th_i u_0^2 - u_0 -u_j u_k}.
\end{equation}
In a neighborhood of $L_j - \{ p_i,p_k \}$ in $\ol{\Sol}(\th)$, 
using $u_j = -(u_k + 1/u_k) u_0 + O(u_0^2)$, we have 
\[
\th_i u_0^2 - u_0 -u_j u_k = u_0 \{ u_k^2 + O(u_0) \},
\] 
which is substituted into (\ref{eqn:siu}) to yield 
\[
u_0' = \frac{u_0}{u_k^2 + O(u_0)} = \frac{u_0}{u_k^2} + O(u_0^2), 
\qquad 
u_k' = \frac{u_k}{u_k^2 + O(u_0)} = \frac{1}{u_k} + O(u_0).
\]
In particular putting $u_0 = 0$ leads to $u_0' = 0$ and 
$u_k' = 1/u_k$. 
This means that $\si_i$ restricts to an automorphism of 
a neighborhood of $L_j - \{p_i,p_k\}$ in $\ol{\Sol}(\th)$ 
which induces a unique automorphism of $L_i$ fixing 
$q_j$ and exchanging $p_i$ and $p_k$. 
This proves assertion (2) and also shows that 
there is no indeterminacy point on $L_j-\{p_i,p_k\}$. 
Assertion (3) and the nonexistence of indeterminacy point 
on $L_k-\{p_i,p_j\}$ are established just in the same manner. 
\par 
From the above argument we have already known that 
there is no indeterminacy point other than $p_i$. 
Then the point $p_i$ is actually an indeterminacy point, 
because $\si_i$ is an involution blowing down 
$L_i$ to $p_i$ and hence blows up $p_i$ to $L_i$ reciprocally. 
This proves assertion (4).  
\hfill $\Box$ \par\medskip 
Later we will need some information about how the involution $\si_i$ 
transforms a line to another curve, which is stated in the following 
lemma. 
\begin{lemma} \label{lem:si} 
For any $\{i,j,k\} = \{1,2,3\}$, the involution $\si_i$ satisfies 
the following properties: 
\begin{enumerate} 
\item $\si_i(E_i)$ intersects $E_i$ at two points counted with 
multiplicity. 
Similarly, $\si_i(E_{i+3})$ intersects $E_{i+3}$ at two points 
counted with multiplicity. 
\item $\si_i(E_i)$ intersects $E_{i+3}$ at one point counted 
with multiplicity. 
Similarly, $\si_i(E_{i+3})$ intersects $E_i$ at one point 
counted with multiplicity. 
\item $\si_i$ exchanges the lines $E_j$ and $G_{j+3}$; 
$E_{j+3}$ and $G_j$; 
$E_k$ and $G_{k+3}$; $E_{k+3}$ and $G_k$, respectively. 
\end{enumerate}
\end{lemma} 
{\it Proof}. 
By Table \ref{tab:line} we may put $E_i = L_i(b_i,b_4;b_j,b_k)$ 
and $E_{i+3} = L_i(b_j,b_k;b_i,b_4)$. 
Assertion (1) of Lemma \ref{lem:L} implies that $\si_i(E_i)$ 
does not intersect $E_i$ nor $E_{i+3}$ at any point at infinity. 
So we can work with the inhomogeneous coordinates 
$x = (x_1,x_2,x_3)$. 
By (\ref{eqn:line2}) the line $E_i$ is given by 
\begin{equation} \label{eqn:Ei} 
x_i = b_ib_4 + (b_ib_4)^{-1}, \qquad 
x_j + (b_ib_4) x_k =  a_k b_i +  a_j b_4. 
\end{equation}
In a similar manner, by exchanging $(b_i,b_4)$ and $(b_j,b_k)$ in 
(\ref{eqn:line2}), the line $E_{i+3}$ is given by 
\begin{equation} \label{eqn:Ej} 
x_i = b_jb_k + (b_jb_k)^{-1}, \qquad 
x_j + (b_jb_k) x_k =  a_4 b_j +  a_i b_k. 
\end{equation}
Moreover, by applying formula (\ref{eqn:si}) to (\ref{eqn:Ei}), 
the curve $\si_i(E_i)$ is expressed as 
\begin{equation} \label{eqn:si(Ei)}
\th_i-x_i-x_jx_k = b_ib_4+(b_ib_4)^{-1}, \qquad 
x_j + (b_ib_4) x_k =  a_k b_i +  a_j b_4. 
\end{equation}
Note that the second equations of 
(\ref{eqn:Ei}) and (\ref{eqn:si(Ei)}) are the same. 
\par 
In order to find out the intersection of $\si_i(E_i)$ with $E_i$, 
let us couple (\ref{eqn:Ei}) and (\ref{eqn:si(Ei)}) together. 
Eliminating $x_i$ and $x_j$ we obtain 
a quadratic equation for $x_k$, 
\[
(b_ib_4) x_k^2 -(a_k b_i + a_j b_4)x_k + 
\th_i-2\{b_ib_4 + (b_ib_4)^{-1}\} = 0. 
\]
For a simple root of this equation we have a simple intersection 
point of $\si_i(E_i)$ with $E_i$ and for a double root we have 
an intersection point of multiplicity two. 
This proves assertion (1) for the pair $\si_i(E_i)$ and $E_i$. 
The assertion (1) for $\si_i(E_{i+3})$ and $E_{i+3}$ is 
proved in a similar manner. 
\par 
Next, in order to find out the intersection of $\si_i(E_i)$ 
with $E_{i+3}$, let us couple (\ref{eqn:Ej}) and 
(\ref{eqn:si(Ei)}). 
From the first equation of (\ref{eqn:Ej}) the 
$x_i$-coordinate is already fixed. 
The second equations of (\ref{eqn:Ej}) and (\ref{eqn:si(Ei)}) 
are coupled to yield a linear system for $x_j$ and $x_k$, 
whose determinant 
\[
b_jb_k-b_ib_4= b_ib_4(b_i^{-1}b_jb_kb_4^{-1}-1)
\]
is nonzero by the assumption that $\ol{\Sol}(\th)$ is smooth, 
that is, the discriminant $\vD(\th)$ in (\ref{eqn:vD}) is nonzero. 
Then the linear system is uniquely solved to determine 
$x_j$ and $x_k$. 
Now we can check that the first equation of (\ref{eqn:si(Ei)}) 
is redundant, that is, automatically satisfied. 
Therefore $\si_i(E_i)$ and $E_{i+3}$ has a simple intersection, 
which implies assertion (2) for the pair $\si_i(E_i)$ and $E_{i+3}$. 
The assertion (2) for $\si_i(E_{i+3})$ and $E_i$ is proved in a 
similar manner. 
\par 
Finally we see that $\si_i$ exchanges $E_j$ and $G_{j+3}$. 
We may put $E_j = L_j(b_j,b_4;b_k,b_i)$ and 
$G_{j+3} = L_j(1/b_j,1/b_4;b_k,b_i)$. 
By formula (\ref{eqn:line2}) (with indices suitably permuted), 
these lines are given by 
\begin{align} 
x_j &= b_jb_4 + (b_jb_4)^{-1}, & 
x_k + (b_jb_4) x_i &= a_ib_j + a_kb_4, \label{eqn:Ek} 
\\[2mm] 
x_j &= b_jb_4 + (b_jb_4)^{-1}, &
x_k + (b_jb_4)^{-1} x_i &= a_i b_j^{-1} + a_k b_4^{-1}. 
\label{eqn:Gl} 
\end{align} 
Using formula (\ref{eqn:si}) we can check that equations 
(\ref{eqn:Ek}) and (\ref{eqn:Gl}) are transformed to each other 
by $\si_i$. 
This together with similar argument for the other lines 
establishes assertion (3). \hfill $\Box$ 
\section{Cohomological Action} \label{sec:cohomology} 
A general theory of the dynamical system for a 
bimeromorphic map of a surface is developed in \cite{DF}. 
The basic strategy employed there is to consider the 
induced action of the map on the $(1,1)$-cohomology group, 
taking into account the influence of its exceptional set 
and indeterminacy set. 
In this section we shall use this technique in our context. 
\par 
Let $S$ be a compact complex surface, $f : S \carl$ 
a bimeromorphic map. 
Then $f$ is represented by a compact complex surface $\vG$, 
called the desingularized graph of $f$, together with proper 
modifications $\pi_1 : \vG \to S$ and 
$\pi_2 : \vG \to S$ such that $f = \pi_2 \ci \pi_1^{-1}$ 
on a dense open subset. 
For $i = 1, 2$, let $\E(\pi_i) := \{\, x \in \vG \,:\, 
\mathrm{\#} \, \pi_i^{-1}(\pi_i (x)) = \infty \,\}$ 
be the exceptional set for the projection $\pi_i$. 
The images $\E(f) := \pi_1 (\E(\pi_2))$  and 
$I(f) := \pi_1 (\E(\pi_1))$ are called the exceptional set 
and the indeterminacy set of $f$ respectively. 
Between these sets there is a useful relation 
\begin{equation}  \label{eqn:EI} 
f(\E(f)) = I(f^{-1}). 
\end{equation} 
If $S = \ol{\Sol}(\th)$ and $f = \si_i$, then 
Lemma \ref{lem:L} readily leads to the following lemma. 
\begin{lemma} \label{lem:IE} 
For each $i \in \{1, 2, 3\}$, we have 
$\E(\si_i) = L_i$, $\si_i(\E(\si_i)) = \{p_i\}$ and 
$I(\si_i) = \{p_i\}$.
\end{lemma}
\par 
Given any element $\si \in G$ other than the unit element, 
we can write 
\begin{equation} \label{eqn:reduced3}
\si = \si_{i_1} \si_{i_2} \cdots \si_{i_n},  
\end{equation} 
for some $n \in \N$ and some $n$-tuple of indices 
$(i_1,\dots,i_n) \in \{1,2,3\}^n$ such that every neighboring 
indices $i_{\nu}$ and $i_{\nu+1}$ are distinct. 
It is not yet clear at this stage whether the expression 
(\ref{eqn:reduced3}) is unique or not, though the uniqueness 
will be established later (see Theorem \ref{thm:coxeter}). 
In any case, we begin with the determination of the exceptional 
set and the indeterminacy set of $\si$. 
\begin{lemma} \label{lem:EI2} 
For the expression $(\ref{eqn:reduced3})$ we have 
\begin{equation} \label{eqn:EI2}
\E(\si) = \bigcup_{\nu=1}^n L_{i_{\nu}}, \qquad 
\si(\E(\si)) = \{p_{i_1}\}, \qquad I(\si) = \{p_{i_n}\}. 
\end{equation}
\end{lemma}
{\it Proof}. Let us prove the first formula of (\ref{eqn:EI2}) 
by induction on the length $n$. 
For $n = 1$ the assertion immediately follows from 
Lemma \ref{lem:IE}. 
Assume that the assertion holds when the length is $n-1$ and 
consider the element $\si' = \si_{i_2} \cdots \si_{i_n}$ of 
length $n-1$. 
Since $p_{i_1}$ and $p_{i_2}$ are distinct, 
we have $I(\si_{i_1}) \cap \si'(\E(\si')) = 
\{p_{i_1}\} \cap \{p_{i_2}\} = \emptyset$ and hence 
$\E(\si') \subset \E(\si)$. 
Therefore, 
\begin{equation} \label{eqn:induction}
\bigcup_{\nu=2}^n L_{i_{\nu}} \subset 
\E(\si) \subset \bigcup_{\nu=1}^n L_{i_{\nu}}, 
\end{equation} 
where the first inclusion follows from the induction hypothesis 
and the second inclusion is easily seen from Lemma \ref{lem:L}. 
If $i_1 \in \{i_2,\dots,i_n\}$, then the leftmost and rightmost 
sets in (\ref{eqn:induction}) are the same and hence all the 
three coincide. 
If $i_1 \not\in \{i_2,\dots,i_n\}$, then Lemma \ref{lem:L} 
implies that $\si'$ maps $L_{i_1}$ isomorphically onto itself 
and then $\si_{i_1}$ blows down $L_{i_1}$ to the single point 
$p_{i_1}$. 
This means that $L_{i_1} \subset \E(\si)$ and hence the 
second inclusion in (\ref{eqn:induction}) becomes equality. 
Thus the assertion is verified for length $n$ and the induction 
is complete.  
\par 
The second formula in (\ref{eqn:EI2}) is also proved by induction 
on the length $n$. 
For $n = 1$ the assertion immediately follows from 
Lemma \ref{lem:IE}. 
Assume that the assertion holds when the length is $n-1$. 
Then we have $\si'(\E(\si')) = \{p_{i_2}\}$ by induction hypothesis 
and hence $\si(\E(\si)) = \si_{i_1}(\E(\si_{i_1}) \cup 
\si'(\E(\si'))) 
= \si_{i_1}(L_{i_1} \cup \{p_{i_2}\}) = \si_{i_1}(L_{i_1}) = 
\{p_{i_1}\}$, since $p_{i_2} \in L_{i_1}$. 
This shows that the assertion is verified for length $n$ and 
hence the induction is complete. 
\par 
Next we prove the last formula of (\ref{eqn:EI2}). 
Instead of $\si$ we consider its inverse $\si^{-1}$. 
Since $\si^{-1} = \si_{i_n} \cdots \si_{i_2} \si_{i_1}$, 
the second formula of (\ref{eqn:EI2}) yields 
$\si^{-1}(\E(\si^{-1})) = \{p_{i_n}\}$. 
Then applying formula (\ref{eqn:EI}) to $f = \si^{-1}$, 
we have $I(\si) = \{p_{i_n}\}$. 
Thus the lemma is established. \hfill $\Box$ \par\medskip
If $S$ is a K\"ahler surface, two natural actions of $f$, 
pull-back and push-forward, on the Dolbeault 
cohomology group $H^{1,1}(S)$ are defined in the following 
manner: 
A smooth $(1,1)$-form $\omega$ on $S$ can be pulled back as 
a smooth $(1,1)$-form $\pi_2^* \omega$ on $\vG$ and then 
pushed forward as a $(1,1)$-current $\pi_{1*}\pi_2^* \omega$ 
on $S$. 
Hence we define the pull-back 
$f^* \omega := \pi_{1*}\pi_2^* \omega$ and also the 
push-forward 
$f_* \omega = (f^{-1})^* \omega := \pi_{2*}\pi_1^* \omega$. 
The operators $f^*$ and $f_*$ commute with the exterior 
differential $d$ and the complex structure of $S$ and 
so descend to linear actions on $H^{1,1}(S)$. 
For general bimeromorphic maps $f$ and $g$, the 
composition rule $(f\ci g)^* = g^* \ci f^*$ is not 
necessarily true. 
However a useful criterion under which this rule becomes 
true is given in \cite{DF}. 
\begin{lemma} \label{lem:compo}
If $I(f) \cap f(\E(g)) = \emptyset$, then 
$(f \ci g)^* = g^* \ci f^* : H^{1,1}(S) \carl$.
\end{lemma}
\par 
We shall apply this lemma to our biratinal transformation $\si$ 
in (\ref{eqn:reduced3}). 
\begin{lemma} \label{lem:compo2} 
For the expression $(\ref{eqn:reduced3})$ we have 
$\si^{*} = \si_{i_n}^{*} \cdots \si_{i_2}^{*} \si_{i_1}^{*} 
: H^{1,1}(\ol{\Sol}(\th)) \carl$. 
\end{lemma}
{\it Proof}. 
We prove the lemma by induction on the length $n$. 
It is trivial when $n = 1$. 
Assume that the lemma holds when the length is $n-1$. 
If we put $\si' = \si_{i_2} \cdots \si_{i_n}$, 
then the induction hypothesis implies that 
$(\si')^* = \si_{i_n}^{*} \cdots \si_{i_2}^{*}$. 
Lemma \ref{lem:EI2} shows that 
$I(\si_{i_1}) \cap \si'(\E(\si')) = 
\{p_{i_1}\} \cap \{ p_{i_2}\} = \emptyset$, 
since $p_{i_1}$ and $p_{i_2}$ are distinct. 
We now apply Lemma \ref{lem:compo} to $f = \si_{i_1}$ and 
$g = \si'$ to obtain $\si^{*} = (\si_{i_1} \si')^{*} = 
(\si')^{*} \si_{i_1}^{*} = \si_{i_n}^{*} \cdots \si_{i_2}^{*} 
\si_{i_1}^{*}$. 
Thus the lemma is true for length $n$. 
\hfill $\Box$ \par\medskip
By Lemma \ref{lem:compo2} the calculation of the action 
$\si^{*} : H^{1,1}(\ol{\Sol}(\th)) \carl$ is reduced to 
that of the actions $\si_1^*$, $\si_2^*$, 
$\si_3^* : H^{1,1}(\ol{\Sol}(\th)) \carl$, 
which is now set forth. 
Since the cubic surface $\ol{\Sol}(\th)$ is rational, 
we have $H^{1,1}(\ol{\Sol}(\th)) = H^2(\ol{\Sol}(\th),\C)$, 
where the latter group is described in (\ref{eqn:basis}). 
\begin{table}[t]
\[
\begin{array}{rl}
\si_1^* = \left(\begin{array}{rrrrrrr}
                 6 &  3 &  2 &  2 &  3 &  2 &  2 \\ 
                -3 & -2 & -1 & -1 & -1 & -1 & -1 \\ 
                -2 & -1 & -1 & -1 & -1 &  0 & -1 \\ 
                -2 & -1 & -1 & -1 & -1 & -1 &  0 \\ 
                -3 & -1 & -1 & -1 & -2 & -1 & -1 \\ 
                -2 & -1 &  0 & -1 & -1 & -1 & -1 \\ 
                -2 & -1 & -1 &  0 & -1 & -1 & -1 
            \end{array}\right) 
\quad & 
\si_2^* = \left(\begin{array}{rrrrrrr}
                 6 &  2 &  3 &  2 &  2 &  3 &  2 \\
                -2 & -1 & -1 & -1 &  0 & -1 & -1 \\ 
                -3 & -1 & -2 & -1 & -1 & -1 & -1 \\ 
                -2 & -1 & -1 & -1 & -1 & -1 &  0 \\ 
                -2 &  0 & -1 & -1 & -1 & -1 & -1 \\ 
                -3 & -1 & -1 & -1 & -1 & -2 & -1 \\
                -2 & -1 & -1 &  0 & -1 & -1 & -1  
            \end{array}\right) 
\\ &  \\
\si_3^* = \left(\begin{array}{rrrrrrr}
                6 &  2 &  2 &  3 &  2 &  2 &  3 \\ 
               -2 & -1 & -1 & -1 &  0 & -1 & -1 \\
               -2 & -1 & -1 & -1 & -1 &  0 & -1 \\
               -3 & -1 & -1 & -2 & -1 & -1 & -1 \\
               -2 &  0 & -1 & -1 & -1 & -1 & -1 \\
               -2 & -1 &  0 & -1 & -1 & -1 & -1 \\
               -3 & -1 & -1 & -1 & -1 & -1 & -2  
            \end{array}\right)
\quad & 
 c^* = \left(\begin{array}{rrrrrrr}
               12 &  6 &  4 &  3 &  6 &  4 &  3 \\
               -3 & -2 & -1 & -1 & -1 & -1 & -1 \\ 
               -4 & -2 & -2 & -1 & -2 & -1 & -1 \\ 
               -6 & -3 & -2 & -2 & -3 & -2 & -1 \\ 
               -3 & -1 & -1 & -1 & -2 & -1 & -1 \\ 
               -4 & -2 & -1 & -1 & -2 & -2 & -1 \\ 
               -6 & -3 & -2 & -1 & -3 & -2 & -2  
        \end{array}\right)
\end{array}
\]
\caption{Matrix representations of $\si_1^*$, $\si_2^*$, 
$\si_3^*$, $c^* : H^2(\ol{\Sol}(\th), \Z) \carl$, 
where $c = \si_1 \si_2 \si_3$}
\label{tab:matrix} 
\end{table} 
\begin{lemma} \label{lem:in2} 
The linear operators $\si_1^*$, $\si_2^*$, $\si_3^* :
 H^2(\ol{\Sol}(\th), \Z) \carl$ have matrix 
representations as in Table $\ref{tab:matrix}$ 
with respect to the basis in $(\ref{eqn:basis})$.  
\end{lemma} 
{\it Proof}. 
First we shall find the matrix representation of $\si_1^*$. 
If $\xi_{ab}$ denotes the $(a,b)$-th entry of the matrix to 
be found, where $0 \le a, b \le 6$, then (\ref{eqn:in1}) 
implies that 
\[
\si_1^* E_b = \sum_{a=0}^6 \xi_{ab} \, E_a 
= \sum_{a=0}^6 \delta_a (\si_1^* E_b, E_a) \, E_a, 
\]
where we put $\delta_a = 1$ for $a = 0$ and $\delta_a = -1$ for 
$a \neq 0$. 
Now we claim that 
\begin{equation} \label{eqn:xi}
\xi_{ab} = \delta_a (\si_1^*E_b, E_a), \qquad 
\xi_{ab} = \delta_a \delta_b \xi_{ba}. 
\end{equation}
The first formula in (\ref{eqn:xi}) is obvious and the second 
formula is derived as follows: 
\[
\xi_{ab} = \delta_a (\si_1^*E_b, E_a) = 
\delta_a (E_b, \si_{1*}E_a) = \delta_a (E_b, \si_1^*E_a) 
= (\delta_a \delta_b) \cdot \delta_b (\si_1^*E_a, E_b) = 
(\delta_a \delta_b) \xi_{ba}, 
\]
where in the third equality we have used the fact that 
$\si_1$ is an involution; $\si_{1*} = (\si_1^{-1})^* = 
\si_1^*$. 
By assertions (1) and (2) of Lemma \ref{lem:si} we have 
$(\si_1^*E_1,E_1) = 2$ and $(\si_1^*E_1, E_4) = 1$ and 
likewise $(\si_1^*E_4,E_4) = 2$ and 
$(\si_1^*E_4, E_1) = 1$. 
Then the first formula of (\ref{eqn:xi}) yields 
\begin{equation} \label{eqn:xi2} 
\xi_{11} = \xi_{44} = -2, \qquad 
\xi_{14} = \xi_{41} = -1. 
\end{equation}
The assertion (3) of Lemma \ref{lem:si} together 
with the second formula of (\ref{eqn:FG}) yields 
\begin{equation} \label{eqn:si3} 
\left\{
\begin{array}{rcl}
\si_1^*E_2 &=&2 E_0-E_1-E_2-E_3-E_4 \phantom{- E_5\,}-E_6, \\[1mm]
\si_1^*E_3 &=&2 E_0-E_1-E_2-E_3-E_4-E_5,                   \\[1mm]
\si_1^*E_5 &=&2 E_0-E_1 \phantom{- E_2\,}-E_3-E_4-E_5-E_6, \\[1mm]
\si_1^*E_6 &=&2 E_0-E_1-E_2 \phantom{- E_3\,}-E_4-E_5-E_6, 
\end{array}
\right.
\end{equation}
It follows from (\ref{eqn:xi2}) and (\ref{eqn:si3}) that the 
matrix representation for $\si_1^*$ takes the form 
\begin{equation} \label{eqn:si4} 
\si_1^* = 
\left(\begin{array}{rrrrrrr}
          * &       * &  2 &  2 &       * &  2 &  2 \\
          * &      -2 & -1 & -1 &      -1 & -1 & -1 \\
    \bullet & \bullet & -1 & -1 & \bullet &  0 & -1 \\
    \bullet & \bullet & -1 & -1 & \bullet & -1 &  0 \\
          * &      -1 & -1 & -1 &      -2 & -1 & -1 \\
    \bullet & \bullet & 0  & -1 & \bullet & -1 & -1 \\
    \bullet & \bullet & -1 &  0 & \bullet & -1 & -1 
\end{array}\right), 
\end{equation}
where the entries denoted by $\bullet$ and $*$ are yet to be 
determined. 
The entries denoted by $\bullet$ are easily determined 
by the second formula in (\ref{eqn:xi}). 
The final ingredient taken into account is the fact that 
$\si_1$ blows down $L_1 = E_0-E_1-E_4$ to the point $p_1$ 
(see Lemma \ref{lem:L}), which leads to 
\[
\si_1^* E_0 -\si_1^* E_1 - \si_1^* E_4 = 0. 
\]
This means that the $0$-th column is the sum of the 
first and fourth columns in the matrix (\ref{eqn:si4}). 
Using the second formula in (\ref{eqn:xi}) repeatedly, we 
see that (\ref{eqn:si4}) becomes the first matrix of 
Table \ref{tab:matrix}. 
The matrix representations of $\si_2^{*}$ and $\si_3^{*}$ 
are obtained just in the same manner. 
\hfill $\Box$ \par\medskip
In order to make Lemma \ref{lem:in2} more transparent, we 
consider the direct sum decomposition 
\begin{equation} \label{eqn:decomp}
H^2(\ol{\Sol}(\th),\C) = V \oplus V^{\perp}, 
\end{equation} 
where $V$ is the subspace spanned by the lines 
$L_1$, $L_2$, $L_3$ at infinity and $V^{\perp}$ is the 
orthogonal complement to it with respect to the 
intersection form. 
In view of (\ref{eqn:lineinf}) and (\ref{eqn:FG}), we have 
\[
L_1 = E_0-E_1-E_4, \qquad 
L_2 = E_0-E_2-E_5, \qquad
L_3 = E_0-E_3-E_6. 
\]
On the other hand, it is easily seen that the subspace 
$V^{\perp}$ is spanned by the vectors 
\[
2 E_0-E_1-E_2-E_3-E_4-E_5-E_6, \qquad 
E_1-E_4, \qquad E_2-E_5, \qquad E_3-E_6. 
\]
A little calculation in terms of the new basis shows that 
Lemma \ref{lem:in2} can be restated as follows. 
\begin{lemma} \label{lem:in3} 
The linear operators $\si_1^*$, $\si_2^*$, $\si_3^* :
H^2(\ol{\Sol}(\th), \Z) \carl$ preserve the subspaces 
$V$ and $V^{\perp}$. 
They act on these subspaces in the following manner. 
\begin{enumerate}
\item 
The operators $\si_1^*$, $\si_2^*$, $\si_3^*$ restricted 
to $V$ are represented by the matrices 
\begin{equation} \label{eqn:matrices}
s_1 = \begin{pmatrix} 0 & 1 & 1 \\ 0 & 1 & 0 \\ 0 & 0 & 1 
\end{pmatrix}, \qquad 
s_2 = \begin{pmatrix} 1 & 0 & 0 \\ 1 & 0 & 1 \\ 0 & 0 & 1 
\end{pmatrix}, \qquad 
s_3 = \begin{pmatrix} 1 & 0 & 0 \\ 0 & 1 & 0 \\ 1 & 1 & 0 
\end{pmatrix}, 
\end{equation}
respectively, with respect to the basis $L_1$, $L_2$, $L_3$. 
\item The operators $\si_1^*$, $\si_2^*$, $\si_3^*$ act on 
$V^{\perp}$ as the negative of identity $-1$. 
\end{enumerate}
\end{lemma} 
It should be noted that each matrix in (\ref{eqn:matrices}) 
has eigenvalues $0$, $1$, $1$, counted with multiplicities, 
and in particular has vanishing determinant. 
\begin{theorem} \label{thm:coxeter}
The group $G = \la \si_1, \si_2,\si_3 \ra$ is a universal 
Coxeter group of rank three over the basic involutions 
$\si_1$, $\si_2$, $\si_3$, that is, there are no 
relations other than $\si_1^2 = \si_2^2 = \si_3^2 = 1$. 
In particular the expression $(\ref{eqn:reduced3})$ is 
unique for any given element $\si \in G$. 
\end{theorem} 
{\it Proof}. 
Assume the contrary that there exists a nontrivial 
relation $\si_{i_1} \si_{i_2} \cdots \si_{i_n} = 1$ 
in $G$ such that each neighboring indices 
$i_{\nu}$ and $i_{\nu+1}$ are distinct. 
Then it follows from Lemma \ref{lem:compo2} that 
$\si_{i_n}^* \cdots \si_{i_2}^* \si_{i_1}^* = 1^* = 1$ as a 
linear endomorphism on $H^2(\ol{\Sol}(\th),\C)$. 
But this is impossible because each factor $\si_{i_{\nu}}^*$ 
has vanishing determinant. 
This contradiction establishes the theorem. 
\hfill $\Box$ \par\medskip 
\begin{remark} \label{rem:coxeter} 
Recall that we have introduced the universal Coxeter group 
$G$ of rank three abstractly in \S \ref{sec:result}. 
Theorem \ref{thm:coxeter} yields a concrete 
realization of it as a group of birational transformations 
on the cubic surface $\ol{\Sol}(\th)$. 
Hereafter the former group will be identified with the latter. 
In this context the $3$-dimensional abstract linear space 
$V$ for the geometric representation $\mathrm{GR} : 
G \to \mathrm{O}_B(V)$ in \S \ref{sec:result} 
is realized as the subspace of $H^2(\ol{\Sol}(\th),\C)$ 
spanned by the lines at infinity $L_1$, $L_2$, $L_3$. 
Here we should put $e_1 = L_1$, $e_2 = L_2$, $e_3 = L_3$ 
in accordance with the notation in \S \ref{sec:result}. 
The symmetric bilinear form $B$ in (\ref{eqn:B}) is now given 
by the negative of the intersection form on 
$H^2(\ol{\Sol}(\th),\C)$ restricted to the subspace $V$. 
The basic reflections in (\ref{eqn:ri}) are then represented 
by the matrices 
\[
r_1 = \left(\begin{array}{rrr} 
-1 & 2 & 2 \\ 
 0 & 1 & 0 \\ 
 0 & 0 & 1 
\end{array}\right), 
\qquad 
r_2 = \left(\begin{array}{rrr} 
 1 &  0 & 0 \\ 
 2 & -1 & 2 \\ 
 0 &  0 & 1 
\end{array}\right), \qquad 
r_3 = \left(\begin{array}{rrr} 
1 & 0 &  0 \\ 
0 & 1 &  0 \\ 
2 & 2 & -1 
\end{array}\right). 
\]
It is easy to see that the linear operators $s_1$, $s_2$, 
$s_3$ in (\ref{eqn:esi}) have matrix representations 
as in (\ref{eqn:matrices}) and hence correspond to 
the operators $\si_1^*$, $\si_2^*$, 
$\si_3^*$ restricted to $V$. 
So the trace $\a(\ga)$ in (\ref{eqn:alpha}) can be 
calculated practically by using the matrix representations 
(\ref{eqn:matrices}). 
\end{remark} 
\par 
Next we shall calculate the characteristic polynomial 
of the linear map $\si^* : H^2(\ol{\Sol}(\th), \C) \carl$. 
In general the characteristic polynomial of a linear 
endomorphism $A$ is denoted by 
\[
P(\l; A) = \det(\l I - A). 
\]
For the reduced expression (\ref{eqn:reduced3}) of the element 
$\si$, we put $s_{\si} := s_{i_n} \cdots s_{i_2} s_{i_1}$ and 
define 
\begin{equation} \label{eqn:alpha2}
\a(\si) := \Tr [\,s_{\si} : V \to V \,]. 
\end{equation} 
\begin{lemma} \label{lem:char} 
The map $\si^*$ preserves the direct sum decomposition 
$(\ref{eqn:decomp})$ and hence factors as 
$\si^* = (\si^*|_V) \oplus (\si^*|_{V^{\perp}})$. 
The characteristic polynomial of the first component 
$\si^*|_V$ is given by 
\begin{equation} \label{eqn:char}
P(\l;\si^*|_V) = \left\{
\begin{array}{ll} 
\l \{\l^2- \a(\si) \l + (-1)^{n-1} \} \qquad & 
(\mbox{if} \,\,\, i_1 = i_n), \\[2mm]
\l \{\l^2- \a(\si) \l + (-1)^n \} \qquad & 
(\mbox{if} \,\,\, i_1 \neq i_n). 
\end{array} \right.
\end{equation}
The second component $\si^*|_{V^{\perp}}$ is just a scalar 
operator $(-1)^n$ having the characteristic polynomial 
\[
P(\l;\si^*|_{V^{\perp}}) = \{\l-(-1)^n \}^4. 
\] 
\end{lemma}
{\it Proof}. 
By Lemma \ref{lem:compo2} we have  
$\si^* = \si_{i_n}^* \cdots \si_{i_2}^* \si_{i_1}^*$. 
Hence the map $\si^*$ preserves the decomposition 
(\ref{eqn:decomp}), because each factor $\si^*_{i_{\nu}}$ 
does so by Lemma \ref{lem:in3}. 
Thus there are factorizations 
$\si^* = (\si^*|_V) \oplus (\si^*|_{V^{\perp}})$ and  
$P(\l;\si^*) = P(\l;\si^*|_V) P(\l;\si^*|_{V^{\perp}})$. 
The second component $\si^*|_{V^{\perp}}$ is found 
\[
\si^*|_{V^{\perp}} = (\si_{i_n}^*|_{V^{\perp}}) \cdots 
(\si_{i_1}^*|_{V^{\perp}}) = (-1)^n. 
\]
since each factor $\si_{i_{\nu}}^*$ restricted to $V^{\perp}$ 
is the scalar operator $-1$ by assertion (2) of Lemma \ref{lem:in3}. 
\par 
It remains to consider the first component $\si^*|_V$, which 
is represented by the three-by-three matrix 
$s_{\si} = s_{i_n} \cdots s_{i_2} s_{i_1}$. 
The argument will be based on the general fact that the 
characteristic polynomial of a three-by-three matrix $A$ 
is given by 
\begin{equation} \label{eqn:charpoly}
P(\l;A) = \l^3 - (\Tr\, A) \l^2 + (\Tr\, \wt{A}) \l - \det A, 
\end{equation}
where $\wt{A}$ is the adjugate matrix of $A$, namely, 
the matrix $\wt{A}$ such that $A \wt{A} = \wt{A} A = 
(\det A)I$. 
Let us apply this formula to $A = s_{\si}$. 
First we have $\Tr(s_{\si}) = \a(\si)$ by definition 
(\ref{eqn:alpha2}). 
Secondly we have $\det(s_{\si}) = 0$, since each 
factor $\si_{i_{\nu}}$ has vanishing determinant. 
Finally we wish to calculate the trace $\Tr(\tilde{s}_{\si})$. 
The general formula $\wt{AB} = \wt{B} \wt{A}$ for the product 
of adjugate matrices yields 
$\tilde{s}_{\si} = \tilde{s}_{i_1} \tilde{s}_{i_2} 
\cdots \tilde{s}_{i_n}$. 
Now it follows from (\ref{eqn:matrices}) that 
\begin{equation} \label{eqn:adjugate} 
\tilde{s}_1 = 
\left(\begin{array}{rrr} 
1 & -1 & -1 \\ 
0 &  0 &  0 \\
0 &  0 &  0
\end{array}\right), 
\qquad 
\tilde{s}_2 = 
\left(\begin{array}{rrr} 
 0 & 0 &  0 \\ 
-1 & 1 & -1 \\
 0 & 0 &  0
\end{array}\right), 
\qquad 
\tilde{s}_3 = 
\left(\begin{array}{rrr} 
 0 &  0 & 0 \\ 
 0 &  0 & 0 \\
-1 & -1 & 1 
\end{array}\right). 
\end{equation}
Note that among the three rows of the matrix $\tilde{s}_i$, 
only the $i$-th row does not vanish. 
Thus the only row of $\tilde{s}_{\si}$ that can be nonzero 
is the $i_1$-th row, so that the trace $\Tr(\tilde{s}_{\si})$ 
is just given by the $(i_1,i_1)$-th entry of $\tilde{s}_{\si}$. 
Now the latter quantity is calculated as 
\[
(\tilde{s}_{i_1})_{i_1i_2} (\tilde{s}_{i_2})_{i_2i_3} 
\cdots (\tilde{s}_{i_{n-1}})_{i_{n-1}i_n} 
(\tilde{s}_{i_n})_{i_ni_1}, 
\]
where $(\tilde{s}_i)_{ij}$ denotes the 
$(i,j)$-th entry of the matrix $\tilde{s}_i$. 
It follows from (\ref{eqn:adjugate}) that 
$(\tilde{s}_i)_{ij}$ is $+1$ or $-1$ according 
as the indices $i$ and $j$ are equal or not. 
Since $i_{\nu}$ and $i_{\nu+1}$ are 
distinct for every $\nu \in \{1, \dots, n-1\}$, 
we have $\Tr(\tilde{s}_{\si}) = (-1)^{n-1}$ or 
$\Tr(\tilde{s}_{\si}) = (-1)^n$ according as 
$i_n$ and $i_1$ are equal or not. 
Putting all these considerations into (\ref{eqn:charpoly}) 
yields formula (\ref{eqn:char}). \hfill $\Box$ \par\medskip 
\section{Ergodic Properties} \label{sec:ergodic} 
We continue to study the dynamical properties of each individual 
transformation $\si \in G$. 
The main concern in this section is the investigation into 
the ergodic properties of this map, where the notions of 
dynamical degree, entropy and invariant measure play 
important roles. 
It is a good application of the fundamental methods and 
techniques in bimeromorphic (or birational) surface dynamics, 
recently developed by \cite{BD,DF,DS,Dujardin}. 
Since they are not so familiar in the circle of 
Painlev\'e equations, we shall develop our discussion 
upon reviewing some rudiments of them. 
\par 
We begin with the concept of first dynamical degree \cite{DF}. 
Given a bimeromorphic map $f$ of a compact K\"ahler surface 
$S$, its {\sl first dynamical degree} $\l_1(f)$ is defined by 
\[
\l_1(f) := \lim_{N\to\infty} || (f^N)^* ||^{1/N}, 
\]
where $||\cdot||$ is an operator norm on 
$\End \, H^{1,1}(S)$. 
It is known that the limit certainly exists, 
independent of the norm $||\cdot||$ chosen, $\l_1(f) \ge 1$, 
and $\l_1(f)$ is invariant under bimeromorphic conjugation. 
It is usually difficult to evaluate this quantity in a simple mean. 
However there is a distinguished class of maps whose first 
dynamical degree can be equated to a more tractable quantity.  
A bimeromorphic map $f : S \carl$ is said to be 
{\sl analytically stable} (AS for short) if 
the condition $(f^n)^* = (f^*)^n : H^{1,1}(S) \carl$ holds  
for every $n \in \N$. 
Evidently, if $f$ is AS then 
\begin{equation} \label{eqn:s-radius}
\l_1(f) = \mathrm{SR}(f^*), 
\end{equation} 
where $\mathrm{SR}(f^*)$ is the spectral radius of 
the linear endomorphism $f^* : H^{1,1}(S) \carl$. 
It is known that any bimeromorphic map is bimeromorphically 
conjugate to an AS map. 
It is also known that a bimeromorphic map $f$ is AS if and 
only if 
\begin{equation} \label{eqn:AS}
 \bigcup_{N \ge 0} f^{-N} I(f)  
 \cap \bigcup_{N \ge 0} f^N I(f^{-1})  = \emptyset. 
\end{equation} 
This condition may be viewed as a separation between the 
obstructions to forward and backward dynamics. 
Back to our context, it is natural to ask when a given element 
$\si \in G$ is AS. 
\begin{lemma} \label{lem:as2} 
An element $\si \in G$ is AS if and only 
if the initial index $i_1$ and the terminal index $i_n$ are 
distinct in the reduced expression $(\ref{eqn:reduced3})$ of 
$\si$. 
\end{lemma} 
{\it Proof}. 
If $\si$ is AS then it follows from condition (\ref{eqn:AS}) 
that $I(\si) \cap I(\si^{-1}) = \emptyset$. 
On the other hand, Lemma \ref{lem:EI2} implies that 
$I(\si) = \{p_{i_n}\}$ and $I(\si^{-1}) = \{p_{i_1}\}$. 
Hence the points $p_{i_1}$ and $p_{i_n}$ must be distinct, 
that is, the indices $i_1$ and $i_n$ must be distinct. 
Conversely, assuming that the indices $i_1$ and $i_n$ are 
distinct, we shall show that for every $N \ge 0$, 
\begin{equation} \label{eqn:AS2}
\si^{-N} I(\si) = \{p_{i_n}\}, \qquad 
\si^N I(\si^{-1}) = \{p_{i_1}\}. 
\end{equation}
It suffices to verify the first formula of (\ref{eqn:AS2}), 
since the second formula is obtained from the first one 
by replacing $\si$ with $\si^{-1}$. 
Since $I(\si) = \{p_{i_n}\}$ by Lemma \ref{lem:EI2}, we have 
only to show that $\si^{-1}(p_{i_n}) = p_{i_n}$, 
namely, that the indeterminacy point $p_{i_n}$ of $\si$ is a 
fixed point of $\si^{-1} = \si_{i_n} \cdots \si_{i_2} \si_{i_1}$. 
By Lemma \ref{lem:L}, if two indices $i$ and $j$ are distinct, then 
the point $p_i$ lies on the line $L_j$ and hence is sent to $p_j$ 
by the map $\si_j$. 
Using this fact repeatedly, we see that 
\[
p_{i_n} \stackrel{\si_{i_1}}{\longmapsto} p_{i_1} 
        \stackrel{\si_{i_2}}{\longmapsto} p_{i_2} 
        \longmapsto \cdots \longmapsto p_{i_{n-1}} 
        \stackrel{\si_{i_n}}{\longmapsto} p_{i_n},  
\]
because every neighboring indices are distinct. 
Now it follows from formula (\ref{eqn:AS2}) that $\si$ satisfies 
condition (\ref{eqn:AS}) and hence is AS as desired. 
\hfill $\Box$ \par\medskip 
\begin{definition} \label{def:AS} 
We introduce two simple examples of AS transformations in $G$. 
\begin{enumerate}
\item An AS element $\si \in G$ is said to be {\sl elementary} 
if $\si = (\si_i \si_j)^m$ for some 
$\{i,j,k\} = \{1,2,3\}$ and $m \in \N$; otherwise, $\si$ is 
said to be {\sl non-elementary}. 
\item An element $\si \in G$ is called a {\sl Coxeter 
element} if $\si = \si_i \si_j \si_k$ for some 
$\{i,j,k\} = \{1,2,3\}$. 
\end{enumerate}
\end{definition}
\par
We may assume without loss of generality that $\si$ is AS, 
since if $\si$ is not AS then it can be replaced with its 
conjugate $\si' := \tau^{-1} \si \tau = \si_{i_{\nu+1}} 
\cdots \si_{i_{n-\nu}}$ which is AS, where 
$\tau = \si_{i_1}\cdots\si_{i_{\nu}}$ with $\nu$ being the 
index such that $i_1 = i_n, \, i_2 = i_{n-1}, \, \dots, \, 
i_{\nu} = i_{n-\nu+1}$ but $i_{\nu+1} \neq i_{n-\nu}$. 
Under this assumption we can apply formula (\ref{eqn:s-radius}) 
to conclude that the first dynamical degree of $\si$ is 
equal to the spectral radius of the linear map 
$\si^* : H^2(\ol{\Sol}(\th), \C) \carl$. 
On the other hand, Lemma \ref{lem:char} implies that the 
eigenvalues of $\si^*$ are $0$, $(-1)^n$ and the 
roots of the quadratic equation 
\begin{equation} \label{eqn:quadratic}
\l^2-\a(\si) \l + (-1)^n = 0, 
\end{equation}
so that the spectral radius of $\si^*$ is the largest absolute 
value of the roots of equation (\ref{eqn:quadratic}). 
This observation leads us to investigate the value distribution 
of $\a(\si)$. 
\begin{lemma} \label{lem:alpha} 
Assume that $\si \in G$ is AS. 
Then $\a(\si)$ is an even positive integer. 
Moreover, 
\begin{enumerate} 
\item $\a(\si) = 2$ if and only if $\si$ is elementary 
in the sense of Definition $\ref{def:AS}$, 
\item $\a(\si) = 4$ if and only if $\si$ is a Coxeter element, 
\item $\a(\si) = 6$ if and only if $\si = \si_i \si_j \si_k \si_j$ 
or $\si = \si_j \si_i \si_j \si_k$ for some $\{i,j,k\} = \{1,2,3\}$. 
\end{enumerate} 
\end{lemma} 
{\it Proof}.  
Let $\si = \si_{i_1} \si_{i_2} \cdots \si_{i_n}$ be the 
reduced expression of $\si$ as in (\ref{eqn:reduced3}). 
For $\nu = 1, \dots, n$, we put 
$A_{\nu} := s_{i_{\nu}} \cdots s_{i_2} s_{i_1}$ and denote 
its $(i,j)$-th entry by $(A_{\nu})_{ij}$. 
By definition (\ref{eqn:alpha2}) we have $\a(\si) = \Tr\, A_n$. 
We may assume that $i_1 = 1$, since the other cases can be 
treated in a similar manner. 
In this case, if we put $M_{\nu} := \min \{\, (A_{\nu})_{ij} 
\,:\, i = 1,2,3, \, j = 2,3 \,\}$, then 
\begin{equation} \label{eqn:order} 
M_{\nu+1} \ge M_{\nu} \qquad (\nu = 1,\dots,n-1). 
\end{equation} 
Moreover, if the index $j_{\nu+1}$ is defined by 
$\{j_{\nu+1}\} = \{1,2,3\}-\{i_{\nu}, i_{\nu+1}\}$, then  
\begin{equation} \label{eqn:recurrence} 
\Tr \, A_{\nu+1} = \Tr \, A_{\nu} + 
2(A_{\nu})_{j_{\nu+1}, i_{\nu+1}} 
\qquad  (\nu = 1,\dots,n-1).  
\end{equation} 
Indeed it is easy to see from formula 
(\ref{eqn:matrices}) that when $i_1 = 1$, the matrix 
$A_{\nu}$ takes the form 
\[
A_{\nu} = 
\begin{pmatrix}
0 & a_{12}^{\nu} & a_{13}^{\nu} \\
0 & a_{22}^{\nu} & a_{23}^{\nu} \\
0 & a_{32}^{\nu} & a_{33}^{\nu} 
\end{pmatrix}, 
\]
where $a_{ij}^{\nu}$, $i = 1,2,3$, $j = 2,3$, are nonnegative 
integers, and $A_{\nu+1} = s_{i_{\nu+1}} A_{\nu}$ is given by 
\[
\begin{array}{rcll}
A_{\nu+1} &=&  
\begin{pmatrix}
0 & a_{22}^{\nu} + a_{32}^{\nu} & a_{23}^{\nu} + a_{33}^{\nu} \\
0 & a_{22}^{\nu} & a_{23}^{\nu} \\
0 & a_{32}^{\nu} & a_{33}^{\nu} 
\end{pmatrix} \qquad & (\mbox{if} \,\, i_{\nu+1} = 1), \\[7mm]
A_{\nu+1} &=& 
\begin{pmatrix}
0 & a_{12}^{\nu} & a_{13}^{\nu} \\
0 & a_{12}^{\nu} + a_{32}^{\nu} & a_{13}^{\nu} + a_{33}^{\nu} \\
0 & a_{32}^{\nu} & a_{33}^{\nu} 
\end{pmatrix} \qquad & (\mbox{if} \,\, i_{\nu+1} = 2), \\[7mm]
A_{\nu+1} &=& 
\begin{pmatrix}
0 & a_{12}^{\nu} & a_{13}^{\nu} \\
0 & a_{22}^{\nu} & a_{23}^{\nu} \\
0 & a_{12}^{\nu} + a_{22}^{\nu} & a_{13}^{\nu} + a_{23}^{\nu} 
\end{pmatrix} \qquad  & (\mbox{if} \,\, i_{\nu+1} = 3). 
\end{array}
\] 
Inequality (\ref{eqn:order}) readily follows from these 
observations and formula (\ref{eqn:recurrence}) is 
verified by a case-by-case check. 
Indeed, if $i_{\nu} = 1$ and $i_{\nu+1} = 2$, then  
$j_{\nu+1} = 3$ and $a_{12}^{\nu} = 
a_{22}^{\nu} + a_{32}^{\nu}$ so that 
\[ 
\Tr \, A_{\nu+1} = 
a_{12}^{\nu} + a_{32}^{\nu}+ a_{33}^{\nu} = 
(a_{22}^{\nu} + a_{32}^{\nu}) + a_{32}^{\nu}+ a_{33}^{\nu} 
= a_{22}^{\nu} + a_{33}^{\nu} + 2 a_{32}^{\nu} = 
\Tr\, A_{\nu} + 2 (A_{\nu})_{j_{\nu+1}, i_{\nu+1}}. 
\] 
If $i_{\nu} = 2$ and $i_{\nu+1} = 1$, then 
$j_{\nu+1} = 3$ and $(A_{\nu})_{j_{\nu+1}, i_{\nu+1}} = 0$ 
so that 
\[
\Tr \, A_{\nu+1} = a_{22}^{\nu} + a_{33}^{\nu} = 
\Tr \, A_{\nu} = \Tr \, A_{\nu} + 
2 (A_{\nu})_{j_{\nu+1}, i_{\nu+1}}. 
\]
The remaining cases can be treated in similar manners. 
Note that (\ref{eqn:recurrence}) yields an inequality 
$\Tr\, A_{\nu+1} \ge \Tr\, A_{\nu}$, since 
$(A_{\nu})_{j_{\nu+1}, i_{\nu+1}}$ is nonnegative. 
A repeated use of formula (\ref{eqn:recurrence}) shows 
that $\a(\si) = \Tr \, A_n$ is an even integer not 
smaller than $2$, because $\Tr\, A_1 = \Tr\, s_1 = 2$. 
\par 
Next we observe that $\a(\si) = 2$ if 
$\si = (\si_i \si_j)^m$ for any $\{i,j,k\} = \{1,2,3\}$ and 
$m \in \N$. 
Indeed, since we are assuming that $i_1 = 1$, we have only 
to check the two cases where $A_n = (s_2 s_1)^m$ and 
$A_n = (s_3 s_1)^m$ with $m \in \N$. 
In either case we have $\a(\si) = \Tr\, A_n = 2$ because 
\[
(s_2 s_1)^m = 
\begin{pmatrix} 
0 & 1 & 2m-1 \\ 0 & 1 & 2m \\ 0 & 0 & 1 
\end{pmatrix}, 
\qquad 
(s_3 s_1)^m = 
\begin{pmatrix} 
0 & 2m-1 & 1 \\ 0 & 1 & 0 \\ 0 & 2m & 1 
\end{pmatrix}. 
\]
\par 
From now on we assume that $\si$ is not of the form 
$(\si_i \si_j)^m$ for any $\{i,j,k\} = \{1,2,3\}$ and 
$m \in \N$. 
Then the length $n$ must be not less than $3$ and there 
exists an index $\nu$ such that 
$\{i_{\nu},i_{\nu+1},i_{\nu+2}\} = \{1,2,3\}$. 
Here we may assume without loss of generality that $\nu = 1$, 
since the quantity $\a(\si)$ is invariant under any cyclic 
permutation of the indices $(i_1,\dots,i_n)$, provided that 
$\si$ is an AS element. 
Since moreover we are assuming that $i_1 = 1$, we have 
\begin{equation} \label{eqn:A3}
A_3 = s_3 s_2 s_1 = 
\begin{pmatrix} 
0 & 1 & 1 \\ 0 & 1 & 2 \\ 0 & 2 & 3 
\end{pmatrix}, 
\qquad \mbox{or} \qquad 
A_3 = s_2 s_3 s_1 = 
\begin{pmatrix} 
0 & 1 & 1 \\ 0 & 3 & 2 \\ 0 & 2 & 1 
\end{pmatrix}. 
\end{equation}
\par 
If $n = 3$ we have $\si = \si_1 \si_2 \si_3$ or 
$\si = \si_1 \si_3 \si_2$. 
Then formula (\ref{eqn:A3}) yields 
$\a(\si) = \Tr\, A_3 = 4$ in either case. 
If $n = 4$ we have $\si = \si_1 \si_2 \si_3 \si_2$ 
or $\si = \si_1 \si_3 \si_2 \si_3$. 
Since 
\begin{equation} \label{eqn:A4}
A_4 = s_2 s_3 s_2 s_1 = 
\begin{pmatrix} 
0 & 1 & 1 \\ 0 & 3 & 4 \\ 0 & 2 & 3 
\end{pmatrix}, 
\qquad \mbox{or} \qquad 
A_4 = s_3 s_2 s_3 s_1 = 
\begin{pmatrix} 
0 & 1 & 1 \\ 0 & 3 & 2 \\ 0 & 4 & 3 
\end{pmatrix},  
\end{equation}
we have $\a(\si) = \Tr\, A_4 = 6$ in either case. 
Finally we assume that $n \ge 5$. 
Then $A_4$ is either (\ref{eqn:A4}) or 
$A_4 = s_1 s_j s_k s_1$ for some $\{j,k\} = \{2,3\}$. 
In the latter case we must have 
$A_5 = s_j s_1 s_j s_k s_1$ or 
$A_5 = s_k s_1 s_j s_k s_1$. 
Here we can eliminate the last term $s_1$ by taking a 
cyclic permutation of the indices $(i_1,\dots,i_n)$ and 
obtain $A_4 = s_j s_1 s_j s_k$ or $A_4 = s_k s_1 s_j s_k$. 
By relabeling the indices, the matrix $A_4$ can be reduced 
to the form (\ref{eqn:A4}). 
So we have only to consider the former case (\ref{eqn:A4}). 
Since $\si$ is assumed to be AS, the index $i_n$ is 
different from $i_1 = 1$ so that $(A_{n-1})_{j_n, i_n} 
\ge M_{n-1}$ by the definition of $M_{\nu}$. 
Then it follows from (\ref{eqn:order}) and 
(\ref{eqn:recurrence}) that $\a(\si) = \Tr\, A_n$ is 
estimated as 
\[
\a(\si) = \Tr\, A_{n-1} + 2 (A_{n-1})_{j_n, i_n} 
\ge \Tr\, A_{n-1} + 2 M_{n-1} 
\ge \Tr\, A_4 + 2 M_4 
= 6 + 2 \times 1 = 8. 
\]
Putting all these arguments together we establish 
the lemma. \hfill $\Box$ \par\medskip 
\begin{lemma} \label{lem:fdd} 
If $\si \in G$ is AS then the first dynamical degree of 
$\si$ is given by 
\begin{equation} \label{eqn:fdd}
\l_1(\si) = \dfrac{1}{2} \left\{\a(\si) + 
\sqrt{\a(\si)^2 + 4(-1)^{n+1}}\right\}, 
\end{equation}
where $n = \ell_G(\si)$ is the length of the element $\si$. 
Moreover, 
\begin{enumerate}
\item if $\si$ is elementary then $\l_1(\si) = 1$, 
\item if $\si$ is a Coxeter element then 
$\l_1(\si) = 2 + \sqrt{5}$, 
\item if $\si = \si_i \si_j \si_k \si_j$ or 
$\si = \si_j \si_i \si_j \si_k$ for some 
$\{i,j,k\} = \{1,2,3\}$, then $\l_1(\si) = 3 + 2 \sqrt{2}$, 
\item otherwise, we have $\l_1(\si) \ge 4 + \sqrt{15}$. 
\end{enumerate}
\end{lemma} 
{\it Proof}. 
Since $\a(\si) \ge 2$ by Lemma \ref{lem:alpha}, the 
quadratic equation (\ref{eqn:quadratic}) has the real 
roots 
\begin{equation} \label{eqn:roots}
\l_{\pm}(\si) = \dfrac{1}{2} 
\left\{\a(\si) \pm \sqrt{\a(\si)^2+4(-1)^{n+1}} \right\}, 
\end{equation} 
where $\l_+(\si) \ge 1$ and 
$|\l_-(\si)| = \l_+(\si)^{-1} \le 1$. 
Therefore the root $\l_+(\si)$ gives the spectral radius 
$\mathrm{SR}(\si^*)$ of $\si^*$ and hence the first dynamical 
degree $\l_1(\si)$ 
of $\si$ by formula (\ref{eqn:s-radius}). 
Assertions (1), (2), (3) can be checked directly by using 
Lemma \ref{lem:alpha}. 
Finally we shall show assertion (4). 
In this case, since $\a(\si) \ge 8$ by Lemma \ref{lem:alpha}, 
formula (\ref{eqn:fdd}) implies that 
\[
\l_1(\si) \ge \dfrac{1}{2} 
\left\{\a(\si) + \sqrt{\a(\si)^2-4} \right\} 
\ge \dfrac{1}{2}\left(8 + \sqrt{8^2-4}\right) = 4 + \sqrt{15}. 
\]
Hence the lemma is proved. \hfill $\Box$ \par\medskip 
\par 
We proceed to the construction of a natural $\si$-invariant 
measure for an AS element $\si \in G$. 
Again let us start with the general situation where 
$f : S \carl$ is an AS bimeromorphic map on a compact 
K\"ahler surface $S$. 
If $\l_1(f) = 1$ then either $f$ is a dynamically trivial 
automorphism or $f$ preserves a rational or elliptic 
fibration and exhibits an essentially $1$-dimensional 
dynamic \cite{DF}. 
In our case where $f = \si$ and $S = \ol{\Sol}(\th)$, 
the condition $\l_1(\si) = 1$ means that $\si$ is elementary 
by Lemma \ref{lem:fdd}. 
If so, the existence of a $\si$-invariant rational fibration 
on $\ol{\Sol}(\th)$ can be seen easily 
(see Remark \ref{rem:proof}). 
So we are not interested in the case $\l_1(f) = 1$ and 
assume hereafter that 
\begin{equation} \label{eqn:fdd2}
\l_1(f) > 1. 
\end{equation}
In this case it is known \cite{DF} that there are positive 
closed $(1,1)$-currents $\mu^{\pm}$ on $S$ such that 
\[
(f^{\pm1})^* \mu^{\pm} = \l_1(f) \, \mu^{\pm}, 
\]
where $\mu^+$ and $\mu^-$ are called the stable and unstable 
currents for $f$. 
A natural strategy to obtain an $f$-invariant measure $\mu$ on 
$S$ is to take the wedge product 
\begin{equation} \label{eqn:invmeas}
\mu = \mu^+ \wedge \mu^-. 
\end{equation}
However the main issue here is whether the operation of 
wedge product is feasible or not. 
If the stable and unstable currents are expressed as 
$\mu^{\pm} = dd^c g^{\pm}$ in terms of local potentials 
$g^{\pm}$, then the wedge product (\ref{eqn:invmeas}) 
may be interpreted as the complex Monge-Amp\`{e}re operator 
$dd^c g^+ \wedge dd^c g^-$. 
In order for this operation to be well-defined, a quantitative 
condition 
\begin{equation} \label{eqn:AS3} 
\sum_{N = 0}^{\infty} \l_1(f)^{-N} 
\log \, \mathrm{dist}(f^N I(f^{-1}), I(f)) > - \infty, 
\end{equation} 
is introduced in \cite{BD}, where $\mathrm{dist}$ is the distance 
on $S$ induced from a Riemannian metric on it. 
This condition is slightly stronger than (\ref{eqn:AS}) and 
a map enjoying this condition might be called quantitatively AS. 
Under these settings the following theorem is established 
in \cite{BD}. 
\begin{theorem} \label{thm:BD} 
If $f : S \carl$ satisfies conditions $(\ref{eqn:fdd2})$ and 
$(\ref{eqn:AS3})$, then the wedge product $\mu$ of the stable 
and unstable currents $\mu^{\pm}$ in $(\ref{eqn:invmeas})$ is 
well defined and, after a suitable renormalization, $\mu$ gives 
an $f$-invariant Borel probability measure such that all 
the conditions in Definition $\ref{def:chaos}$ are satisfied. 
Moreover the measure $\mu$ puts no mass on any algebraic 
curve on $S$. 
\end{theorem} 
\par 
Applying this theorem to our situation, we obtain the 
following theorem. 
\begin{theorem} \label{thm:chaos2} 
For any non-elementary AS map $\si \in G$ there 
exists the wedge product 
$\mu_{\si} = \mu_{\si}^+ \wedge \mu_{\si}^-$ of 
the stable and unstable currents $\mu_{\si}^{\pm}$ for 
$\si$ and, after a suitable renormalization, $\mu_{\si}$ 
gives a $\si$-invariant Borel probability measure such that 
all the conditions in Definition $\ref{def:chaos}$ are satisfied. 
Moreover the measure $\mu_{\si}$ puts no mass on any algebraic 
curve on $\ol{\Sol}(\th)$. 
\end{theorem}
{\it Proof}. It is enough to check that any non-elementary AS 
map $\si \in G$ satisfies conditions $(\ref{eqn:fdd2})$ 
and $(\ref{eqn:AS3})$. 
Lemma \ref{lem:fdd} implies that $\l_1(\si) > 1$ 
if and only if $\si$ is non-elementary, so that condition 
(\ref{eqn:fdd2}) is satisfied. 
In order to check condition (\ref{eqn:AS3}) 
let $\si = \si_{i_1} \cdots \si_{i_n}$ be the reduced expression 
of $\si$. 
Since $\si$ is assumed to be AS, the indices $i_1$ and $i_n$ are 
distinct and hence $\mathrm{dist}(p_{i_1}, p_{i_n}) > 0$. 
On the other hand, by formula (\ref{eqn:AS2}), we have  
$I(\si) = \{p_{i_n}\}$ and $\si^N I(\si^{-1}) = \{p_{i_1}\}$ 
independently of $N \ge 0$. 
Therefore we have  
\begin{eqnarray*}
\sum_{N=0}^{\infty} \l_1(\si)^{-N} \log\,
\mathrm{dist}(\si^N I(\si^{-1}), I(\si)) 
&=& \log\, \mathrm{dist}(p_{i_1}, p_{i_n}) \sum_{N=0}^{\infty} 
\l_1(\si)^{-N} \\ 
&=&  \dfrac{\l_1(\si) \log\, 
\mathrm{dist}(p_{i_1}, p_{i_n})}{\l_1(\si)-1} > - \infty, 
\end{eqnarray*}
which shows that condition (\ref{eqn:AS3}) is satisfied. 
The theorem then follows from Theorem \ref{thm:BD}. 
\hfill $\Box$  
\begin{remark} \label{rem:lyapunov2} 
Under the setting of Theorem \ref{thm:BD} it is shown in 
\cite{BD} that the Lyapunov exponents $L_{\pm}(f)$ of $f$ 
with respect to the ergodic measure $\mu$ satisfy the 
estimate 
\[
L_-(f) \le - \dfrac{\log \l_1(f)}{8} < 0 < 
\dfrac{\log \l_1(f)}{8} \le L_+(f), 
\]
which applies to the mapping $\si : \Sol(\th) \carl$ 
in Theorem \ref{thm:chaos2}. 
On the other hand, we have $L_-(\si) = - L_+(\si)$ since 
$\si$ is area-preserving with respect to the Poincar\'e 
residue $\omega(\th)$ in (\ref{eqn:poincare}). 
It follows from the above estimate that 
$L_+(\si) \ge \frac{1}{8} \log \l_1(\si)$. 
\end{remark}
\par
Finally we shall calculate the entropy of a non-elementary AS 
map $\si \in G$. 
For a birational map $f : S \carl$ of a projective 
surface $S$ and an $f$-invariant Borel probability measure 
$\mu$ on $S$, there are two concepts of entropies: one is the 
{\sl measure-theoretic entropy} $h_{\mu}(f)$ with respect to 
the invariant measure $\mu$ and the other is the {\sl topological 
entropy} $h_{\mathrm{top}}(f)$. 
In general these quantities and the first dynamical degree 
$\l_1(f)$ are related as 
\begin{equation} \label{eqn:inequalities}
h_{\mu}(f) \le h_{\mathrm{top}}(f) \le \log \l_1(f), 
\end{equation}
where the first inequality is the so-called variational 
principle and the second inequality is a consequence of a 
main result of \cite{DS}. 
Moreover, if $f$ satisfies conditions (\ref{eqn:fdd2}) and 
(\ref{eqn:AS3}) and if $\mu$ is the invariant measure mentioned 
in Theorem \ref{thm:BD}, then it is proved in \cite{Dujardin} 
that the leftmost and rightmost terms in 
(\ref{eqn:inequalities}) are equal and consequently all 
the three terms in (\ref{eqn:inequalities}) coincide. 
Applying this triple coincidence to our situation we obtain 
the following theorem. 
\begin{theorem} \label{thm:entropy} 
For any non-elementary AS map $\si \in G$, we have  
\begin{equation} \label{eqn:trinity}
h_{\mu_{\si}}(\si) = h_{\mathrm{top}}(\si) = \log \l_1(\si).  
\end{equation}
where $\mu_{\si}$ is the $\si$-invariant probability measure 
mentioned in Theorem $\ref{thm:chaos2}$. 
The value of $(\ref{eqn:trinity})$ is not smaller than 
$\log(2 + \sqrt{5})$ with equality if and only if 
$\si$ is a Coxeter element. 
\end{theorem} 
{\it Proof}. 
The proof is already finished in the above argument. 
The assertion that (\ref{eqn:trinity}) takes its minimum 
precisely when $\si$ is a Coxeter element follows from 
Lemma \ref{lem:fdd}. 
\hfill $\Box$ 
\begin{remark} \label{rem:affine} 
Theorems \ref{thm:chaos2} and \ref{thm:entropy} are results 
for an element $\si \in G$ viewed as a birational map of the 
projective surface $\ol{\Sol}(\th)$. 
However, since the invariant measure $\mu_{\si}$ put no mass 
on any algebraic curve on $\ol{\Sol}(\th)$, the lines 
$L = L_1 \cup L_2 \cup L_3$ at infinity can be neglected as far 
as the ergodic properties of $\si : \ol{\Sol}(\th) \carl$ 
relative to the measure $\mu_{\si}$ are concerned. 
So Theorems \ref{thm:chaos2} and \ref{thm:entropy} lead 
to results for the biregular map $\si' := \si|_{\Sol(\th)}$ 
of the affine surface $\Sol(\th) = \ol{\Sol}(\th) - L$. 
Namely $\mu_{\si}$ can be restricted without losing any mass 
to an $\si'$-invariant Borel probability measure $\mu_{\si'}$ 
on $\Sol(\th)$ such that the conditions in 
Definition \ref{def:chaos} are satisfied, and one has an equality 
$h_{\mu_{\si'}}(\si') = \log \l_1(\si)$. 
Here we do not refer to $h_{\mathrm{top}}(\si')$, because 
the concept of topological entropy, usually 
defined on a compact space, is not very clear on the 
affine surface $\Sol(\th)$. 
In what follows $\si'$ and $\mu_{\si'}$ will be written 
$\si$ and $\mu_{\si}$ for the simplicity of notation. 
\end{remark} 
\section{Number of Periodic Points} \label{sec:periodic} 
Given any non-elementary AS element $\si \in G$, we are 
interested in the number of periodic points of the birational 
map $\si : \ol{\Sol}(\th) \carl$. 
For each positive integer $N \in \N$ we shall consider the set 
of all periodic points of period $N$ on the projective cubic 
surface $\ol{\Sol}(\th)$, 
\[
\ol{\Per}_N(\si;\th) := 
\{\, X \in \ol{\Sol}(\th) - I(\si^N) \,:\, \si^N(X) = X \, \}, 
\]
as well as the corresponding set on the affine cubic surface 
$\Sol(\th)$, 
\[
\Per_N(\si;\th) := \{\, x \in \Sol(\th) \,:\, \si^N(x) = x \, \}. 
\]
Our tasks are then to count the cardinality of 
$\ol{\Per}_N(\si;\th)$ and to relate it with the cardinality 
of $\Per_N(\si;\th)$. 
The first task is based on the Lefschetz fixed point formula, 
while the second one is by a careful inspection of the 
behavior of the map $\si$ around the lines $L 
= L_1 \cup L_2 \cup L_3$ at infinity. 
In order to apply the Lefschetz fixed point formula, 
we need the following lemma. 
\begin{lemma} \label{lem:percurve} 
Assume that $\si \in G$ is AS and non-elementary. 
Then for any $N \in \N$, the birational map 
$\si : \ol{\Sol}(\th) \carl$ admits no curves of periodic 
points of period $N$. 
\end{lemma}
{\it Proof}. 
The lemma is proved by contradiction. 
Assume that $\si$ admits a curve (an effective divisor) 
$D \subset \ol{\Sol}(\th)$ of periodic points of some 
period $N$. 
Since $\si^N$ fixes $D$ pointwise, we have 
$(\si^N)^* D = D$ in $H^2(\ol{\Sol}(\th),\Z)$. 
Moreover, since $\si$ is assumed to be AS, we have 
$(\si^*)^N = (\si^N)^*$ and hence $(\si^*)^N D = D$, which 
means that $(\si^*)^N$ has an eigenvalue $1$ with an 
eigenvector $D$. 
On the other hand, by Lemma \ref{lem:char}, there is a direct 
sum decomposition $\si^* = (\si^*|_V) \oplus (\si^*|_{V^{\perp}})$ 
with $H^2(\ol{\Sol}(\th),\Z) = V \oplus V^{\perp}$ as in 
(\ref{eqn:decomp}) such that $\si^*|_V$ has the eigenvalues 
$0$ and $\l_{\pm}(\si)$ as in (\ref{eqn:roots}), 
while $\si^*|_{V^{\perp}}$ is the scalar operator $(-1)^n$ on 
$V^{\perp}$, where $n = \ell_G(\si)$ is the length of $\si$. 
By Lemma \ref{lem:fdd} we have $\l_+(\si) = \l_1(\si)> 1$ and 
$|\l_-(\si)| = \l_1(\si)^{-1} < 1$, since $\si$ is assumed to 
be non-elementary. 
Therefore the eigenvector $D$ of $\si^*$ must belong to 
the subspace $V^{\perp}$ and its eigenvalue $1$  must arise as 
the $N$-th power of the scalar operator 
$(-1)^n = \si^*|_{V^{\perp}}$, where the integer 
$n N$ must be even. 
Since $L_1, \, L_2, \, L_3 \in V$ and $D \in V^{\perp}$, 
we have 
\begin{equation} \label{sec:DLi}
(D, L_1) = (D, L_2) = (D, L_3) =0. 
\end{equation} 
We now write $D = D' + m_1 L_1 + m_2 L_2 + m_3 L_3$, 
where $D'$ is either empty or an effective divisor 
not containing $L_1$, $L_2$, $L_3$ as an irreducible 
component of it and $m_1$, $m_2$, $m_3$ are nonnegative 
integers. 
Since $(L_i, L_j) = -1$ for $i = j$ and 
$(L_i, L_j) = 1$ for $i \neq j$, 
formula (\ref{sec:DLi}) yields
\[
\begin{array}{rclcl}
0 &=& (D, L_1) &=& (D', L_1)-m_1+m_2+m_3, \\[1mm]
0 &=& (D, L_2) &=& (D', L_2)+m_1-m_2+m_3, \\[1mm]
0 &=& (D, L_3) &=& (D', L_3)+m_1+m_2-m_3, 
\end{array}
\]
which sum up to 
\begin{equation} \label{eqn:sum}
(D', L_1)+(D', L_2)+(D', L_3)+m_1+m_2+m_3 = 0.
\end{equation} 
Since none of the lines $L_1$, $L_2$, $L_3$ is an irreducible 
component of $D'$, the intersection number 
$(D', L_i)$ must be nonnegative for every $i = 1, 2, 3$. 
Since the numbers $m_1$, $m_2$, $m_3$ are also nonnegative, 
formula (\ref{eqn:sum}) implies that 
$(D',L_1) = (D',L_2)= (D',L_3) =0$ 
and $m_1 = m_2 = m_3 = 0$. 
Hence $D = D'$ and $(D,L_1) = (D,L_2) = (D,L_3) =0$. 
It follows that $D$ is an effective divisor such that 
$(D, L_i) = 0$ and $L_i$ is not an irreducible component 
of $D$ for every $i = 1,2,3$. 
This means that the compact curve $D$ does not intersect 
$L = L_1 \cup L_2 \cup L_3$ and hence must lie in 
the affine cubic surface $\Sol(\th) = \ol{\Sol}(\th)-L$. 
However no compact curve can lie in any affine variety. 
By this contradiction the lemma is established. 
\hfill $\Box$ \par\medskip 
Now we shall apply the Lefschetz fixed point formula to 
the iterates of a non-elementary AS element $\si \in G$. 
For each $N \in \Z$ let $\vG_N \subset 
\ol{\Sol}(\th) \times \ol{\Sol}(\th)$ be the graph of the 
$N$-th iterate $\si^N : \ol{\Sol}(\th) \carl$, and 
$\vD \subset \ol{\Sol}(\th) \times \ol{\Sol}(\th)$ be 
the diagonal. 
Note that $\vG_N = \vG_{-N}^{\vee}$, where 
$\vG_{-N}^{\vee}$ is the reflection of 
$\vG_{-N}$ around the diagonal $\vD$. 
Moreover let $I_N \subset \ol{\Sol}(\th)$ denote the 
indeterminacy set of $\si^N$. 
Then the Lefschetz fixed point formula consists of 
two equations concerning the intersection number 
$(\vG_N, \vD)$ of the cycles $\vG_N$ and $\vD$ in 
$\ol{\Sol}(\th) \times \ol{\Sol}(\th)$,  
\begin{eqnarray} 
(\vG_N, \vD) &=& \sum_{q=0}^4 (-1)^q \, \Tr\,[\, 
(\si^N)^* : H^q(\ol{\Sol}(\th),\Z) \carl\,], 
\label{eqn:lfpf1} \\
(\vG_N, \vD) &=& \mathrm{\#} \, \ol{\Per}_N(\si;\th) 
+ \sum_{p \in I_N} \mu((p,p), \vG_N \cap \vD), 
\label{eqn:lfpf2} 
\end{eqnarray} 
where $\mu((p,p), \vG_N \cap \vD)$ denotes the multiplicity 
of intersection between $\vG_N$ and $\vD$ at $(p,p)$. 
Lemma \ref{lem:percurve} assures that all terms involved in 
(\ref{eqn:lfpf1}) and (\ref{eqn:lfpf2}) are well defined and 
finite. 
\begin{lemma} \label{lem:lfpf1} 
Let $n = \ell_G(\si)$ be the length of $\si$. 
Then formula $(\ref{eqn:lfpf1})$ becomes 
\[
(\vG_N, \vD) = \l_1(\si)^N + (-1)^{nN}\l_1(\si)^{-N} + 
4 (-1)^{nN} + 2. 
\]
\end{lemma} 
{\it Proof}. 
We put 
$T_N^q = \Tr\,[\, (\si^N)^* : H^q(\ol{\Sol}(\th),\Z) \carl\,]$. 
Because $\ol{\Sol}(\th)$ is a smooth rational surface, 
\[
H^q(\ol{\Sol}(\th), \Z) \cong \left\{ 
\begin{array}{cl}
\Z \qquad & (q = 0,4), \\[1mm] 
0  \qquad & (q = 1,3). 
\end{array}
\right. 
\]
Trivially we have $T_N^0 = 1$ and $T_N^1 = T_N^3 = 0$. 
Since $\si$ and so $\si^N$ are birational, we have $T_N^4 = 1$. 
Since the map $\si$ is assumed to be AS, we have 
$(\si^N)^* = (\si^*)^N : H^2(\ol{\Sol}(\th),\Z) \carl$. 
By Lemmas \ref{lem:char} and \ref{lem:fdd}, 
$\si^*$ has three simple eigenvalues $0$, 
$\l_+(\si) = \l_1(\si)$, $\l_-(\si) = (-1)^n \l_1(\si)^{-1}$ 
and a quadruple eigenvalue $(-1)^n$. 
Thus we have $T_N^2 = 0^N + \l_1(\si)^N + 
(-1)^{nN} \l_1(\si)^{-N} + 4(-1)^{nN}$. 
Substituting these data into (\ref{eqn:lfpf1}) yields 
the assertion of the lemma. \hfill $\Box$ \par\medskip 
\begin{lemma} \label{lem:lfpf2} 
Formula $(\ref{eqn:lfpf2})$ becomes 
\[
(\vG_N, \vD) = \mathrm{\#} \, \ol{\Per}_N(\si;\th) + 1 
= \mathrm{\#} \, \Per_N(\si;\th) + 2. 
\] 
\end{lemma} 
{\it Proof}. 
Let $\si = \si_{i_1} \cdots \si_{i_n}$ be the reduced 
expression of $\si$. 
Since $\si$ is assumed to be AS, for any $N \in \N$ 
the reduced expression of $\si^N$ is given by 
$\si^N = \overbrace{\si_{i_1} \cdots \si_{i_n}} \cdots 
\overbrace{\si_{i_1} \cdots \si_{i_n}}$ ($N$-times). 
Moreover, since $\si$ is assumed to be non-elementary, 
the indices $\{i_1, \dots, i_n\}$ range the entire 
index set $\{1,2,3\}$. 
By Lemma \ref{lem:EI2} the exceptional set of $\si^N$ is 
given by 
\[
\E(\si^N) = \bigcup_{\nu=1}^n L_{i_{\nu}} 
= L_1 \cup L_2 \cup L_3 = L, 
\]
whose $\si^N$-image is 
$\si^N(L) = \si^N(\E(\si^N)) = \{p_{i_1}\}$. 
This means that $p_{i_1}$ is the unique fixed point of 
the map $\si^N$ on the lines $L$ at infinity. 
Lemma \ref{lem:EI2} also implies that $p_{i_n}$ is the 
unique indeterminacy point of $\si^N$. 
Therefore we have 
$\ol{\Per}_N(\si;\th) = \Per_N(\si;\th) \cup \{p_{i_1}\}$ 
and $I_N = \{p_{i_n}\}$, which implies that formula 
(\ref{eqn:lfpf2}) is rewritten as 
\begin{equation} \label{eqn:lfpf3}
\begin{array}{rcl} 
(\vG_N, \vD) &=& \mathrm{\#} \, \ol{\Per}_N(\si;\th) 
+ \mu((p_{i_n},p_{i_n}), \vG_N \cap \vD),   \\[2mm]  
\mathrm{\#} \, \ol{\Per}_N(\si;\th) &=& 
\mathrm{\#} \, \Per_N(\si;\th) + \nu(p_{i_1}, \si^N), 
\end{array}
\end{equation} 
where $\nu(p_{i_1}, \si^N)$ is the local index of the map 
$\si^N$ around the fixed point $p_{i_1}$. 
If $j$ and $k$ are defined by $\{j, k\} = \{1,2,3\}-\{i_1\}$, 
then $L_j$ and $L_k$ are linearly independent lines passing 
through the point $p_{i_1}$. 
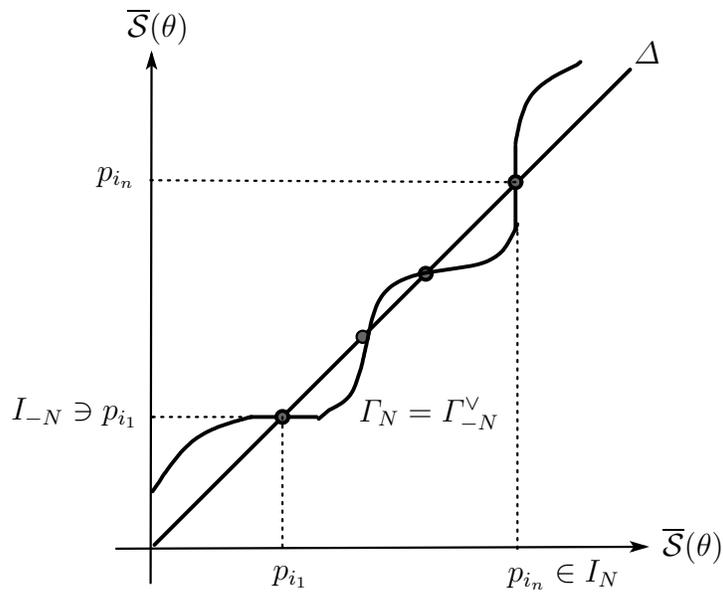
\begin{figure}[t]
\begin{center}
\unitlength 0.1in
\begin{picture}(35.70,29.90)(-0.20,-30.80)
%
\special{pn 13}%
\special{pa 690 2900}%
\special{pa 3460 2890}%
\special{fp}%
\special{sh 1}%
\special{pa 3460 2890}%
\special{pa 3393 2870}%
\special{pa 3407 2890}%
\special{pa 3393 2910}%
\special{pa 3460 2890}%
\special{fp}%
%
\special{pn 13}%
\special{pa 870 3080}%
\special{pa 876 310}%
\special{fp}%
\special{sh 1}%
\special{pa 876 310}%
\special{pa 856 377}%
\special{pa 876 363}%
\special{pa 896 377}%
\special{pa 876 310}%
\special{fp}%
%
\special{pn 20}%
\special{pa 890 2880}%
\special{pa 3380 390}%
\special{fp}%
\put(35.5000,-29.8000){\makebox(0,0)[lb]{$\ol{\Sol}(\th)$}}%
\put(7.5000,-2.6000){\makebox(0,0)[lb]{$\ol{\Sol}(\th)$}}%
\put(34.0000,-3.6000){\makebox(0,0)[lb]{$\vD$}}%
\put(27.4000,-29.6000){\makebox(0,0)[lt]{$p_{i_n} \in I_N$}}%
%
\special{pn 20}%
\special{pa 2780 790}%
\special{pa 2780 1230}%
\special{fp}%
\special{pa 2780 1230}%
\special{pa 2780 1230}%
\special{fp}%
%
\special{pn 20}%
\special{sh 0.600}%
\special{ar 2780 980 32 32  0.0000000 6.2831853}%
%
\special{pn 20}%
\special{pa 1390 2210}%
\special{pa 1740 2210}%
\special{fp}%
%
\special{pn 20}%
\special{sh 0.600}%
\special{ar 1560 2210 32 32  0.0000000 6.2831853}%
%
\special{pn 13}%
\special{pa 2790 1200}%
\special{pa 2790 2900}%
\special{dt 0.045}%
\special{pa 2790 2900}%
\special{pa 2790 2899}%
\special{dt 0.045}%
%
\special{pn 13}%
\special{pa 1390 2210}%
\special{pa 880 2210}%
\special{dt 0.045}%
\special{pa 880 2210}%
\special{pa 881 2210}%
\special{dt 0.045}%
%
\special{pn 20}%
\special{sh 0.600}%
\special{ar 2310 1460 32 32  0.0000000 6.2831853}%
%
\special{pn 13}%
\special{sh 0.600}%
\special{ar 1980 1790 32 32  0.0000000 6.2831853}%
%
\special{pn 20}%
\special{pa 2790 1200}%
\special{pa 2776 1230}%
\special{pa 2761 1260}%
\special{pa 2745 1289}%
\special{pa 2728 1315}%
\special{pa 2710 1340}%
\special{pa 2688 1361}%
\special{pa 2665 1379}%
\special{pa 2640 1394}%
\special{pa 2612 1407}%
\special{pa 2583 1417}%
\special{pa 2553 1426}%
\special{pa 2521 1434}%
\special{pa 2489 1440}%
\special{pa 2455 1445}%
\special{pa 2421 1450}%
\special{pa 2387 1455}%
\special{pa 2352 1461}%
\special{pa 2318 1467}%
\special{pa 2284 1473}%
\special{pa 2250 1482}%
\special{pa 2218 1491}%
\special{pa 2186 1503}%
\special{pa 2156 1516}%
\special{pa 2129 1531}%
\special{pa 2103 1549}%
\special{pa 2080 1569}%
\special{pa 2060 1592}%
\special{pa 2043 1618}%
\special{pa 2028 1646}%
\special{pa 2015 1676}%
\special{pa 2005 1706}%
\special{pa 1995 1738}%
\special{pa 1987 1771}%
\special{pa 1979 1803}%
\special{pa 1972 1835}%
\special{pa 1965 1866}%
\special{pa 1958 1897}%
\special{pa 1950 1928}%
\special{pa 1942 1959}%
\special{pa 1934 1989}%
\special{pa 1924 2020}%
\special{pa 1914 2050}%
\special{pa 1902 2080}%
\special{pa 1888 2110}%
\special{pa 1871 2137}%
\special{pa 1850 2160}%
\special{pa 1825 2180}%
\special{pa 1798 2196}%
\special{pa 1769 2211}%
\special{pa 1750 2220}%
\special{sp}%
%
\special{pn 20}%
\special{pa 1400 2210}%
\special{pa 1368 2219}%
\special{pa 1337 2228}%
\special{pa 1305 2237}%
\special{pa 1274 2248}%
\special{pa 1244 2258}%
\special{pa 1214 2270}%
\special{pa 1185 2283}%
\special{pa 1157 2297}%
\special{pa 1130 2313}%
\special{pa 1105 2331}%
\special{pa 1080 2350}%
\special{pa 1057 2371}%
\special{pa 1035 2393}%
\special{pa 1013 2416}%
\special{pa 993 2441}%
\special{pa 973 2466}%
\special{pa 953 2493}%
\special{pa 934 2519}%
\special{pa 915 2547}%
\special{pa 897 2574}%
\special{pa 880 2600}%
\special{sp}%
%
\special{pn 20}%
\special{pa 2780 780}%
\special{pa 2787 747}%
\special{pa 2794 714}%
\special{pa 2801 682}%
\special{pa 2810 650}%
\special{pa 2820 619}%
\special{pa 2831 590}%
\special{pa 2844 561}%
\special{pa 2860 535}%
\special{pa 2877 510}%
\special{pa 2898 488}%
\special{pa 2920 467}%
\special{pa 2945 447}%
\special{pa 2971 429}%
\special{pa 2999 412}%
\special{pa 3028 396}%
\special{pa 3058 381}%
\special{pa 3088 366}%
\special{pa 3119 351}%
\special{pa 3120 350}%
\special{sp}%
%
\special{pn 13}%
\special{pa 2780 980}%
\special{pa 870 970}%
\special{dt 0.045}%
\special{pa 870 970}%
\special{pa 871 970}%
\special{dt 0.045}%
%
\special{pn 13}%
\special{pa 1560 2210}%
\special{pa 1560 2890}%
\special{dt 0.045}%
\special{pa 1560 2890}%
\special{pa 1560 2889}%
\special{dt 0.045}%
\put(7.9000,-9.0000){\makebox(0,0)[rt]{$p_{i_n}$}}%
\put(1.4000,-21.1000){\makebox(0,0)[lt]{$I_{-N}\ni p_{i_1}$}}%
\put(15.1000,-29.9000){\makebox(0,0)[lt]{$p_{i_1}$}}%
\put(19.6000,-21.1000){\makebox(0,0)[lt]{$\vG_N = \vG_{-N}^{\vee}$}}%
\end{picture}%
\end{center}
\caption{The indeterminacy point $p_{i_n}$ of $\si^N$ 
is a superattracting fixed point of $\si^{-N}$}
\label{fig:cubic5}
\end{figure}
These two lines are mapped onto the single point $p_{i_1}$ by 
$\si^N$ since $\si^N(L) = \{p_{i_1}\}$. 
This implies that $p_{i_1}$ is a superattracting 
fixed point of $\si^N$, namely, 
\[
\nu(p_{i_1}, \si^N) = 
\det(I- (d \si^N)_{p_{i_1}}) = \det(I-O) = 1. 
\]
Likewise $p_{i_n}$ is a superattracting fixed 
point of $\si^{-N} = (\si^{-1})^N$ where $\si^{-1} = 
\si_{i_n} \cdots \si_{i_1}$ is the reduced expression 
of $\si^{-1}$ (see Figure \ref{fig:cubic5}), so that the same 
reasoning as above with $\si$ replaced by $\si^{-1}$ 
yields $\nu(p_{i_n}, \si^{-N}) = 1$. 
Therefore we have  
\[
\mu((p_{i_n},p_{i_n}), \vG_N \cap \vD) = 
\mu((p_{i_n},p_{i_n}), \vG_{-N}^{\vee} \cap \vD ) = 
\mu((p_{i_n},p_{i_n}), \vG_{-N} \cap \vD) = 
\nu(p_{i_n}, \si^{-N}) = 1. 
\]
These arguments imply that (\ref{eqn:lfpf3}) is equivalent 
to the assertion of the lemma. 
\hfill $\Box$ \par\medskip 
Putting Lemmas \ref{lem:lfpf1} and \ref{lem:lfpf2} together, 
we have established the following theorem. 
\begin{theorem} \label{thm:periodic} 
Let $\si \in G$ be any non-elementary AS map with length 
$n = \ell_G(\si)$. 
For any $N \in \N$ the cardinalities of periodic points 
of period $N$ are finite and explicitly given by 
\begin{equation} \label{eqn:carperiod}
\begin{array}{rcl}
\mathrm{\#} \, \ol{\Per}_N(\si;\th) &=& 
\l_1(\si)^N + (-1)^{nN} \l_1(\si)^{-N} + 4 (-1)^{nN} + 1, \\[2mm]
\mathrm{\#} \, \Per_N(\si;\th) &=& 
\l_1(\si)^N + (-1)^{nN} \l_1(\si)^{-N} + 4 (-1)^{nN}. 
\end{array} 
\end{equation}
The numbers grow exponentially as the period $N$ tends to 
infinity, with the growth rate $\l_1(\si)$. 
\end{theorem}
\section{Back to Painlev\'e VI} 
\label{sec:back} 
Back to the space of initial conditions for $\PVI(\k)$ 
through the Riemann-Hilbert correspondence, 
we are now able to deduce the dynamical properties of 
the Poincar\'e return map for $\PVI(\k)$ from the already 
established properties of the dynamical system on the 
affine cubic surface $\Sol(\th)$. 
This deduction is based on the following lemma. 
\begin{lemma} \label{lem:conjugacy} 
Assume that $\k \in \K-\Wall$. 
Given any loop $\ga \in \pi_1(Z,z)$, let $\si \in G(2)$ be the 
corresponding element via the isomorphism $(\ref{eqn:isom})$. 
Then the Poincar\'e return map $\ga_* : M_z(\k) \carl$ along 
the loop $\ga$ is strictly conjugated to the biregular map 
$\si : \Sol(\th) \carl$ via the Riemann-Hilbert correspondence 
$(\ref{eqn:RHtk2})$ and the commutative diagram 
$(\ref{cd:reduction})$. 
\end{lemma}
{\it Proof}. 
By Theorem \ref{thm:RH} the Riemann-Hilbert correspondence 
(\ref{eqn:RHtk2}) is biholomorphic under the assumption that 
$\k \in \K-\Wall$. 
Hence, for $i = 1,2,3$, the half-Poincar\'e map 
$\b_{i*}$ in (\ref{eqn:PSi}) is strictly conjugate to the 
transformation $g_i$ in (\ref{eqn:gi}). 
Being squared, $\b_{i*}^2$ is strictly conjugate to $g_i^2$. 
On the other hand, using formulas (\ref{eqn:gi}) and 
(\ref{eqn:si}), one can easily check that 
\begin{equation} \label{eqn:gi2} 
g_i^2 = \si_i \si_{i+1}, 
\end{equation}
where the index should be considered modulo $3$. 
Furthermore, in view of formula (\ref{eqn:r*}), 
the Poincar\'e return map $\ga_{i*} : M_z(\k) \carl$ is strictly 
conjugate to $\beta_{i*}^2 : \M_t(\k) \carl$ via the commutative 
diagram (\ref{cd:reduction}). 
Then the lemma is established by combining all these observations. 
\hfill $\Box$ \par\medskip 
The above conjugacy principle stands on the isomorphism of 
groups $\pi_1(Z,z) \to G(2)$ in (\ref{eqn:isom}), where 
the abstract group $G(2)$ in \S\ref{sec:result} is identified 
with its concrete realization as a group of birational maps 
on $\ol{\Sol}(\th)$ (see Remark \ref{rem:coxeter}). 
In order to utilize the results on cubic surface, 
we need to establish certain relations between the above two 
groups, e.g., between the minimality of a loop in $\pi_1(Z,z)$ 
and the analytic stability of an element in $G(2)$, etc. 
\begin{lemma} \label{lem:mini-as} 
Let $\si \in G(2)$ be the image of a loop $\ga \in \pi_1(Z,z)$ 
under the isomorphism $(\ref{eqn:isom})$. 
\begin{enumerate} 
\item If $\ga$ is minimal in the sense of 
Definition $\ref{def:minimal}$, then $\si$ is AS. 
\item If moreover $\ga$ is non-elementary in the sense of 
Definition $\ref{def:elementary}$, then $\si$ is 
non-elementary in the sense of Definition $\ref{def:AS}$.
\end{enumerate}
\end{lemma} 
{\it Proof}. Let (\ref{eqn:reduced}) and (\ref{eqn:reduced2}) 
be the reduced expressions of $\ga$ and $\si$ respectively, 
where $n = 2m$. 
Assume the contrary that $\si$ is not AS, namely, that 
$i_1 = i_n$. 
The argument is separated into two cases: 
Case 1 where $i_2 = i_{n-1}$ and Case 2 where 
$i_2 \neq i_{n-1}$. 
If we define $\si'$ and $\tau$ by 
\[
\si' := \left\{
\begin{array}{l}
\si_{i_3} \si_{i_4} \cdots \si_{i_{n-2}} \\[2mm]
\si_{i_{n-1}} \si_{i_2} \cdots \si_{i_{n-2}}
\end{array} \right. 
\qquad 
\tau := \left\{
\begin{array}{ll}
\si_{i_1} \si_{i_2} \qquad & 
(\mbox{Case 1}), \\[2mm]
\si_{i_1} \si_{i_{n-1}} \qquad & 
(\mbox{Case 2}). 
\end{array} \right. 
\]
then one has $\si = \tau \si' \tau^{-1}$ and the length 
of $\si'$ is given by 
\[
\ell_G(\si') = \left\{
\begin{array}{ll}
n-4 = 2(m-2) \qquad & (\mbox{Case 1}), \\[2mm]
n-2 = 2(m-1) \qquad & (\mbox{Case 2}). 
\end{array} \right. 
\]
Let $\ga'$, $\delta \in \pi_1(Z,z)$ be the loops corresponding 
to $\si'$, $\tau \in G(2)$. 
Then $\ga = \delta \ga' \delta^{-1}$ and 
\[
\ell_{\pi_1}(\ga') = \left\{
\begin{array}{ll}
m-2 \qquad & (\mbox{Case 1}), \\[2mm]
m-1 \qquad & (\mbox{Case 2}). 
\end{array} \right. 
\] 
In either case $\ga'$ is conjugate to $\ga$ and 
the length of $\ga'$ is smaller than that of $\ga$. 
This contradicts the minimality of $\ga$ and hence 
$\si$ must be AS, which proves assertion (1). 
Assertion (2) easily follows from Definitions 
\ref{def:elementary} and \ref{def:AS} and 
the translation rule (\ref{eqn:transl}). 
\hfill $\Box$ \par\medskip 
We are now in a position to establish our main results,  
Theorems \ref{thm:chaos}, \ref{thm:per} and 
\ref{thm:algorithm}, together with the related statements 
in Remarks \ref{rem:mostelem} and \ref{rem:lyapunov}. 
Let $\ga \in \pi_1(Z,z)$ be any non-elementary loop and  
$\si \in G(2)$ be the corresponding element under the 
isomorphism (\ref{eqn:isom}). 
As mentioned in Remark \ref{rem:conjugacy}, we may assume 
without loss of generality that $\ga$ is minimal. 
By Lemma \ref{lem:mini-as} the birational map 
$\si : \ol{\Sol}(\th) \carl$ is AS and non-elementary, 
so that Theorems \ref{thm:chaos2}, \ref{thm:entropy} and 
\ref{thm:periodic} can be applied to the map $\si$. 
Then the concluding arguments of this article proceed 
as follows. 
\par\medskip\noindent
{\bf Proof of Theorem \ref{thm:chaos} and (1) of 
Theorem \ref{thm:algorithm}}. 
Let $\mu_{\si}$ be the $\si$-invariant Borel probability 
measure stated in Theorem \ref{thm:chaos2}. 
As is mentioned in Remark \ref{rem:affine}, the measure 
$\mu_{\si}$ can be restricted to the affine cubic surface 
$\Sol(\th)$ without losing any mass and any ergodic 
properties. 
The resulting measure on $\Sol(\th)$ is also denoted by 
$\mu_{\si}$. 
We pull it back to the space $M_z(\k)$ of initial conditions 
via the Riemann-Hilbert correspondence. 
Let $\mu_{\ga}$ be the resulting measure on $M_z(\k)$. 
It is now clear from Theorems \ref{thm:chaos2} and 
\ref{thm:entropy} that the measure $\mu_{\ga}$ satisfies 
all the requirements in Theorem \ref{thm:chaos} 
and in assertion (1) of Theorem \ref{thm:algorithm}. 
Here note that formula (\ref{eqn:fdd}) leads to 
(\ref{eqn:lambda}), since the length $n = \ell_G(\si)$ 
of $\si \in G(2)$ is an even integer. 
\hfill $\Box$ \par\medskip\noindent 
{\bf Proof of Theorem \ref{thm:per} and 
(2) of Theorem \ref{thm:algorithm}}. 
We have defined in (\ref{eqn:per}) the set $\Per_N(\ga; \k)$ 
of periodic points of period $N$ for the Poincar\'e return 
map $\ga_*$. 
By Lemma \ref{lem:conjugacy} the Riemann-Hilbert 
correspondence (\ref{eqn:RHtk2}) maps $\Per_N(\ga;\k)$ 
bijectively onto $\Per_N(\si;\th)$ and hence 
\[
\mathrm{\#}\,\Per_N(\ga;\k) = \mathrm{\#}\,\Per_N(\si;\th). 
\]
Then Theorem \ref{thm:per} and assertion (2) of Theorem 
\ref{thm:algorithm} are an immediate consequence of the 
above equality and Theorem \ref{thm:periodic}, where 
we note that $n = \ell_G(\si)$ is even in the formula 
(\ref{eqn:carperiod}). \hfill $\Box$ 
\begin{remark} \label{rem:proof} 
Detailed explanations of Remarks \ref{rem:mostelem} 
and \ref{rem:lyapunov} are in order at this stage. 
\begin{enumerate} 
\item The first half of Remark \ref{rem:mostelem} follows 
from Lemma \ref{lem:fdd}. 
Indeed one has $\l_1(\si) \ge 3 + 2 \sqrt{2}$ for every 
even non-elementary map $\si \in G(2)$. 
Here one has the equality if and only if 
$\si = \si_i\si_j\si_k\si_j$ or 
$\si = \si_j \si_i \si_j \si_k$ for some 
$\{i,j,k\} = \{1,2,3\}$. 
As is easily seen, this occurs precisely when $\si$ 
comes from an eight-loop in Example \ref{ex:per} 
through the isomorphism (\ref{eqn:isom}). 
For example, if $(i,j,k) = (1,2,3)$ then 
$\si = \si_1\si_2\si_3\si_2$ comes from 
the eight-loop $\ga_1\ga_2^{-1}$. 
\item An inspection of formula (\ref{eqn:gi}) shows 
that the transformation $g_i^2 = \si_i \si_{i+1}$ in 
(\ref{eqn:gi2}) preserves the fibration 
$\Sol(\th) \to \C$, $x = (x_1,x_2,x_3) \mapsto x_k$, 
where $(i,j,k)$ is the cyclic permutation of $(1,2,3)$. 
Pull it back to the space $M_z(\k)$ via the 
Riemann-Hilbert correspondence. 
Then the resulting fibration 
$M_z(\k) \to \C$ is preserved by the 
Poincar\'e return map $\ga_{i*} : M_z(\k) \carl$ along 
the $i$-th basic loop $\ga_i$. 
This explains the second half of Remark \ref{rem:mostelem}. 
\item As is mentioned in Remark \ref{rem:area2}, 
the Riemann-Hilbert correspondence is an area-preserving 
biholomorphism between $(M_z(\k), \omega(\k))$ and 
$(\Sol(\th), \omega(\th))$ which intertwines 
the invariant measures $\mu_{\ga}$ and $\mu_{\si}$. 
Therefore Remark \ref{rem:lyapunov} readily follows 
from Remark \ref{rem:lyapunov2}. 
\item In connection with Example \ref{ex:per} we give the 
relation between Pochhammer loops in $\pi_1(Z,z)$ and 
Coxeter elements in $G$. 
For any cyclic permutation $(i,j,k)$ of $(1,2,3)$, 
the Pochhammer loop $\wp =[\ga_i,\ga_j^{-1}]$ 
corresponds to the square $c^2$ of the Coxeter element 
$c = \si_i\si_j\si_k$ via the isomorphism (\ref{eqn:isom}). 
Hence one has 
$\l(\wp) = \l_1(c^2) = \l_1(c)^2 = (2+\sqrt{5})^2 = 
9+4 \sqrt{5}$, 
which yields the formula in (2) of Example \ref{ex:per}. 
\end{enumerate} 
\end{remark}
\par 
In this article we have observed that the geometry of cubic 
surfaces and dynamical systems on them play important 
parts in understanding an aspect of the global structure 
of the sixth Painlev\'e equation. 
Their relevance to other aspects will be discussed 
elsewhere (e.g.~\cite{Iwasaki4}). 
\par\medskip\noindent
{\bf Acknowledgment}. 
The authors are grateful to Yutaka Ishii and Yasuhiko Yamada for 
their comments on an earlier version \cite{IU} of this article, 
which were quite helpful in revising the manuscript up to 
the present version. 

\end{document}